\documentclass{article}
\usepackage{arxiv}
\usepackage[utf8]{inputenc} 
\usepackage[T1]{fontenc}    
\usepackage{hyperref}       
\usepackage{url}            
\usepackage{booktabs}       
\usepackage{amsfonts}       
\usepackage{nicefrac}       
\usepackage{microtype}      
\usepackage{lipsum}		
\usepackage{graphicx}
\usepackage{natbib}
\usepackage{doi}
\usepackage{tabularx}
\usepackage{multicol}
\usepackage{multirow}
\usepackage{pdflscape}
\usepackage{array}
\usepackage{tikz}
\usepackage{tabularx}
\usepackage{multicol}
\usepackage{multirow}
\usepackage{pdflscape}
\usepackage{array}
\usepackage{threeparttable}
\usepackage{tikz}
\usetikzlibrary{arrows, shapes, positioning, calc, shadows}
\usepackage{ragged2e}
\usepackage{amsmath}

\usepackage{microtype}

\usepackage[para]{footmisc}
\usepackage{setspace}
\usepackage[format=hang,singlelinecheck=false]{caption}
\usepackage{ragged2e}

\usepackage{longtable}

\let\svthefootnote\thefootnote
\newcommand\blankfootnote[1]{%
    \let\thefootnote\relax\footnotetext{#1}%
    \let\thefootnote\svthefootnote%
}

\title{A Review on Intermodal Transportation and Decarbonization: An Operations Research Perspective}

\author{ \href{https://orcid.org/0000-0001-7334-9500}{\includegraphics[scale=0.06]{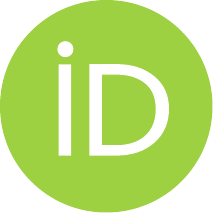}\hspace{1mm}Madelaine ~Martinez Ferguson}\\
	Department of Industrial and Systems Engineering\\
	University of Tennessee\\
	Knoxville, TN 37996 \\
	\texttt{mmart199@vols.utk.edu} \\
    \And
	\href{https://orcid.org/0000-0002-2285-845X}{\includegraphics[scale=0.06]{orcid.pdf}\hspace{1mm}Aliza ~Sharmin}\\
	Department of Industrial and Systems Engineering\\
	University of Tennessee\\
	Knoxville, TN 37996 \\
	\texttt{asharmin@vols.utk.edu} \\
    \And
	\href{https://orcid.org/0000-0001-7465-7783}{\includegraphics[scale=0.06]{orcid.pdf}\hspace{1mm}Mustafa C.~Camur} \\
	Department of Industrial and Systems Engineering\\
	University of Tennessee\\
	Knoxville, TN 37996 \\
	\texttt{mcamur@utk.edu} \\
	\And
	\href{https://orcid.org/0000-0003-1990-0159}{\includegraphics[scale=0.06]{orcid.pdf}\hspace{1mm}Xueping ~Li}\\
	Department of Industrial and Systems Engineering\\
	University of Tennessee\\
	Knoxville, TN 37996 \\
	\texttt{Xueping.Li@utk.edu} \\
}

\renewcommand{\shorttitle}{A Review on Intermodal Transportation and Decarbonization}

\hypersetup{
pdftitle={A template for the arxiv style},
pdfsubject={q-bio.NC, q-bio.QM},
pdfauthor={David S.~Hippocampus, Elias D.~Striatum},
pdfkeywords={First keyword, Second keyword, More},
}

\begin{document}
\maketitle



\begin{abstract}
This paper reviews intermodal transportation systems and their role in decarbonizing freight networks from an operations research perspective, analyzing over a decade of studies (2010-2024). We present a chronological analysis of the literature,
illustrating how the field evolved over time while highlighting the emergence of new research avenues.  We observe a significant increase in research addressing decarbonization since 2018, driven by regulatory pressures and technological advancements. Our integrated analysis is organized around three themes: a) modality, b) sustainability, and c) solution techniques. Key recommendations include the development of multistage stochastic models to better manage uncertainties and disruptions within intermodal transportation systems. Further research could leverage innovative technologies like machine learning and blockchain to improve decision-making and resource use through stakeholder collaboration. Life cycle assessment models are also suggested to better understand emissions across transportation stages and support the transition to alternative energy sources.
 
\end{abstract}


\section{Introduction}


Freight transportation (FT) is a vital component of the global supply chain (SC), facilitating commercial relationships among stakeholders (e.g., suppliers, distributors) from different regions. However, its broad application to improve efficiency along the SC leads to environmental consequences, including increased greenhouse gas (GHG) emissions \citep{zhang2022acyclic}. The global freight transportation industry is responsible for up to 11\% of the world's GHG emissions \citep{international2018co2}.



The pursuit of environmentally conscious solutions prompted a paradigm shift in FT towards diverse transportation configurations. Intermodal transportation (IMT),  multimodal transportation (MMT), and synchromodal transportation (SMT) are examples of this endeavor, which involve the integration of multiple modes (e.g., road-rail, rail-waterways) of transportation \citep{tavasszy2017intermodality}.
In short, IMT involves multiple carriers and contracts for different segments; MMT uses a single contract with one provider managing multiple modes; and SMT dynamically adjusts transport modes in real-time based on current conditions. For example, in \cite{craig2013estimating}, IMT is found to have a 46\% lower average carbon intensity than road transport, resulting in significant environmental benefits. Our study primarily focuses on IMT, incorporating relevant aspects of MMT and SMT where applicable.

Governments worldwide are expanding programs and funding opportunities to minimize carbon emissions during IMT operations. For instance, the European Union (EU) set modal shift goals by 2030 in its ``White Paper on Transport" \citep{ec2011roadmap}. Similarly, in the United States (US),  the Advanced Research Projects Agency-Energy (ARPA-E)  has started initiatives to enhance the IMT  sector by promoting low-carbon emission modes through technology and data-driven solutions \citep{ARPAE2023}.

Researchers are integrating sustainability dimensions and decarbonization strategies into the study of IMT to enhance FT efforts in minimizing GHG impacts and align the sector with changing regulatory requirements \citep{ghisolfi2024dynamics}. Over the past decade, these challenges and developments have sparked interest in the Operations Research (OR) community. There has been extensive research, including survey and review papers on these topics within the OR and transportation communities.

In early studies,  \citet{macharis2004opportunities} and \citet{merrina2007intermodal} condense
 studies with OR modeling techniques during the emergence of the concept of IMT, emphasizing the diverse actors involved and system complexity due to differing objectives.   However, these studies lack consideration of sustainability and carbon emissions. \cite{crainic2018simulation} propose a taxonomy to group studies utilizing simulation modeling between 2007 and 2017 in IMT, while \citet{Archetti20221ORMultimodal} focus on optimization models for long-haul transportation in IMT, without specifically addressing sustainability and decarbonization or grouping papers by various OR methodologies employed.  Therefore, integrating both aspects would constitute a significant contribution and address a crucial gap in the literature.

Several literature reviews in the broader context of decarbonizing FT   synthesized specific approaches and created analytical frameworks. \cite{emodi2022transportdeca} examine key drivers for decarbonizing FT in the Global South from 2010 to 2021 qualitatively.  Similarly, \cite{Aminzadegan2022decarintermodal} adopt a qualitative approach to assess environmental impacts and challenges in IMT for FT decarbonization. While these reviews contribute valuable insights, there is a need for a comprehensive synthesis that includes a broader range of OR techniques applied to IMT for decarbonization during freight operations.

Our research aims to enhance understanding of how OR methods contribute to sustainability in IMT networks. We systematically review and synthesize studies from the last 15 years (2010-2024) that apply OR techniques with a particular focus on advancing decarbonization efforts. Additionally, we reference foundational research from as far back as 2008 to provide context for our analysis. Our analysis categorizes these studies according to modality mix, decision level, emission consideration, and OR techniques used (see Section \ref{Integrated Analysis of Literature}), assessing their impact on sustainable IMT operations.  We also propose an agenda for future research.  Lastly, we note a significant increase in research activity and publications since 2018 (see Fig.~\ref{fig: CO$_2 $}), which highlights the timeliness and relevance of our work.


The remainder of our paper is organized as follows. We outline our review methodology in \ref{Review_Methodology}, provide an in-depth analysis of the studies in Sections \ref{Analysis} and \ref{Integrated Analysis of Literature}, discuss research gaps and future directions in Section \ref{FutureResearch}, and summarize our findings in Section \ref{Conclusions}.

\section{Chronological Analysis of Literature} \label{Analysis}
Since the early 21st century, efforts to evaluate transportation's environmental impact have focused on transitioning to lower-carbon modes like rail \citep{forkenbrock2001comparison}. Initially, these evaluations centered on unimodal systems, but the integration of intermodal transportation (IMT) networks into analyses has spurred the development and application of OR techniques. A timeline of emerging topics is presented in \ref{Timeline}.

\subsection{Early Interventions} 

\noindent \textbf{Early Studies on Cost and Emissions Trade-offs}: A foundational study in decarbonizing intermodal FT using OR  by \cite{kim2009trade}  investigates the  cost-CO$_2$ emissions trade-offs in intermodal and truck-only freight networks. Their findings reveal a nearly linear relationship: higher freight costs correlate with significant reductions in CO$_2$ emissions, demonstrated through Pareto optimal solutions.

Within the same year,  \citeauthor{kim2009assessment} further extend their analysis on CO$_2$ emissions for truck only and intermodal freight transport systems through two additional studies: \citep{kim2009assessment} and \citep{kim2009assessmentprevious}. 
The authors adopt a semi-Life Cycle Assessment (LCA) approach by focusing solely on exhaust and production emissions, excluding construction and vehicle manufacturing emissions. We refer the reader to \cite{horvath2006environmental} for further details regarding LCA in IMT.

\cite{kim2009assessmentprevious} first mathematically model emissions (i.e., direct and production emissions) across three freight modes (i.e., truck, diesel,  and electric locomotive) and associated processes (drayage, long-haulage, and transshipment), asserting that intermodal options typically emits less CO$_2$. Building on this, their subsequent study incorporates geographical data, diverse loading units, variable speeds, and vessel-based short sea shipping, validating their prior findings \citep{kim2009assessment}.

\noindent \textbf{Introduction of Geospatial Models:} Prior to \cite{kim2009trade}, \cite{winebrake2008assessing} propose the Geospatial Intermodal FT (GIFT) model, which finds optimal routes for a given origin-destination (O/D) pair. However, we believe that the study focuses on a generic network flow problem using the Dijkstra algorithm, utilizing connection arcs for modal shift and assuming equal costs and times for intermodal transfers,   potentially overlooking realistic constraints. Thus, we exclude their study from our analysis. Nevertheless, \cite{comer2010marine} improve the GIFT model by using a “hub-and-spoke” intermodal connection feature to investigate shifts from heavy-duty trucks to marine vessels to improve environmental efficiency in FT. 

\noindent \textbf{Incorporating Cost of Emission in Planning:} Traditional intermodal planning prioritizes minimizing travel time or costs until \cite{bauer2010minimizing} introduce GHG emissions, including CO$_2$, into the cost considerations of intermodal FT planning. GHG and carbon emissions differ; while carbon emissions specifically refer to carbon compounds like CO$_2$, GHG emissions include a broader range of gases such as CO$_2$, CH$_4$, N$_2$O, and fluorinated gases that contribute to the greenhouse effect \citep{epa_ghg}.

Considering that CO$_2$ accounts for approximately 80\% of GHG emissions \citep{liao2009comparing}, \cite{bauer2010minimizing} assess the impact of CO$_2$ emissions on optimizing IMT for fleet management (i.e., service schedules). Concurrently, \cite{chang2010optimization} consider the inclusion of mode-specific environmental costs in IMT planning with short sea shipping (SSS) to reduce congestion at major ports and promote eco-friendly transportation policies. In the following year, \cite{zhang2011mode} optimized the choice of transportation mode in IMT by accounting for the cost of carbon emission generated by transportation, mode transfer, and inventory. Later, \cite{le2013model} explore the effect of transportation mode selections (truck, sea, air) with vehicle-specific emission costs for a global SC.

\noindent  {\textbf{Dry Ports in Reducing Emissions:}} 
Dry ports, linking seaports to inland intermodal terminals via rail, offer potential CO$_2$ reductions through modal shifts creating a port-hinterland system.  A simulation model developed by \cite{roso2007evaluation} is one of the first to demonstrate a 25\% reduction in emissions using dry ports; however, the model's specifications are not thoroughly outlined in the study.  Building on this, \cite{henttu2011financial}  explore optimizing the network design (ND) of seaport-dry port systems, examining how the number and location of dry ports affect transportation costs using macro gravitational distribution models that consider both internal costs and external factors like noise, congestion, and CO$_2$ emissions as a mean to incorporate environmental and societal aspects into the analysis of port-hinterland systems. This gravitational model approach is then utilized by \cite{lattila2013hinterland} to assess dry port feasibility through discrete-event simulation (DES), comparing scenarios involving unimodal (road) networks and seaport-dry port networks. Around the same time,  \cite{iannone2012private} build upon the multi-commodity capacitated dry port model introduced in \cite{thore2010economic} and further expanded by \cite{iannone2012model} (latter two excluded from further analysis since emissions are not considered) by adding environmental cost (e.g., emissions,  accidents, congestion) to the generalized costs (e.g., in transit holding cost, container leasing cost) in the objective function.


\noindent \textbf{Carbon Pricing Mechanisms:} CO$_2$ emission pricing mechanisms, such as carbon taxes and cap-and-trade systems, influence transportation networks \citep{goulder2013carbon}.  Carbon taxes directly assign a price to CO$_2$ emissions through a predetermined tax rate, whereas cap-and-trade systems establish a total emissions limit and allow the market to set the price via tradable allowances, ensuring a consistent price across all participants. First, \cite{chaabane2011designing} and \cite{chaabane2012design} improve the mixed-integer linear programming (MILP) SCND model first introduced by \cite{ramudhin2008carbon}, incorporating strategies to minimize environmental impacts and integrate LCA for GHG estimation.

Continuing this exploration, \cite{pishvaee2012credibility} 
pioneer the use of stochastic programming (SP) in SCND to handle uncertainty by managing variables such as demand and facility capacity. Concurrently, \cite{holmgren2012tapas} introduce the first agent-based simulation model (ABM) to evaluate the repercussions of CO$_2$ taxes, infrastructure investments, and policy shifts on modal choices within freight networks. Subsequent studies delve into the broader implications of carbon pricing from a SC design perspective. \cite{fahimnia2013impact} evaluate the effects of carbon emission pricing on the planning and operations of closed-loop SCs, making it one of the first studies to analyze its impact on both forward and reverse SCs' carbon footprints. Similarly, \cite{zhang2013optimization} examine how emission pricing affects the optimal configuration of terminal networks by considering terminal capacity, location, and economies of scale, with a focus on long-term infrastructural supply modeling.

\cite{hoen2014switching} investigate how carbon emission regulations affect transportation mode selection, focusing on the cost-effectiveness of emissions reductions under variable demand and carbon constraints. This is expanded by \cite{rezaee2017green}, who develop a two-stage SP model that differentiates strategic decisions, such as location allocation, from tactical considerations like material flow, amid the uncertainties introduced by carbon pricing and fluctuating demand. Additionally, \cite{fahimnia2015tactical} and \cite{zhang2017multimodal} assess organizational impacts. \cite{fahimnia2015tactical} explore how carbon tax policies affect SC financial outcomes and emissions strategies tactically, while \cite{zhang2017multimodal} evaluate the influence of different CO$_2$ taxes on transportation mode selection and infrastructure investments, with a specific focus on customer service windows.


\noindent \textbf{Incorporating Societal Impacts:} Decarbonizing IMT involves addressing a range of societal impacts beyond just environmental concerns, such as noise pollution and accident risks. In this context in 2012, \cite{sawadogo2012reducing}  first integrate societal impacts (e.g., noise, accident risk) into IMT freight routing. Their approach provides routes that minimize combined environmental, social, and economic impacts, which offer policymakers flexible options. Subsequently, \cite{kim2013multimodal} optimize facility investment decisions in multimodal container freight networks by minimizing social costs related to GHG emissions, while adhering to maximum allowable emission limits. Although earlier research \citep{ricci2005social, macharis2010decision} considers societal impacts, it does not employ any OR methods.

\noindent \textbf{Balancing  Priorities and Cost-Emission Trade-Offs:}  In  IMT and MMT  management, decision-making involves balancing different stakeholder priorities: managers might prioritize reducing carbon emissions, carriers mainly focus on economic efficiency and speed, and shippers seek quick deliveries at reasonable costs.  Customizing priorities is crucial for addressing these multifaceted challenges effectively. 

To this end,  \cite{kengpol2014development} propose a decision support system (DSS) that minimizes transportation costs,  time, risk, and CO$_2$ emissions. 
They analyze transportation modes: air, truck, rail, and sea, both individually (unimodal) and in combination (multimodal), revealing that the most cost-effective route is the unimodal sea route. Conversely, when considering all other criteria combined, the unimodal air route emerges as the most favorable option. Yet, these findings contradict the common belief that MMT generally outperforms unimodal options across all factors. In addition, air transportation typically faces limitations due to high costs, weight and volume restrictions, and handling challenges. 

Continuing in this vein, \cite{soysal2014modelling} examine the trade-off between transport costs and CO$_2$ emissions in food SCs, illustrating a Pareto frontier similar to \cite{kim2009trade}. Their study, while focused on truck-based logistics, extends the analysis by incorporating factors such as road design, vehicle and fuel types, vehicle loads, product perishability, and backhauls. Although the research includes various factors and modes, it does not claim to cover multimodal aspects. \cite{pan2017multimodal}  then integrate demand allocation, transportation mode selection, and routing optimization to examine the trade-offs between costs and CO$_2$  emissions, finding that costs rise sharply beyond a nearly 50\% reduction in emissions.

It is worth mentioning that very few studies consider air transportation in the IMT literature. This is noteworthy since air cargo accounts for around 35\% of the total value of world trade \citep{air_35}. More importantly, aviation, which includes both passenger and cargo flights, contributes approximately 2.5\% of global CO$_2$ emissions and about 4\% of global warming to date \citep{owid-global-aviation-emissions}.

\noindent \textbf{Seaport Management and Routing Strategies:} In a series of studies, \cite{rodrigues2014assessing}, \cite{chen2014design} and \cite{rodrigues2015uk} explore strategies to optimize emission reductions in FT through strategic seaport expansion and selection. \cite{rodrigues2014assessing} study alternative seaport gateways and freight network restructuring to cut transport-related emissions.  \cite{chen2014design} propose a bi-level programming model where the upper level designs the IMT (road and waterway) service network factoring carbon emissions, state subsidies, and profit while the lower level manages traffic flow through a user-equilibrium assignment model. Further, \cite{rodrigues2015uk} assess the CO$_2$ impacts of re-routing containers, factoring in emission reductions from road use, container handling costs at seaports and rail terminals, emissions from seaport expansions, and traffic congestion effects.

\noindent \textbf{Integration of Sourcing and Production Decisions:} 
Decarbonizing IMT requires integrating sourcing and production decisions to address interconnected emissions and inefficiency factors. \cite{chaabane2011designing} extend the model from \cite{ramudhin2008carbon} by integrating Bill of Materials (BOM) constraints, enabling simultaneous optimization of sourcing, production, and transportation decisions. \cite{liotta2015optimisation} further explore this issue by assuming transportation emissions are linearly proportional to the number of items transported over an arc, directly linking emissions to transport volume.

However, it is crucial to note a discrepancy in the literature: the same team later proposed a multi-stage optimization and simulation approach to integrate aggregate production, MMT, and distribution plans, taking into account demand uncertainty and transport operators' coalitions \citep{liotta2014optimization}. These coalitions are based on multiple pre-negotiated agreements to share capacity and collectively operate transportation services across various modes. This approach does not align perfectly with a typical multimodal setting, where a single contract often governs operations. We also note a discrepancy in the publication year for this paper: Google Scholar lists 2014, while other studies cite 2016.

\subsection{Development of New Concepts}
\noindent \textbf{Government Policy Measures:} In 2015, we observe a shift toward including government intervention in research coincided with the adoption of the 2030 Agenda for Sustainable Development by the EU and all UN  states \citep{chasek2016getting}. This agenda prompted the study of tax implications for sustainability goals. Previously, \cite{maia2013strategic} proposed an optimization model to assist policy decisions in regional and national planning for intermodal road-rail networks, aiming to minimize both generalized costs and CO$_2$ emissions. However, their model did not explicitly address the government's perspective.

\cite{duan2015carbon}  pioneer research into how IMT operations change under government tax policies. Similarly, \cite{wang2015effects} examine tax implications using a two-stage Stackelberg gaming model including the government and a company–third-party logistics alliance as players, marking the first study to apply game theory (GT) in this context. We humbly believe that Wang's approach may be more realistic since the government does not fully control logistics operations.   

\cite{santos2015impact} also explore IMT through government policies but omit sustainability, leading us to exclude their study from our analysis. In the same year, 
\cite{bouchery2015cost} incorporate government incentives such as train subsidies and road taxes in their gravity model to analyze how transportation costs, emissions, and modal shifts in terminal locations impact IMT networks, considering both hub-and-spoke and direct road transit systems. From an MMT context, \cite{zhang2015freight} extend the analysis, formulating a bi-level optimization model 
that balances governmental decisions on operational and emission costs with service network decisions through cost-effective routing.

Further, \cite{lin2017modeling} use bi-level programming to improve railway infrastructure to facilitate the shift from road to rail, balancing investment, transportation, and emission costs by selecting government investment projects at the upper level and managing freight flow at the lower level. \cite{mostert2017road} examine the impact of road tax policy and external human health costs on the modal split between road, intermodal rail, and inland waterways. In a follow-up study, they introduce a model that incorporates economies of scale in IMT to assess the effects of policy measures (i.e., subsidies, internalization of external costs) \citep{mostert2018intermodal}. This research is significant for incorporating a door-to-door scheme into the traditional hub system, providing a more realistic view of decarbonization scenarios.

\cite{haddadsisakht2018closed} develop a hybrid robust-stochastic program considering demand and carbon tax uncertainty, finding that adjusting transportation capacities is more effective than altering product flows or building additional facilities in response to carbon tax uncertainty. Next, \cite{zhang2018joint} expand on the approach of \cite{zhang2015freight} and integrate infrastructure investments and subsidies in their bi-level programming to meet CO$_2$ reduction targets. In a similar direction, \cite{choi2019system} advance the field by introducing the first system dynamics (SD) model to assess how governmental policies (i.e., road taxes and containerization measures) impact modal shift rates, revealing that new policies have a more pronounced effect.

Moving forward to 2020, 
\cite{jiang2020regional} incorporate government and logistics carriers' decisions in a regional logistic network. Their model emphasizes logistic parks' decisions regarding the number, location, and capacity, integrating demand uncertainty and the determination of railway subsidies while considering a national carbon reduction target using a robust optimization approach. Then, while examining freight operations on railroads, \cite{li2020integrated} introduce consigners as a third stakeholder. Their model focuses on government carbon emission policies to encourage a modal shift from road-only to road-rail, while reducing subsidies. It integrates rail freight decisions (e.g., service frequency) with pricing strategies for operators and transportation mode choices for consigners, based on a utility function considering freight rates, carbon taxes, and other factors.

Later, \cite{yang2021coastal} take a step further and study an optimal policy-network joint design problem to determine the level of capacity expansion in a road-water IMT network. Concurrently, \cite{qian2021decision}  pursue an optimal modal split for regional freight, maximizing efficiency while balancing economic and environmental impacts, similar to \cite{bouchery2015cost}, while \cite{gallardo2021sequential} combine MMT planning with DES to guide ND and Service ND (SND) decisions and achieve a net-zero freight system. Expanding on governmental strategies to encourage modal shifts, \cite{pedinotti2022freight} analyze modal shifts from heavy trucks to trains using a TIMES-type energy model, assessing the effects of both endogenous variables like technology and modal choices, and exogenous variables such as price elasticities.
Recently, \cite{halim2023boosting} evaluate modal shifts across varying ambition levels using the ``avoid, shift, improve" framework. They highlight essential policies such as enhancing service quality and capacity and securing dedicated funding for IMT rail infrastructure to reduce emissions and costs.

\noindent  {\textbf{Optimizing Service ND:}} The importance of terminal operations in IMT has led to their increased inclusion in OR formulations since 2016.  \cite{qu2016sustainability} and \cite{rudi2016freight} incorporate previously neglected operational elements in multi-commodity SND problems. \cite{qu2016sustainability} incorporate terminal operational costs for intermodal unit handling into the objective function, whereas \cite{rudi2016freight} include the transshipment process resulting from carrier replacement at nodes in their objective function. In terms of emission estimations, \cite{qu2016sustainability} use an activity-based function to estimate transportation CO$_2$ emissions \citep{mckinnon2010measuring}, while \cite{rudi2016freight} account for GHG emissions from both transportation and transshipment, using emission factors and energy demand.

Subsequent studies leverage multi-objective optimization to include more comprehensive elements in SND analysis. 
\cite{baykasouglu2016multi} develop a
multi-period network flow model for sustainable load planning, minimizing overall transport cost and emissions while maximizing customer satisfaction. 
Their model considers simultaneous import and export flows, including consolidation operations and load exchanges at transshipment stations, under both crisp and fuzzy decision-making environments. Next, \cite{ji2017hybrid} address a multi-source single-destination flow problem, focusing on transshipment costs and time. They consider capacity constraints on intermediate terminals to mitigate traffic congestion and reduce transportation pressure. Further, \cite{zhou2018capacitated} refine MMT network flow planning by addressing realistic features such as costs associated with specific combinations of transportation modes and vehicle types. They include flow capacity per node, consolidation, rerouting of shipments, delivery time constraints, and emissions from both vehicle types and handling operations during transshipment. \cite{sun2019bi} develop a bi-objective model to manage the risk of population exposure to hazardous materials in a road-rail MMT system with a hub-and-spoke network configuration. Their approach aims to minimize total costs (e.g., transportation, handling, and storage costs), while also addressing exposure risk, with CO$_2$ emissions included as a constraint under uncertainty.

\noindent {\textbf{Incorporating Stochasticity in SND:} } Incorporating stochasticity into models enhances the realism of SND by addressing uncertainties such as network disruptions. \cite{demir2016green} pioneer this approach in SND by incorporating stochastic components such as service travel times and demand. Their weighted MILP formulation, including constraints to ensure service sequencing and manage plan reliability, is solved using sample average approximation (SAA). Excess demand is managed by adding additional services, which increases costs in the objective function, with emissions also accounted for as CO$_2$ costs.

Building on this foundation, \cite{hruvsovsky2018hybrid} focus on the uncertainty of service travel times, modifying the MILP model to include in-transit inventory costs and enhancing capacity management for larger instances through the first hybrid simulation-optimization model in the SND domain. \cite{layeb2018simulation} address the same uncertainties using a simulation optimization approach aiming for more accurate travel time data. Further, \cite{sun2018time} expand the scope by considering multiple uncertainties (i.e. railway service capacity and loading times) in a multi-objective fuzzy MINLP model for a rail-road network. This model, which utilizes both hub-and-spoke and point-to-point configurations, charges CO$_2$ emissions as a separate objective. We note that while hub-and-spoke configurations are not new in the SND field \citep{bouchery2015cost}, this study does not cite earlier works. Later, \citet{baykasouglu2019fuzzy} study several sources of uncertainties (e.g., cost, vehicle capacity, transit time) by developing a hybrid fuzzy and SP model for a fleet planning model in the IMT context.

\noindent{}{\textbf{Sustainability and Stakeholder Dynamics in Port-hinterland Systems:}} In the domain of ND and dry ports, \cite{tran2017container} emphasize the importance of integrating inland connections into the global design of container shipping networks. Their model focuses on reducing costs associated with maritime segments (i.e., ship and port expenses), hinterland connections (i.e., inland/feeder transportation and inventory costs), and carbon emissions across the entire container flow. Expanding on this concept, \cite{tsao2018seaport} develop a model using a GT approach for designing a multi-echelon seaport-dry port network. The model accounts for multiple stakeholder decisions, including storage pricing and port service areas, while shippers decide between an all-road network or a combined road-rail network connecting dry ports to seaports. However, despite referencing \cite{qiu2015bilevel} and stating that \textit{``In that study, carbon emissions were calculated for shared as well as independent transportation systems..''}, no evidence of carbon emissions assessment is found in \cite{qiu2015bilevel}.

Following a similar methodology in \cite{tsao2018seaport} but with a different focus, \cite{xu2018modelling} model port competition, where ports manage dry ports to maximize their payoff, and shippers make informed routing decisions. Their model considers emissions from various transportation modes and transfer operations at dry ports and seaports. Further, \cite{liu2019system} analyze the port-hinterland system focusing on emissions from FT, utilizing a system dynamics (SD) model that considers fuel consumption influenced by dynamic factors such as fleet configuration and cargo volume.

Other studies bring into aspects of dry ports with the potential to increase the efficiency and sustainability of freight flow.  \cite{tsao2019multi}  integrate key sustainability components (economy, society, and environment) along with uncertainty factors (demand, cost, capacity, CO$_2$ emissions, social cost) into dry port ND, combining robust optimization with fuzzy programming. Their model expands the social costs addressed earlier \citep{henttu2011financial, lattila2013hinterland} by adding unemployment, immigration, and traffic congestion. Decisions of both government (i.e., quantity, location, and capacity of dry ports), and shippers (i.e., routing, and quantity of freight) are considered while CO$_2$ emissions are calculated from dry port opening and operations and modal transportation.

Building on the global SC analyses similar to \cite{tran2017container}, \cite{wei2019import} optimize freight flows within an integrated cross-border logistics network, connecting supply points and destinations through a network of dry ports and seaports. Their bi-objective MIP model is developed assuming that the logistics integrator is responsible for meeting customer and government requirements and focused on calculating emissions per container and setting emission limits based on local government regulations.  \cite{li2019carbon} also examine freight flow by extending the study by \cite{rodrigues2015uk} on utilizing dry ports as a strategy to reduce carbon emissions. They assess MMT's cost and environmental impacts on seaport selection, aiming to minimize total costs, including transportation, handling, and emissions, where emission costs are derived from container handling and transportation.

In port-hinterland IMT, \cite{wang2020integrated} address port-centric SCs (see \cite{mangan2008port} for details) by integrating production and transportation scheduling for a made-to-order SC under uncertainty with multiple modes of transportation.  
\cite{dai2020distributionally} emphasize container transportation ND with waterways and intermodal terminals under port-hinterland demand uncertainty, contrasting with the road-rail focus of previous studies. Additionally, \cite{abu2020optimization} enhance IMT sustainability by optimizing container terminal layouts and hinterland transportation, considering both transportation and storage costs, as well as emissions. \cite{yin2021interrelations} optimize hinterland division in port-center MMT, using a multi-objective model to minimize transportation costs, emissions, and travel time.

\noindent{}{\textbf{Empty Container Repositioning :}} 
Within IMT SND, empty container management remains a pivotal area. Empty container repositioning problems (ECRP) involve decisions on the quantity, route, and location of empty containers \citep{choong2002empty}. 
Since 2002, the importance of ECRP for IMT efficiency has been recognized, but it was not until 2016 that \cite{lam2016market} included it in the context of IMT decarbonization. They consider waterway-road, rail-road, and road-only combinations for inland transportation, integrating ECRP to minimize transportation costs and limit emissions. Later,  \cite{zhao2018stochastic} employ a two-phase model to first generate a distribution and leasing plan for empty containers, followed by routing plans back to demand terminals, incorporating emissions and uncertainties related to container demand and supply. Further, \citet{castrellon2023assessing} utilize a combination of ABM and DES to evaluate ECR strategies at varying levels of container substitution using dry ports.

\noindent  {\textbf{Refining GHG Emission Estimation:}} Recently, various models have been developed to estimate GHG or CO$_2$ emissions. Emission estimation models are classified into energy-based or activity-based methods and further differentiated by their output precision into microscopic or macroscopic models \citep{demir2014review}. Addressing the gap between the detailed accuracy of microscopic models and the data efficiency of macroscopic models, \cite{kirschstein2015ghg} introduce mesoscopic models for road-rail IMT, which utilize parameters such as speed, weight, and traffic conditions to assess energy demand.

These mesoscopic models are further expanded by \cite{heinold2018emission} to obtain large-scale emission estimations, incorporating high shipment volumes and novel factors like vehicle and shipment characteristics, as well as route-specific variables such as gradient. On the activity-based front, \cite{pizzol2019deterministic} incorporates error propagation into LCA analysis to enhance carbon footprint estimations for sea-road transportation. In contrast, \cite{tao2021energy} improve estimates by including emissions from well-to-wheel (WTW) transshipment operations in hinterland transport, highlighting the importance of considering ECRP in these calculations. Furthermore, \cite{guo2022modeling} broaden this approach by taking into account emissions from electric modes across the hinterland network, acknowledging that earlier studies of port-hinterland systems \citep{wang2020integrated} did not fully account for certain complexities of the emission process (e.g., emissions at node level).

\subsection{Advancement of Concepts}

\noindent  {\textbf{Innovations in Food SCs for Emission Reduction:}}
In 2018, new strategies within food SCs aimed at reducing emissions were explored, particularly focusing on aspects like postponing packaging to minimize freight weight.  \cite{harris2018impact}  examine the impact of different components on SC emissions, using IMT for wine distribution. Their analysis models various scenarios involving combinations of routes, modes of transport (road and waterway), and packaging types, assessing both CO$_2$ and sulfate emissions, which are significant pollutants from the maritime sector.

Subsequently, \cite{ma2018optimization} address emission reduction from another perspective within the food industry by developing a novel scheduling optimization model from the shipper's viewpoint. This model aims to reduce costs, mitigate quality degradation of food products, and cut CO$_2$ emissions, with a particular focus on cold chains. Cold chains are especially significant as they consume considerable energy mainly due to refrigeration during transport, a factor often overlooked in studies focusing on non-perishable goods \citep{fitzgerald2011energy}. Additionally, \cite{maiyar2019modelling} propose a model to optimize hub locations and IMT decisions in food grain SCs. Their multi-period MINLP model addresses constraints related to hub storage capacity, handling capacity, vehicle capacity, and availability.

\noindent{\textbf{Advancements in Freight Planning Tools:}}
In 2019, significant progress was made in developing freight planning tools and gathering insights for stakeholders' modal shift decisions. These tools are categorized as online or offline, depending on whether they utilize historical or real-time data and are used before or during service.  For instance, \cite{demir2019green} focus on an offline planning tool (i.e., using deterministic data before transportation service) to provide timely environmental and economic solutions. Building upon their earlier work on the stochastic ND problem in green IMT \citep{demir2016green}, they eliminate constraints related to travel time, demand uncertainties, and penalties for delivery delays.

In 2020, \cite{laurent2020carbonroadmap} highlight the need for operational transportation planning tools in IMT. Developed in early 2017 but published three years later, this possibly pioneering work provides visual representations of carbon emissions for both unimodal (road only) and multimodal (combining road, rail, and waterways) transportation, introducing the concept of the intermodal carbon-efficient boundary.

\noindent{}{\textbf{Emerging Focus on SMT:}}
The concept of SMT has gained attention for its potential to enhance IMT benefits through dynamic planning and real-time information  \citep{giusti2019synchromodal}. \cite{zhang2016synchromodal} demonstrate that SMT could reduce total carbon emissions by 28\% compared to IMT, although this study was excluded from our analysis since it only considered CO$_2$ emissions as a metric and not as part of the optimization model. In 2018, \cite{farahani2018decision} develop a DSS to aid SC decision-makers in balancing energy consumption and costs in SC design and operation, incorporating air as a fourth mode alongside other shipment factors such as packaging type, quantity, and container sizes. 

Later, \cite{resat2019discrete} expand on this by including previously overlooked benefits of SMT (e.g., vehicle capacity and speed) using the Synchromodal ND and operational IMT model initially created by \cite{resat2015design} (omitted from our study since emissions are not included). Their significant contribution is a time-dependent synchromodal model that facilitates real-time planning in response to vehicle capacity availability or disruptions. In the following year, \cite{guo2020dynamic} develop an online system to match shipment requests with transportation services in SMT, aiming to minimize costs while managing service flexibility and capacity. Later, \cite{hruvsovsky2021real} focus on a DSS that combines offline planning with online re-planning under disruptions in SMT. This approach differs from their previous work \citep{demir2019green} by employing a hybrid simulation-optimization method using the MILP model to create and monitor transportation plans while addressing potential disruptions.

\noindent{}{\textbf{Enhancing Reliability in IMT and MMT:}}  In 2020, advancements in OR models focused on improving the reliability of IMT and MMT networks by addressing uncertainties. \cite{sun2020green} identify travel time as the most frequently studied source of uncertainty in IMT routing. They incorporate uncertainties related to network capacities, departure times, and loading times into a road-rail network using a multi-objective fuzzy MINLP  model. CO$_2$ emissions are included in the objective, calculated using an activity-based method discussed in \cite{sun2018time}. In contrast, \cite{wang2020transportation}  focus on node capacity uncertainty within a road-inland waterways network. They employ a fuzzy MILP model that includes emissions as a cost component in the objective function, subject to a regional carbon tax \citep{zhang2017multimodal}.

In the broader context of SC reliability,
\cite{kabadurmus2020sustainable} integrate supplier risk evaluation into a multi-product multimodal SC network, using carbon cap-and-trade for CO$_2$ mitigation. Their MILP model includes BOM as used in \cite{chaabane2011designing}, to account for emissions from both production and transportation. The model also features emission trading as part of the objective function, employing emission models from \cite{mckinnon2010measuring} for water and rail, and \cite{mckinnon2015green} for air and road. In a similar context, \citet{mousavi2020robust} study SCND decisions at strategic and tactical levels for bioethanol SC under operational disruption.


\noindent{}{\textbf{Emphasizing GHG Emissions in Route Planning:}} In 2020, three studies integrated GHG emission considerations into transport route planning. \cite{maneengam2020bi} criticizes traditional approaches for prioritizing cost or time over GHG emissions in MMT route planning. However, \cite{sun2015modeling} show that CO$_2$ emissions have been part of MMT routing considerations since 2015. In contrast, \cite{wang2020modelling} incorporate CO$_2$ emissions in OR studies to enhance routing, developing a multi-objective MMT framework combining a DES model with operational and environmental efficiency metrics. Additionally, \cite{heinold2020emission} introduce order-specific emission limits requested by shippers, using egalitarian and payload-based allocation schemes within a road-rail IMT network.



\noindent{}{\textbf{Various Logistic Service Provider (LSP) Perspectives:}} FT has been explored from different LSP perspectives to examine factors influencing decision-making and operational strategies. From a stakeholder perspective, \cite{fulzele2019model} assess the influence of various factors (e.g., costs, damage, delays, and GHG) on modal shift decisions conducted by logistics service providers. In a similar context, \cite{wang2021bi} explore competitive pricing strategies between logistic enterprises to balance FT modes through a bi-level multi-objective approach, optimizing revenue and emission reductions while maintaining market equilibrium.
\cite{zhang2021low} address dual uncertainty (demand and time) in  MMT  from an LSP's perspective, using a hybrid robust stochastic formulation to minimize transportation costs, emissions, and time.  
Emissions are considered under a cap-and-trade policy \citep{chaabane2011designing}. 

Similarly, \cite{tiwari2021freight} propose a model for third-party logistic providers (3PLs) that separates long and short transportation decisions under a carbon tax policy. The model addresses decisions on containerization and consolidation with emissions also accounted for through a cap-and-trade policy.  
\cite{wu2021research} go further by focusing on logistic ND in the bulk sector for fourth-party logistic providers (4PLs), creating a platform for efficient freight distribution across multiple supply and demand nodes, where emissions are primarily calculated from transportation activities. Adding to the work of \cite{zhang2021low}, \cite{li2022path} assess the effects of uncertain demand and random carbon trading prices on MMT routing design. They aim to reduce transportation, transfer, time, and emission costs by employing a hybrid fireworks algorithm with a gravitational search operator to solve the model.




\noindent{}{\textbf{Enhancing Travel Time Efficiency on SND:}}
Incorporating travel time into the formulation of SND problems influences not only the scheduling and routing of transportation modes but also customer satisfaction, resource allocation, and network performance. In this direction, \cite{lu2019sustainable} integrate hub location and vehicle routing decisions in MMT along with travel time consideration, analyzing a two-echelon network that includes highways and domestic and international railway systems. \cite{liang2021multi} optimize paths in the MMT context through a multi-objective model, incorporating storage costs due to early arrival and penalty costs for late arrival. Additionally, they include a third objective that considers logistics service satisfaction, measured through indexes that account for timeliness, economic factors, and service quality. Their model allows the selection of any of the three available modes at each node, though this is considered unrealistic in practical scenarios.

In the multi-commodity MMT context, \cite{qi2022transport} extend the SND analysis by including elements like departure dates and carbon emissions in shipment planning. 
Their MIP model focuses on door-to-door logistics, starting the service when an order is ready to depart from the origin warehouse rather than the shipper’s release date. The model analyzes outbound freight from China to Europe via the ``Belt and Road Initiative" \citep{huang2016understanding}.
Similarly, \cite{xie2022research} evaluate cross-border logistics, as seen in \cite{wei2019import}, but within the context of the initiative mentioned above. Their study integrates air transport into China's inland logistics network using a bi-objective MIP model, considering transportation costs (transport, port, and carbon taxes) and time (transportation and port storage).

Further, \cite{sun2022green} improve rail-road IMT network routing by incorporating time-varying parameters (e.g., travel and truck departure times), using fuzzy soft time windows to ensure timeliness. They also assess the impact of uncertain and time-varying carbon emission factors, calculated through an activity-based model \citep{liao2009comparing}. 
Their network strategy uses a hub-and-spoke system as the primary mode, supplemented by direct road transport, reflecting the methodology of \cite{sun2018time}.

\noindent {\textbf{Integrating Eco-Labels into IMT:}}
Building on \cite{heinold2020emission}, \cite{zhang2022collaborative}
develop a bi-objective model for collaborative planning among carriers to address shippers' sustainability demands.  This model employs an eco-label system \citep{galarraga2002use} preferred by shippers, where carriers maximize their fulfilled requests at the upper level while minimizing transportation costs at the lower level, extending methodologies from \cite{guo2020dynamic}. Furthermore,  \cite{heinold2023primal} explore eco-labels' influence in IMT utilizing a multiobjective sequential decision process, which reflects customers' environmental awareness and sets emission limits per order. They incorporate stochasticity in new-order arrivals and use reinforcement learning to address the decision-making process.

\subsection{Recent Developments}
\noindent {\textbf{Current advancements in SMT: }} Within 2022, several studies on synchromodality focus on optimizing schedules, user preferences, and resource utilization. First, \cite{zhang2022preference} examine carriers' preferences regarding cost, time, and emissions, creating a multi-objective model to optimize total costs (i.e., container, fuel, and carbon tax costs), emissions per mode and container, and overall time (travel and waiting times). In another study,
\cite{zhang2022heterogeneous} address shippers' diverse preferences using a planning MIP model coupled with a multiple criteria decision-making (MCDM) approach using fuzzy set theory. Further, \cite{zhang2022flexible} explore service flexibility to enhance resource use, presenting an MILP that includes fixed and flexible routes, and complex schedules. In the following year, \cite{oudani2023combined} 
 introduces the first blockchain technology to SMT, proposing a green blockchain framework coupled with energy-efficient management solutions. Their approach uses weighting and epsilon constraint methods to produce Pareto-optimal solutions assessed via multiple MCDM methods.

\noindent \textbf{Game Theoretic Approaches to  Sustainability:} Promoting low-carbon practices in  MMT  requires addressing the diverse interests of stakeholders. We observe a surge in the use of GT models, particularly in 2023.  \cite{shams2023game} employ a GT  approach to compare three carbon regulations—cap and trade, carbon offset, and carbon tax—considering various governmental perspectives, including social, economic, and environmental dimensions. Similarly, \cite{fallahi2023game} analyze competition in FT networks and sustainability dimensions, examining scenarios with and without government intervention. Unlike \cite{shams2023game}, \cite{fallahi2023game} focus on taxation based on fuel consumption rather than carbon emissions, although both studies treat the government as the Stackelberg leader and the FT systems as followers. However, since \cite{fallahi2023game} do not consider emissions, their study is excluded from our literature analysis figures and tables.

Additionally, \cite{wu2023game} consider carbon tax while optimizing the layout and location of dry ports by developing a GT model between the government and shippers.  This model accounts for capacity constraints, shipper behavior, and empty container transportation, distinguishing dry ports from other transit nodes.  Further, \cite{rahiminia2023adopting}  investigate the interactions between rail operators (leader) and shippers (follower) based on the pricing offered by the rail operator considering the triple bottom line in a sustainable transportation context. Lastly, \cite{chen2023subsidy} develop a subsidy decision model focused on emission reduction. They solve the Nash equilibrium in a tripartite pricing game using graphical methods, aiming to maximize carbon reduction per unit of subsidy while incorporating demurrage charges and estimating emissions reductions through an activity-based approach.

\noindent \textbf{Employing System Dynamics for Policy Making:} The analysis of port-hinterland IMT systems and the impact of carbon taxation on carbon emissions and economic growth is crucial for informed policymaking. Within this context, various studies have employed SD. For example, \cite{zhong2023system} utilize SD to simulate the dynamics of these transport systems and assess the effectiveness of carbon taxation in reducing CO$_2$ emissions while considering the costs of emission reduction. Their study emphasizes the need for local governments to tailor carbon tax policies to regional specifics.

Similarly, \cite{nassar2023system} employ SD to explore fiscal, regulatory, and infrastructure policies aimed at promoting shifts towards lower-emission modes in freight transport. Their findings indicate a gradual transition with diminishing effectiveness of modal shift measures over time, which slows decarbonization and reveals the limitations of relying solely on these policies for significant emission reductions. Furthermore, \cite{guo2023toward} combine SD with Monte Carlo simulation (MCS) to assess policy mixes, including cost-based pricing, road development, railway service improvements, and railway subsidies, while accounting for economic transition uncertainties aimed at promoting modal shifts. The study reveals that the collective impact of these policies is less than the sum of their individual effects, as they all target similar outcomes.


\noindent \textbf{Carbon Peak Times:}
Predicting carbon peak times is crucial for effective climate action and energy management. This peak marks the beginning of a decline in emissions, which is essential for sustainability and policy development. Research in this area typically aims to identify peak timings and develop strategies for managing emissions across various sectors. In this context, \cite{zuo2023using} develop an annual control model to adjust FT and manage energy consumption in China, predicting a peak in freight transport energy consumption by 2029. This aligns with the findings of \cite{yu2018achievement}, which highlight the sensitivity of carbon emission peaks to energy consumption intensity.

In recent years, \cite{guo2023carbon} simulate carbon peak emission times and values in container intermodal networks using MCS. Their case study suggests that meeting carbon peaking targets before 2030 will require careful management of energy technologies and bi-control measures, such as optimizing transport structure and network connectivity. Also, \cite{ke2023optimization} study freight sharing rates to support adjustments in freight structures to meet carbon peak goals. They use an adaptive genetic algorithm (GA) to balance emissions with transportation utility values, such as safety, price, speed, and flexibility.

\noindent \textbf{GA-Based Optimization for MMT:} Both \cite{zhang2023research} and \cite{yang_multi-objective_2023} address multi-objective path optimization for container MMT using GA. \cite{zhang2023research} utilize a modified GA combined with the TOPSIS method to optimize MMT paths, considering factors such as cost, emissions, logistics service quality, and mode change times to improve cargo transport efficiency. In contrast, \cite{yang_multi-objective_2023} implement an improved fuzzy adaptive GA, which adjusts crossover and mutation probabilities based on population variance, to optimize MMT paths for cost, time, and emissions under multi-task conditions from the perspective of MMT operators.

Similarly, \cite{li2023hierarchical} use an adaptive GA to optimize secondary hub locations, cargo flow, and mode selection in a hierarchical multimodal hub-and-spoke network. Their model considers hub capacity, cargo clearance efficiency, and carbon tax implications. As carbon emission taxes increase, the optimized network sees more railway links (offering lower emissions but slower speeds) and additional secondary hubs (which reduce transit time for goods). \cite{shoukat2023sustainable} apply GA to optimize logistics networks linking dry ports and seaports. They introduce a bi-objective MILP model, building on \cite{demir2019green}'s green logistics framework, which initially focused on the intermodal service network within dry ports.  \cite{shoukat2023sustainable} address the nuanced differences between multimodality and intermodality, addressing potential inconsistencies in these terms' usage across different sections of their research.

\noindent \textbf{Uncertainties and Dynamic Constraints: }
Recent studies explore various challenges and strategies for optimizing IMT amidst uncertainties and dynamic conditions.
\cite{sun2023modeling} focus on first-mile pickup and last-mile delivery, using soft time windows and accounting for uncertainties in truck speeds and rail capacities while integrating carbon tax regulations into transportation path planning. \cite{li2023optimum} 
study MMT path planning considering uncertain customer demand, transportation volume, and transshipment times. Their model uses mixed-time window constraints and combines GT with the weighted sum method to dynamically adjust objective weights during optimization under uncertainty.

Expanding on the theme of dynamic factors, \cite{liu2023multimodal} explore cold chain container routing, taking into account variables such as traffic congestion, weather conditions, and vehicle breakdowns that directly impact transport speed. This study emphasizes the importance of managing penalty costs and adhering to stringent customer time windows, with a focus on reducing carbon emissions based on actual travel times. Additionally, \cite{kurtulucs2023green} present a container shipping scheduling model focused on recovering from disruptions (e.g., speed adjustment and port skipping) in port and sea passage. While port skipping is found to be effective, its impact diminishes significantly when shipping lines and terminal operators share real-time information and coordinate custom time windows and handling rates.

\noindent \textbf{Influence of Shippers' Choice:} Recent research highlights the pivotal role of shippers' choices in shaping sustainable transportation strategies.
\cite{sun2023study} optimize river channel upgrades and ship operations, aiming to increase river shipping's share in container port access. The proposed bi-level programming model integrates strategic decisions on infrastructure and ship choices with route optimization, emphasizing the role of carbon taxes in promoting environmentally sustainable transportation solutions and reducing reliance on trucking for port logistics.

\cite{wu2023evaluation} further examine the broader scope of shippers' decision-making processes, including selecting transport modes, inland terminals, and ports. They focus on how transport capacity expansion, carbon taxes, and subsidies influence shifts in transport modes, alongside impacts on carbon emissions and costs. Additionally, \cite{guo2023integrated} analyze a multimodal feeder shipping network, focusing on the impact of shippers' decisions, and propose a nonlinear model to optimize route selection, schedule design, and fleet allocation with consideration for carbon emission costs.

\noindent \textbf{Bulk Cargo Distribution: }Strategic enhancements in bulk cargo distribution are paramount, particularly through the integration of MMT solutions that consider physical constraints and customer preferences. \cite{ko2022stochastic} examine the impact of modal shifts but within the context of biomass SCs under transportation cost stochasticity, where road transport prevails due to limited rail infrastructure. Their model incorporates transportation costs, including railcar leasing, as well as environmental costs from GHG emissions and societal costs. Next, \cite{de2023flexible} introduce a DES model designed to facilitate collaborative decisions among multiple stakeholders (i.e., shippers and consignees). Their model evaluates seven transportation strategies, helping stakeholders choose between using a seaport or dry port, opting for road or multimodal routes (combining truck and rail), and deciding on consolidating or deconsolidating cargo containers. Expanding on the logistical complexities of bulk cargo, \cite{feng2023multimodal} address the distribution challenges within MMT networks, particularly focusing on the transition from bulk to containerized transport. This study integrates customer preferences and includes inland waterway transport as a viable mode, considering physical constraints like bridge heights to optimize mode selection, vehicle routing, depot choice, and cargo containerization strategies.

\noindent \textbf{Most Recent Developments: }
Looking forward to 2024, studies by  \cite{yin2024low} and \cite{zhang2024sparrow} compare different carbon pricing strategies.  \cite{yin2024low} evaluate combined policies, such as rail freight subsidies and a carbon trading mechanism. They point out the limited success in modal shifts and the lack of differentiation between regional carbon pricing mechanisms. Their analysis includes a multi-objective IMT ND that integrates transportation costs, travel times, and carbon trading. As for  \cite{zhang2024sparrow}, they develop a multi-objective approach to optimize MMT planning, to minimize transportation costs and time, considering fuzzy demand and time intervals with various carbon policies.

Additionally, \cite{ghisolfi2024dynamics} and \cite{derpich2024pursuing} contribute to modal and energy strategy shifts.  \cite{ghisolfi2024dynamics} use an SD model to evaluate the impact of electrification, increased biofuel usage, accelerated fleet renewal, and modal shifts on carbon emissions within the Brazilian freight system. Meanwhile, \cite{derpich2024pursuing} tackle an MMT ND problem and adapt the hub-and-spoke strategy to fit regional demand characteristics, focusing on consolidation and concentration in hub selection.

\section{Integrated Analysis of Literature} \label{Integrated Analysis of Literature}

In this section, we examine general trends in the literature related to modality mix, logistics decision levels, and emissions considerations in OR models. We also provide an overview of the characteristics of OR model formulations and the solution techniques developed in this field. We consider studies from 2010 to 2023 for this analysis, intentionally excluding 2024 due to the limited number of studies (only 4) published that year, which could skew trends or insights.

\subsection{Modality Analysis}
\begin{figure}[!htp]
    \centering
    \includegraphics[width=0.80\textwidth]{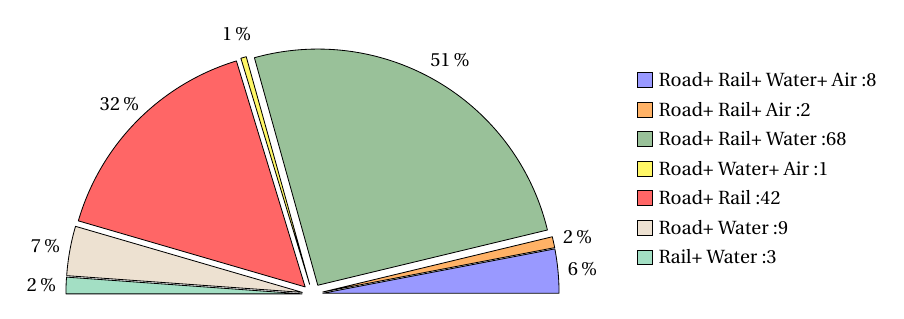}
    \caption{Classification of Publications by Transportation Modes}
    \label{fig:ModalityMix}
\end{figure}

For analyzing the distribution of research across various transportation modalities, it is crucial to understand how different combinations of transport modes are represented in the literature. To provide a comprehensive overview, we examine the modality mix in Fig.~\ref{fig:ModalityMix}, which presents seven distinct combinations of transportation modes, illustrating the focus and frequency of studies within each combination.

The combination of road, rail, and water transport is the most extensively studied, with significant research focused on SND. This combination saw a notable increase in studies in 2023, driven by global efforts to reduce carbon emissions and recent pandemic-related SC disruptions, highlighting the need for more sustainable and resilient transport systems. Following this, the road and rail combination is also frequently explored, particularly in the context of hinterland connections like dry ports, which facilitate modal shifts at inland terminals. This pairing leverages the complementary strengths of each mode: rail offers efficiency and cost savings for long distances, while road provides flexibility and time efficiency for short distances. In contrast, air transportation is less studied, with approximately one study published annually, due to its high costs, weight and volume restrictions, and handling challenges. Additionally, some research addresses Global SC ND, focusing on critical decision-making factors such as travel time and risk, including the potential for freight damage.

\begin{figure}[!htp]
    \centering
    \includegraphics[width=0.85\textwidth]{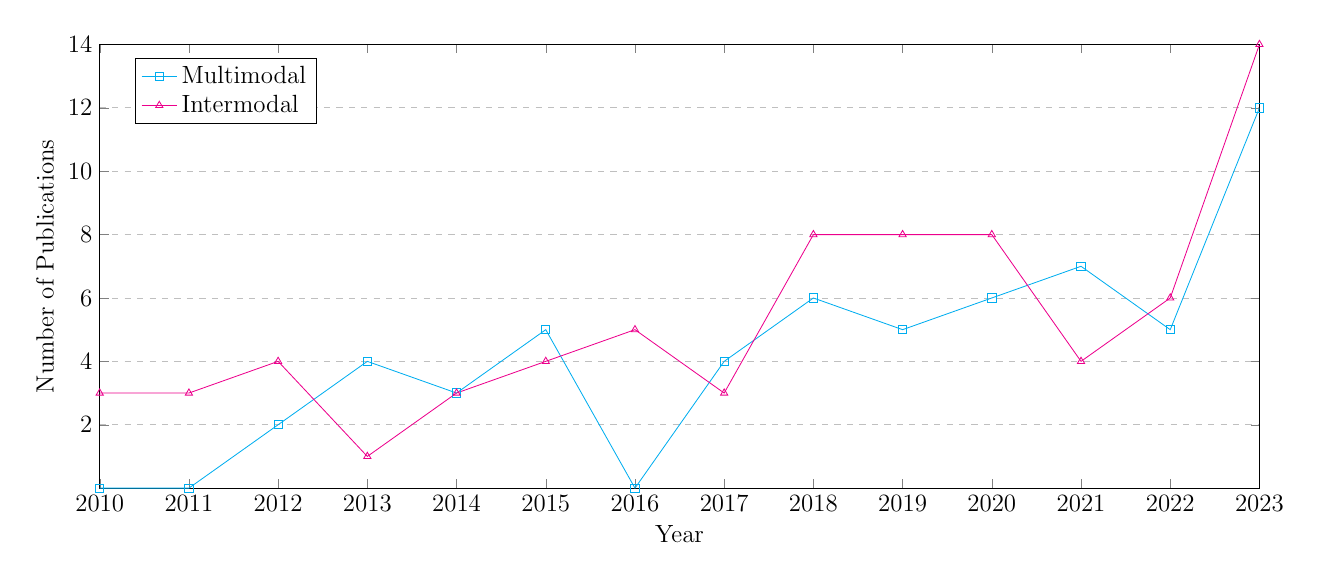}
    \caption{Number of Publications on Intermodal  vs Multimodal Transportation Over the Years}
    \label{fig:IMTvsMMT}
\end{figure}

Next, we present the evolution of IMT and MMT studies in Fig.~\ref{fig:IMTvsMMT}. The volume of publications on IMT has shown a continuous increase throughout the research period, reflecting its growing importance in improving the financial and environmental dimensions of FT.  However, in 2021, the pandemic shifted priorities toward healthcare logistics and emergency response, leading to a temporary dip in IMT research. By 2022 and 2023, attention returned to optimizing SCs, with IMT’s potential to enhance efficiency and resilience gaining renewed interest. The surge in both MMT and IMT research in 2023 can be attributed to several factors, including post-pandemic recovery, a renewed focus on SC resilience, and evolving policy changes. For this figure, research on SMT, considered an extension of IMT, has been incorporated into the IMT publication dataset.

\begin{figure}[!htp]
    \centering
    \includegraphics[width=0.85\textwidth]{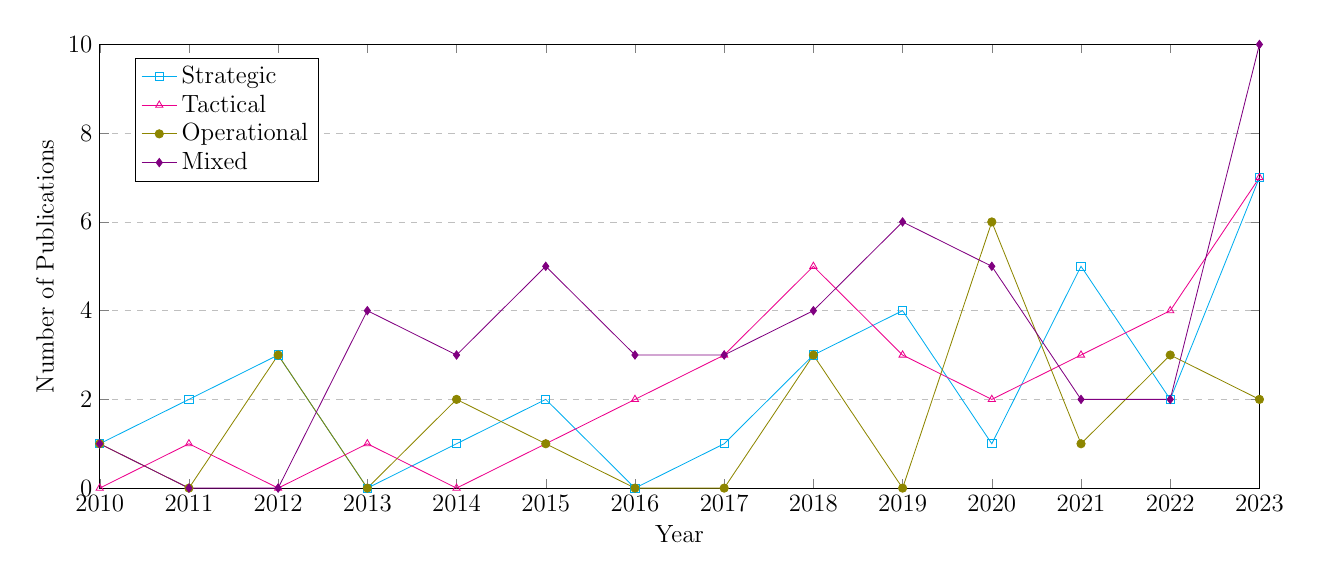}
    \caption{Number of Publications Regarding Planning Decision Levels Over the Years}
    \label{fig:LogisticDecisionsAll}
\end{figure}

\begin{figure}[!htp]
    \centering
    \includegraphics[width=0.85\textwidth]{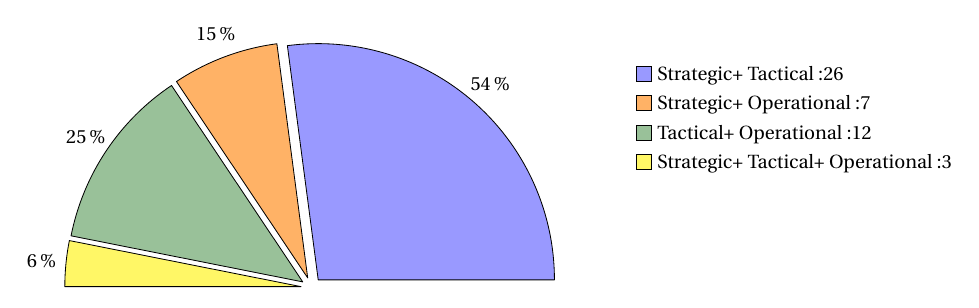}
    \caption{Publication Overview of Mixed Planning Decision Levels}
    \label{fig:MixLogisticDecisions}
\end{figure}

Planning decisions can be divided into strategic, tactical, and operational levels \citep{crainic1997planning}. Based on this classification, we organize the studies as shown in Fig.~\ref{fig:LogisticDecisionsAll} by time. We identify a trend where multiple decision levels are often addressed simultaneously, leading us to introduce a Mixed category. This category encompasses the highest number of studies, therefore we present the various combinations in Fig.~\ref{fig:MixLogisticDecisions}. This trend reflects the growing recognition that isolating decision levels may not adequately capture real-world scenarios' dynamic interactions and trade-offs. 

Interestingly, there has been a marked decline in research focusing solely on operational planning over the past three years. This suggests a shift towards integrating operational planning with other decision levels. Tactical decisions are the most frequently examined, followed by strategic decisions, indicating a strong focus on deterministic and stochastic ND studies. The observed decline in operational planning studies, coupled with the increased emphasis on strategic and tactical research, signals a broader shift toward long-term planning and resilience.

\subsection{Sustainability Analysis} \label{Sustainability}

As sustainability becomes an increasingly central focus in OR, examining how carbon emissions are integrated into OR formulations is essential. This section explores the various approaches used to account for carbon emissions, categorizing them into five main classifications, as shown below, each representing different strategies for incorporating emissions considerations into OR models. 
\begin{enumerate}
    \item \underline{Carbon estimation}: Using emission estimations or factors that are minimized within OR models.
    \item  \underline{Carbon cap}: An emissions limit, often regulated by authorities.
    \item  \underline{Carbon cost}: Assigning a monetary value charged per unit of emissions (e.g., taxes).
    \item  \underline{Cap-and-trade}: Carbon credits that can be bought or sold in a market regulated by a third party.
    \item  \underline{Combination}: Any mix of the above methods.
\end{enumerate}

\begin{figure}[!htp]
    \centering

     \includegraphics[width=0.85\textwidth]{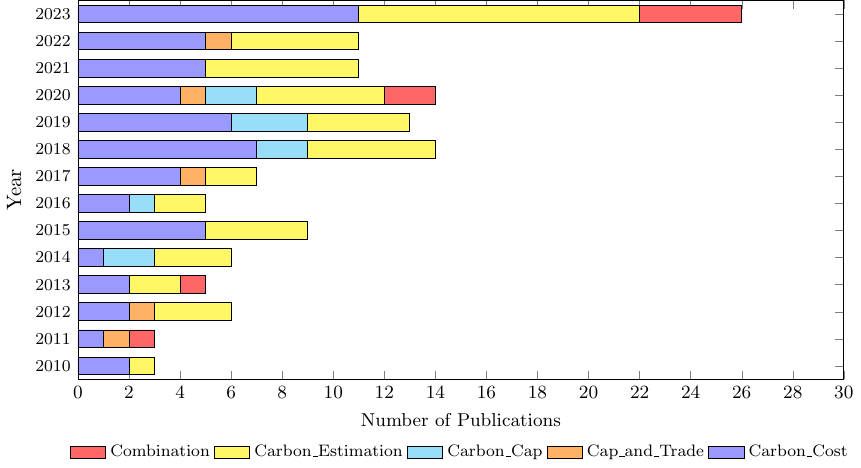}

    \caption{Number of Studies By Carbon Emissions Consideration in OR Models}
    \label{fig: CO$_2 $}
\end{figure}

In Figure \ref{fig: CO$_2 $}, we show the number of studies by each classification over time. Carbon cost (43\%) and carbon estimation (40\%) are the most explored throughout the study period. 
We note that in the early years, the location and selection of hub terminals within the transportation network were largely influenced by carbon cost modeling. However, recent trends indicate a shift toward emphasizing the importance of enhancing network-wide operating efficiency in decision-making.
Carbon estimations played a crucial role in infrastructure expansion decisions within the transportation network in the early years and have recently shifted to assess the environmental impacts of modal shift strategies. 

Another trend we have observed in recent years is the increase in research focused on enhancing the precision of emissions estimation models. Cap-and-trade is primarily used in decisions related to SCND, while carbon caps are initially applied to assess government decisions. However, recently, the focus has shifted toward analyzing a combination of emissions strategies, driven by the increasing need for more impactful approaches from governments to decarbonize the sector.

It is crucial to note that these strategies have varying degrees of real-world application. For instance, carbon estimation is widely used by companies to track and reduce their emissions, often guided by frameworks such as the GHG Protocol \citep{green2010private}. Carbon cost, in the form of carbon taxes, is implemented in 30 jurisdictions as of 2020 \citep{metcalf2021carbon}, encouraging businesses to reduce emissions to avoid financial penalties. Cap-and-trade systems, as seen in regions like California, enable companies to trade carbon credits, providing incentives for emissions reduction \citep{lessmann2024effect}.

While many of these approaches are actively implemented by academia, particularly the OR community, often explore idealized versions focused on optimizing factors like cost, emissions, and operational efficiency. However, real-world implementation is more complex due to challenges such as profit and revenue loss, fluctuating carbon prices, regulatory differences, and technological constraints. Companies may voluntarily adopt these strategies, especially in regions with strong environmental regulations or market incentives, but compliance is more commonly driven by mandatory policies.

\subsection{OR Techniques Analysis}

To examine the studies through the lens of OR techniques, we classify the studies based on two key factors: the problem being addressed and the OR technique employed to tackle that problem. We then explore the distribution of these techniques across different problem categories and summarize the key factors of the primary methods being utilized in tables. Emission estimation models (6 studies) are excluded from this section, as they do not involve OR techniques but were previously analyzed for their relevance to IMT decarbonization. 

Before diving into an extensive analysis,  we present an overview of the mathematical models discussed in the articles. This includes their objective functions, classifications as deterministic or stochastic, single or multiple objectives, and details of the methods used to solve these problems. This information is available in Table \ref{tab:Solution} of \ref{Comprehensive_Analysis}. The table indicates that economic objectives dominate the literature, with about 75\% of studies addressing transportation costs, while emission and handling costs are considered in approximately 40\%. In contrast, transportation time and emissions are objectives in only 17\% and 33\% of studies, respectively. However, there has been a noticeable trend over the last five years towards increasing consideration of transportation time as an objective.

\subsubsection{Classification based on problem addressed:}

To illustrate the application of OR techniques to various problems, we categorize the studies by the problems addressed and the OR techniques used. This section delineates the problems into five categories: a) Network Design (NDP),  b) Service Network Design (SNDP), c) Supply Chain Network Design (SCND), d) Policy Support (PS), and e)  Mixed category representing combinations of these problems. We now provide a brief description of each primary category.  NDP focuses on optimizing transportation infrastructure by the construction or improvement of facilities, as well as integrating decisions related to the allocation of resources with the aim to improve route efficiency, cost, and capacity of the network \citep{liu2014global}. 
Within this category, we identify three main sub-categories: Location Problem (LP), Capacity Expansion Problems (CEP), and Hub Selection Problem (HSP). Next, SNDP involves planning the selection, routing, and scheduling of transportation services, as well as managing terminal operations and freight routing, ensuring the fulfillment of the demand and profitability \citep{crainic2000service}. Four main sub-categories are included in this regard \citep{wieberneit2008service}: Service Selection Problem (SSP), Traffic Distribution Problem (TDP), Asset Management Problem (AMP), and Revenue Management Problem (RMP). While SCND focuses on determining the optimal location and size of facilities and managing the flow of goods through these facilities \citep{moreno2019sustainability}, the PS category encompasses studies that aim to contribute to the decision-making processes of stakeholders by evaluating strategies across various decision levels (e.g., tactical) \citep{banister2011transportation}. 
These studies offer insights and evidence-based assessments for policy formulation, implementation, and optimization in transportation.

\subsubsection{Classification based on OR techniques:}

Correspondingly, we classify the OR techniques into eight categories: Exact (e.g., branch \& bound, simplex), Heuristic, Meta-heuristic, Simulation, Hybrid (i.e., the combination of at least two techniques), Game Theory (GT), Decision Support System (DSS), and Others. The ``Others" category encompasses less commonly used methodologies, such as Markov Decision Processes (MDP) \citep{puterman1990markov}, Lagrangian Relaxation (LR) \citep{lemarechal2001lagrangian}, and Benders Decomposition (BD) \citep{geoffrion1972generalized}. Next, we offer a concise overview of OR categorization. Exact methods, like branch-and-bound and simplex \citep{morrison2016branch}, guarantee optimal solutions but can be time-consuming, especially as the problem size increases. For simplicity, we separate exact solution algorithms such as BD and Column Generation from this category. Heuristics provide quick, problem-specific solutions without ensuring optimality \citep{muller1981heuristics}. Meta-heuristics guide and modify heuristics to explore search spaces more effectively and avoid local optima \citep{boussaid2013survey}. Simulation techniques, including DES \citep{agalianos2020discrete}, ABM \citep{shepherd2014review} and SD \citep{huang2022overview}, are used for evaluating and optimizing systems through various what-if scenarios \citep{crainic2018simulation}.  Each method offers unique strengths: DES in event sequencing, SD in feedback understanding, and ABM in modeling intricate agent interactions. While hybrid methods either combine optimization and simulation to handle complex systems where social factors complicate precise modeling \citep{amaran2016simulation} or employ exact methods and heuristics together, GT models interactions among decision-makers and is often paired with heuristics for large-scale problems  \citep{cachon2006game}. On the other hand, DSSs improve decision-making by enhancing various performance metrics in transportation \citep{caris2013decision}.

\begin{table}[!htbp]
    \centering
    \caption{Distribution of OR Techniques Across Different Problem Categories in IMT}
    \label{tab:ORByProblem}
    \footnotesize

    \begin{tabularx}{\textwidth}{l c *{9}{>{\centering\arraybackslash}X}}
        \toprule
        \textbf{Category} & & \textbf{Exact} & \textbf{Heuristic} & \textbf{Meta-heuristic} & \textbf{Simulation} & \textbf{Hybrid} & \textbf{Game Theory} & \textbf{DSS} & \textbf{Others} \\ 
        \midrule
        \textbf{NDP}      & & 8  & 2  & 2   & 1  &   & 3  &    &   \\
        \textbf{}         & LP & 3  &    & 1   &    &  1 & 3  &    &   \\
        \textbf{}         & CEP & 3  & 1  & 1   &    &   &    &    &   \\
        \textbf{}         & HSP & 2  & 1  &     & 1  &   &    &    &   \\
        \midrule
        \textbf{SNDP}     & & 35 & 1  & 19  & 3  & 6 & 1  & 6  & 3 \\
        \textbf{}         & SSP & 17 &    & 8   & 1  & 5 &    & 2  & 2 \\
        \textbf{}         & TDP & 12 &    & 8   & 1  &   &    & 4  &   \\
        \textbf{}         & AMP & 6  &    & 3   & 1  &   &    &    &   \\
        \textbf{}         & RMP &    & 1  &     &    & 1 & 1  &    & 1 \\
        \midrule
        \textbf{SCND}     & & 9  & 1  & 2   & 1  & 1 &    & 1  & 1 \\
        \midrule
        \textbf{PS}       & & 1  & 1   & 3   & 5  &   & 3  & 3  & 1 \\
        \midrule
        \textbf{Mixed}    & & 1  &    &  1   &    & 1 &    &    &   \\
        \midrule
        \textbf{Total}    & & 54 & 5  & 27  & 10 & 9 & 7  & 10  & 5 \\
        \bottomrule
    \end{tabularx}
    \begin{threeparttable}
        \begin{tablenotes}
            \footnotesize
            \item \textbf{DSS}: Decision Support System, 
            \textbf{NDP}: Network Design Problem, 
            \textbf{LP}: Location Problem, 
            \textbf{CEP}: Capacity Expansion Problem, 
            \textbf{HSP}: Hub Selection Problem, 
            \textbf{SNDP}: Service Network Design Problem, 
            \textbf{SSP}: Service Selection Problem, 
            \textbf{TDP}: Traffic Distribution Problem, 
            \textbf{AMP}: Asset Management Problem, 
            \textbf{RMP}: Revenue Management Problem, 
            \textbf{SCND}: Supply Chain Network Design, 
            \textbf{PS}: Policy Support.
        \end{tablenotes}
    \end{threeparttable}
\end{table}

\subsubsection{Distribution of OR Techniques Across Different Problem Categories:} 
Next, we provide in depth analysis in Table \ref{tab:ORByProblem},  detailing the distribution of studies across the problem types and OR techniques.

\noindent  \textbf{Exact Methods:}  The application of these OR techniques for the decarbonization of FT is primarily led by exact methods. In this context, the majority of the problems are related to SNDP, specifically SSP, where decisions regarding the structure of transportation services are the main focus. Most studies have traditionally concentrated on deterministic problems, with significant emphasis placed on decision-making at multiple levels, particularly strategic and tactical, with tactical decisions receiving secondary focus. However, since the first study to incorporate uncertain factors using exact methods in 2016 \citep{demir2016green}, we note two key aspects. The first is a shift from generalized cost considerations in the objective function (e.g., emission costs, penalty costs) to the inclusion of additional factors, such as transportation time and emissions. The second trend is the significant increase in research output, with 57\% of the studies being developed between 2019 and 2020. 

As the field has evolved, in the last few years, researchers have increasingly incorporated three key aspects: sustainability trade-offs between stakeholders, non-linearity sources (e.g., traffic congestion, pricing decisions), and uncertainty (e.g., demand). Trade-offs are often addressed through multi-objective formulations, with common reformulation methods including the epsilon constraint, weighted sum, and others. Among the 54 studies utilizing exact methods, 27 address multiple objectives, with 11 using the epsilon constraint method. This highlights the clear dominance of this approach in tackling multi-objectivity. For non-linear formulations, a common approach is to add auxiliary variables and linear constraints, allowing these problems to be solved with commercial solvers (e.g., CPLEX).

\begin{table}[!htbp]
    \centering
    \caption{Summary of Papers Using Meta-heuristic Approaches }
    \label{tab:MetaHeuristic}
    \footnotesize

    \resizebox{\textwidth}{!}{
        \begin{tabularx}{\textwidth}{p{1.6cm} c p{2.9cm} c p{1.2cm} p{6.7cm}}
            \toprule
            \textbf{Article} & \textbf{Problem} & \textbf{Decision(s)} & \textbf{OT} &  \textbf{Algorithm} & \textbf{Key Characteristics} \\
            \midrule
            \citet{sawadogo2012reducing} & SNDP & Determine sustainable path & MO  & ACO & Defines heuristic information `a priori' and updates dynamically. Updates pheromone trails based on transportation mode. \\
            \midrule
            \citet{kim2013multimodal} & NDP & Upgrade facility and assign traffic flow  & BLO  & GA-based & Performs feasibility checks with budget constraint and neighborhood improvement based on volume-to-capacity ratios. \\
            \midrule
            \citet{chen2014design} & PS & Define state subsidies and transport service design & BLO  & GFWHA & Alternates selection methods. Includes carbon reduction in fitness function, and shipping routes in chromosome coding. \\
            \midrule
            \citet{fahimnia2015tactical} & SCND & Choose production and distribution plans under carbon tax policies & MO & NICE & Considers adaptive update for binary variables, iterative update for feasible solutions until the ratio drops below a set value. \\
            \midrule
            \citet{ji2017hybrid} & SNDP & Optimize service selection & MO & HEDA & Uses heterogeneous marginal distribution law as probability model. Embeds local search methods. \\
            \midrule
            \citet{zhang2017multimodal} & Mixed & Optimize service selection and node selection & SO & Hybrid GA & Considers two-part chromosome (node selection and mode selection). Applies two-point crossover, two-point mutation, and customized mutation.  \\
            \midrule
            \citet{zhao2018stochastic} & SNDP & Reposition empty containers under uncertainty & SO & Two-phase TS & Uses swap move sample, dynamic tabu duration and aspiration criterion for both, relocation and routing phase. \\
            \midrule
            \citet{zhang2018joint} & PS & Optimize investment, subsidy value, and flow assignment & BLO & GFWHA & Considers two-part chromosome (infrastructure selection and subsidy values). Uses two-point crossover, customized crossover, two-point mutation, and customized mutation. \\
            \midrule 
            \citet{maiyar2019modelling} & SNDP & Choose hub terminal locations and flow assignment under disruptions & SO & PSO with DE & Conducts particle encoding with population size and decision variables. Uses a discretization scheme for boundary violations. \\
            \midrule
            \multirow{2}{*}{\citet{lu2019sustainable}}   & \multirow{2}{*}{SNDP}  & \multirow{2}{*}{Optimize service selection} & \multirow{2}{*}{SO} & DE & Addresses the railway transport echelon. \\
            & &   &  &  CW+LS & Addresses the highway transport echelon. Uses CW for initial solution, LS for optimal solution (with four operators). 
            \\
            \midrule
            \citet{wei2019import} & SNDP & Allocate freight flow in cross-border network & MO  &  awGA & Uses matrix integer encoding (five layers), and adaptive weight fitness function with elitism strategy for offspring generation. \\
            \midrule
            \citet{wang2020integrated} & SCND & Schedule production and transportation & SO  & Improved GA & Uses a multi-heterogeneous coding method with heuristic rules (two-part chromosome: production and transportation).  \\
            \midrule
            \citet{wu2021research} & SNDP & Optimize logistics network under 4PL scheme & SO & Improved PSO & Uses a constraint processing mechanism with three algorithms for flow allocation, capacity limit, and capacity adjustment. \\
            \midrule
            \citet{liang2021multi} & SNDP & Optimize service selection with user satisfaction & MO & NSGA-III & Uses priority-based real-number chromosome encoding, a simulated binary crossover and polynomial mutation operator. \\
            \midrule
            \citet{li2022path}         &  SNDP    & Obtain sustainable path under uncertainty 
            &  SO  & FAGSO        & Considers Gaussian mutation sparks and modulo operation for mapping rules. Uses attractive force search operators. \\
            \midrule                
        \end{tabularx}
 }
    \begin{threeparttable}
        \begin{tablenotes}
            \footnotesize
            \item \textbf{OT}: Optimization Type, \textbf{SO}: Single-objective Optimization, \textbf{MO}: Multi-objective Optimization, \textbf{BLO}: Bi-level Optimization, \textbf{ACO}: Ant Colony Optimization, \textbf{GA}: Genetic Algorithm, \textbf{GFWHA}: Genetic and Frank–Wolfe Hybrid Algorithm, \textbf{NICE}: Nested Integrated Cross-Entropy, \textbf{HEDA}: Hybrid Estimation of Distribution Algorithm, \textbf{TS}: Tabu Search, \textbf{PSO}: Particle Swarm Optimization, \textbf{DE}: Differential Evolution, \textbf{CW}: Clarke–Wright Savings Algorithm, \textbf{LS}: Local Search, \textbf{awGA}: Adaptive-weight GA, \textbf{NSGA-III}: Non-dominated Sorting GA III, \textbf{FAGSO}: Fireworks Algorithm with Gravitational Search Operator, \textbf{NSGA-II}: Non-dominated Sorting GA II, \textbf{ALNS}: Adaptive Large Neighborhood Search Algorithm, \textbf{HEGA}: Hummingbird Evolutionary GA, \textbf{TOPSIS}: Technique for Order Preference by Similarity to an Ideal Solution.
        \end{tablenotes}
   \end{threeparttable}
\end{table}

\begin{table}[!htbp]\ContinuedFloat
    \centering
    \caption{(Cont.) Summary of Papers Using Meta-heuristic Approaches }
    \label{tab:MetaHeuristic_2}
    \footnotesize 

    \resizebox{\textwidth}{!}{
        \begin{tabularx}{\textwidth}{p{1.6cm} c p{2.9cm} c p{1.2cm} p{6.7cm}}
            \toprule
            \textbf{Article} & \textbf{Problem} & \textbf{Decision(s)} & \textbf{OT} & \textbf{Algorithm} & \textbf{Key Characteristics} \\
            \midrule 
            \citet{xie2022research} & SNDP & Locate air hubs in cross-border logistic & MO & NSGA-II & Uses integer encoding,  crowding distance calculation, simulated binary crossover, and polynomial mutation. \\
            \midrule 
            \citet{zhang2022preference} & SNDP & Develop pick-up and delivery SMT plans & MO & ALNS & Customizes operators for node and route removal. Ensures synchronization procedures. Addresses constraint violations. \\
            \midrule 
            \citet{zhang2022flexible} & SNDP & Develop routing plans for flexible SMT service & SO & ALNS & Predefines initial solutions. Customizes swap operator. Utilizes simulated annealing for acceptance criterion. \\
            \midrule 
            \citet{zhang2022collaborative} & SNDP & Develop collaborative SMT & MO & ALNS & Considers preference constraints (higher load assignment to trains and barges). \\
            \midrule 
            \citet{zhang2022heterogeneous} & SNDP & Optimize SMT planning & MO & ALNS & Prioritizes requests based on preference. Conducts synchronization checks for time and preferences. \\
            \midrule
            \citet{ke2023optimization} & PS & Adjust freight structure  & MO & AwGA & Uses floating point for chromosome encoding, roulette selection, adaptive crossover and mutation rates. \\
            \midrule
            \citet{shoukat2023sustainable} & SNDP & Optimize service selection & MO & GA & Uses random generator for initial population, binary chromosome coding. Addresses two objective functions in the fitness function simultaneously. \\
            \midrule
            \citet{liu2023multimodal} & SNDP & Optimize service selection under uncertainty & SO & HEGA & Uses insertion heuristic for initial population and queen-bee evolution method for crossover. Implements elite retention strategy. \\
            \midrule
            \citet{zhang2023research} & SNDP & Optimize service selection & MO & Modified GA & Utilizes hierarchical encoding, adaptive crossover and mutation mechanism. Conducts fitness evaluation using TOPSIS.  \\
            \midrule
            \citet{li2023hierarchical} & NDP & Determine secondary hub locations & MO & Adaptive GA & Considers three-part chromosome coding for freight allocation, with multi-point crossover and single-point mutation. \\
            \midrule
            \citet{yang_multi-objective_2023} & SNDP & Optimize service selection & MO & Improved fuzzy GA & Uses fuzzy controller to dynamically adjust GA parameters based on population variance and fitness values. \\ 
            \midrule
            \citet{guo2023integrated} & SNDP & Optimize service selection & SO & PSO & Considers particle swarm optimization with path selection ratio, as well as inbound and outbound specification. \\
            \midrule
        \end{tabularx}
    }
\end{table}

\noindent  \textbf{Meta-heuristics:} The second most preferred approach, after exact methods, is meta-heuristics. In Table \ref{tab:MetaHeuristic}, we present a comprehensive overview of studies employing this solution approach. It highlights key aspects such as problem type, decision(s), and optimization type (single-objective, multi-objective, or bi-level optimization). Additionally, it provides the name of the algorithm and 
outlines the unique characteristics of each algorithm, emphasizing the specific adaptations or innovations introduced in each study. Researchers have consistently employed meta-heuristics, with a significant uptick in studies over the past two years, accounting for around 50\% of all recent research. This trend indicates a growing and sustained interest in these techniques. We believe this surge is due to the flexibility and adaptability of this technique. Notably, a significant portion of the studies (63\%) focus on solving multi-objective and bi-level optimization problems. This is well-founded, given their ability to efficiently navigate vast solution spaces, balance exploration and exploitation, and address complex, conflicting objectives. 

A unique case is displayed by \citet{lu2019sustainable} where the model is decomposed into two echelons, each one solved through a specific meta-heuristic algorithm. 
Regarding the problem type, we can observe that the majority of the studies are considered SNDP, specifically SSP, and TDP.  
Among the 16 studies conducted between 2020 and 2023, 11 employed GA, PSO, or some modified versions of these techniques, underscoring their relevance in current research. GAs are typically applied in the context of service selection, often including tasks such as determining service frequency or coordinating between services with fixed and flexible schedules. We can also observe that ALNS solves all SNDPs in the context of SMT.  We attribute this to the nature of decisions requiring collaboration between stakeholders and the need for real-time solutions, where the fast convergence of ALNS is particularly well-suited.

\noindent  \textbf{Simulation Methods:} Simulation techniques have gained considerable prominence, with 90\% of studies conducted in the past five years and 50\% occurring in 2023 alone utilizing these methods. Table \ref{tab:Sim} presents a comprehensive overview of the diverse simulation methodologies utilized to analyze IMT systems and policies. The table is organized into six columns—article, approach, objective, decision level, transportation modes, and findings—outlining each study's simulation method, focus, decision level, modes considered, and key outcomes. Based on the table, it is evident that the most commonly employed approach is SD, with all studies utilizing this method focused on solving PS problems at the strategic decision level. Additionally, ABM has consistently been used in combination with DES. In terms of decision levels, SD primarily addresses strategic decisions, while DES focuses on operational decisions. Strategic decisions are more frequently explored through simulation compared to operational and tactical decisions. Interestingly, only one study considers all three decision levels, underscoring the need for further research in this area. Road and rail modes are included in every study, while only half of the studies consider waterways.

\begin{table}[!htbp]
    \centering
    \caption{Summary of Papers Using Simulation }
    \label{tab:Sim}
    \footnotesize

    \resizebox{\textwidth}{!}{
        \begin{tabularx}{\textwidth}{p{1in} p{0.4in} p{1.2in} c c p{1.46in}}
            \toprule
            \textbf{Article} & \textbf{Approach} & \textbf{Decision(s)} & \textbf{Decision Level} & \textbf{Transp. Modes} & \textbf{Findings} \\
            \midrule
            \cite{holmgren2012tapas} & ABM + DES & Transport policy and infrastructure analysis & S & RW, RR, WW & Captures goals, interactions, and time aspects of transport chain actors. \\
            \midrule
            \cite{lattila2013hinterland} & DES & Estimate seaport vs dry port usage effects & O & RW, RR & Dry ports reduce emissions and costs. \\
            \midrule
            \cite{liu2019system} & SD & Evaluate long-term effects of policies on emissions & S & RW, RR & Weight regulations and grade 1 railway construction reduce emissions. \\
            \midrule
            \cite{choi2019system} & SD &  Containerization and tax effect on modal shift & S & RW, RR & Containerization leads to faster modal shifts than taxation. \\
            \midrule
            \cite{wang2020modelling} & DES & Evaluate multimodal transport route efficiency & O & RW, RR, WW & Reorganizing routes and setting up container centers cuts costs. \\
            \midrule
            \cite{nassar2023system} & SD & Modal shift policies for freight decarbonization & S & RW, RR, WW & Stricter policies and investments accelerate shift to rail. \\
            \midrule
            \cite{de2023flexible} & DES & Compare seaport/dry port and multimodal efficiency & S, O & RW, RR & Dry ports are cost-efficient for long-term storage.  \\
            \midrule
            \cite{zhong2023system} & SD &  Carbon taxation impact on port-hinterland transport & S & RW, RR, WW & Medium-high taxes shift traffic to rail and waterways, reducing emissions. \\
            \midrule
            \cite{guo2023toward} & SD + MC & Design policy mixes to enhance sustainability & S, T & RW, RR, WW & Cost-based pricing and better rail services are effective in mixes. \\
            \midrule
            \cite{castrellon2023assessing} & ABM + DES & Assess dry port strategies with different container substitution levels & S, T, O & RW, RR & Full ownership container substitution and extended free storage reduce empty container repositioning. \\
            \bottomrule
        \end{tabularx}
    }
    \begin{threeparttable}
        \begin{tablenotes}
            \footnotesize
            \item \textbf{ABM}: Agent-based model, \textbf{DES}: Discrete event simulation, \textbf{SD}: System dynamics, \textbf{MC}: Monte Carlo, \textbf{RW}: Road, \textbf{RR}: Railway, \textbf{WW}: Waterway, \textbf{S}: Strategic, \textbf{T}: Tactical, \textbf{O}: Operational
        \end{tablenotes}
    \end{threeparttable}
\end{table}

\noindent  \textbf{Game Theory:} GT has been employed by researchers for mostly analyzing NDP and PS problems. Table \ref{tab:GT} summarizes research on GT applications in these areas. A recurring theme in these studies is the interaction between a leader and followers, where the leader sets policies or makes strategic decisions, and the followers respond to these actions. The table highlights this interaction by identifying both the leaders and followers, outlining the leaders' actions, and detailing the subsequent responses of the followers. It also indicates the impact of these decisions across economic, social, and environmental dimensions. In 50\% of the studies, the government takes on the role of the leader, implementing policies or setting prices that influence the behavior of firms, ports, or shippers to adopt more sustainable practices. Another important insight is that nearly all studies examining the government's role as a leader were conducted last year, indicating a clear trend. 
We can see a strong emphasis on economic and environmental decisions, particularly regarding pricing strategies, transportation efficiency, emission reductions, and sustainability. Decisions related to social factors are addressed in only about half of the studies, being relevant primarily in the context of infrastructure planning.

\begin{table}[!htbp]
    \centering
    \caption{Summary of Papers Using Game Theory Approaches }
    \label{tab:GT}
    \footnotesize

    \resizebox{\textwidth}{!}{
        \begin{tabularx}{\textwidth}{p{2 cm} p{3.7cm} p{3.7cm} p{3.7cm} c}
            \toprule
            \textbf{Article} & \textbf{Game Players (Leader, Follower)} & \textbf{Leader Actions} & \textbf{Follower Actions} & \textbf{Decisions} \\
            \midrule
            \citet{wang2015effects} & Government (leader), Firm (follower) & Impose carbon emission tax & Choose shipment transportation mode and selling price & Ec, S, En \\
            \midrule
            \citet{tsao2018seaport} & Seaport (leader), Dry ports and shippers (followers) & Determine storage price & Determine service areas, prices, and delivery schedule & Ec, S, En \\
            \midrule
            \citet{xu2018modelling} & Port operators (both) & Strategy to locate dry ports & Strategy to locate dry ports & Ec, En \\
            \midrule
            \citet{chen2023subsidy} & Government (leader), Carriers and shippers (followers) & Conduct pricing game & Select transportation mode & Ec, En \\
            \midrule
            \citet{shams2023game} & Government (leader), FT systems (followers) & Set emission reduction coefficient and penalties & Bargain for trading permits, determine equilibrium prices, penalize cap violation & Ec, S, En \\
            \midrule
            \citet{rahiminia2023adopting} & Railway operator (leader), Shipper (follower) & Set rail transportation price & Decide on freight volume to be transported & Ec, S, En \\
            \midrule
            \citet{wu2023evaluation} & Government (leader), Shippers (followers) & Determine number, location, and capacity of dry ports & Choose seaports and paths for container export & Ec, En \\
            \bottomrule
        \end{tabularx}
    }

    \begin{threeparttable}
        \begin{tablenotes}
            \footnotesize
            \item \textbf{Ec}: Economic, \textbf{S}: Social, \textbf{En}: Environmental
        \end{tablenotes}
    \end{threeparttable}
\end{table}

\begin{table}[!htbp]
    \centering
    \caption{Summary of Papers Using Bi-level Programming }
    \label{tab:Bi}
    \footnotesize

    \resizebox{\textwidth}{!}{
        \begin{tabularx}{\textwidth}{p{0.75in} c p{1.1in} p{1.1in} p{1.1in} p{0.5in} p{0.65in}}
            \toprule
           \textbf{Article}  & \textbf{Problem} & \textbf{Upper-level Aim} & \textbf{Lower-level Aim} & \textbf{Level Interactions} &  \textbf{Solution Technique} & \textbf{Decision Level}
           \\ 
            \midrule
            \cite{kim2013multimodal}   & NDP    & Improve network by assessing facility capacity and mobility & Determine traffic volume of each facility & LL considers facility capacity set by UL & MH & T \\
            \midrule
            \cite{zhang2013optimization}   & NDP & Best terminal network setup and emission price & Assign multi-commodity flow in the network & LL follows UL to optimize cost and emission & H, MH & S, T \\
            \midrule
            \cite{chen2014design}   & PS    & Select optimal liner route to minimize state subsidies & Assign traffic flow, estimate emissions and profit using UEA model & LL returns profit to UL & MH & O \\
            \midrule
            \cite{zhang2015freight}   & SNDP    & Minimize total costs and  CO$_2 $ emissions & Assign multi-commodity flow & LL returns total link costs and flows to UL & H, MH & S, T \\
            \midrule
            \cite{lin2017modeling}   & NDP  & Choose projects for railway network planning & Assign railway freight flow to maximize profit & UL assigns freight flow via LL model to compute objective & E & S \\
            \midrule
            \cite{zhang2018joint}   & PS   & Select infrastructure investments and subsidies & Describe selected service routes of logistics users using UEA model & LL returns carbon emission to UL for benefit-cost ratio & MH & S, T \\
            \midrule
            \cite{li2020integrated}   & SNDP    & Subsidize rail shift while maximizing profit & Minimize total cost & LL returns equilibrium freight volume to UL & H & S, T, O \\
            \midrule
            \cite{jiang2020regional}   & PS    & Maximize total routing flow & Minimize transportation cost using flow assignment & LL returns assigned flows to UL & E & S, T \\
            \midrule
            \cite{wang2021bi}   & SNDP   & Maximize revenue while minimizing emissions & Assign cargo flow using UEA model & LL returns freight volume and price to UL & H & T \\
            \midrule
            \cite{sun2023study}   & Mixed   & Optimize channel upgrades, liner size, and frequency & Determine optimal container routes using discrete choice theory & LL output feeds UL. Total cost from both LL and UL & H & S, T \\
            \bottomrule
        \end{tabularx}
    }

    \begin{threeparttable}
        \begin{tablenotes}
            \footnotesize
            \item \textbf{LL}: Lower level, \textbf{UL}: Upper level, 
            \textbf{E}: Exact,
 \textbf{H}: Heuristics, \textbf{MH}: Meta-heuristics,
            \textbf{S}: Strategic, \textbf{O}: Operational, \textbf{T}: Tactical, \textbf{UEA}: User equilibrium assignment
        \end{tablenotes}
    \end{threeparttable}
\end{table}

\noindent  \textbf{Bi-level Programming:} The number of studies utilizing a bi-level programming approach has declined recently, with only a 20\% increase in the last two years. This trend reflects a diminishing interest in applying this method. The detailed overview of studies using bi-level programming is presented in Table \ref{tab:Bi}.The table is structured into columns detailing key aspects of each study, including problem classification, the aim of upper-level (UL) and lower-level (LL) analyses, interactions between levels, techniques used, and the decision level addressed.  This approach is particularly beneficial for IMT systems for its ability to address their hierarchical nature. Despite this, we believe that the disinterest in bi-level programming is likely due to a shift toward more flexible modeling approaches, such as simulation, and the emergence of new methodologies, like machine learning, which offer simpler and more effective solutions. From the table, we observe that 8 out of 10 studies employ heuristics or metaheuristics as solution techniques, while only 2 utilize exact methods. Additionally, decision-making at the tactical level is the most explored, represented in 8 studies, whereas only 2 analyze strategic decisions. In the earlier years, NDPs were mainly formulated as bi-level programming models, but recent years have seen a shifted focus on SNDP.

\noindent  \textbf{Stochasticity in IMT:} Approximately 20\% of studies have addressed uncertainty in their analysis, with around 50\% conducted between 2020 and 2023, reflecting a growing interest in stochastic models (Table \ref{tab:Uncertainty}). The overview of research on uncertainty in IMT problems, including various sources of uncertainty and methods for addressing them. Demand uncertainty is the most commonly studied issue, appearing in about 90\% of the research, likely due to its significant impact on the IMT supply chain. Transportation time uncertainty is also a significant focus, ranked second to demand, due to its influence on scheduling and coordination. In contrast, factors such as fixed costs, supply, speed risk, and transshipment costs are less frequently examined, suggesting they are deemed less critical. Most studies focus on single or dual uncertainties, with only 4 out of 26 studies addressing more than two simultaneously. Regarding decision levels, uncertainties are predominantly analyzed at the tactical level, as seen in 19 out of 26 studies. This trend has remained consistent over the years, whether individually or in combination. Notably, transportation time uncertainties have not been explored at the strategic level, and uncertainties related to transportation emissions and carbon pricing remain unaddressed at the operational level. In terms of uncertainty handling techniques, around 40\% of the studies used chance-constrained methods and fuzzy theory across all decision levels. Four studies employed scenario generation for diverse outcomes, while only one quantified scenario likelihood with probability distributions. Additionally, only one study addressed stochasticity across all levels, highlighting a need for further research.




\begin{table}[!htbp]
    \centering
    \caption{Summary of Papers Addressing Uncertainty}
    \label{tab:Uncertainty}
    
    \footnotesize
    \resizebox{\textwidth}{!}{
        \begin{tabularx}{\textwidth}{p{3cm} p{1cm} p{0.15cm} p{0.32cm} p{0.20cm} p{0.32cm} p{0.20cm} p{0.32cm} p{0.15cm} p{0.20cm} p{0.10cm} p{0.20cm} p{0.15cm} p{0.10cm} p{0.15cm} p{3.5cm}}
        \toprule
          

        \textbf{Article}
        & 
        \textbf{Decision Level}
        & \multicolumn{13}{c}{\textbf{Uncertainty Source}} &  \textbf{Uncertainty Handling Technique} \\ \cline{3-15}
          
        && \textbf{CA} & \textbf{TRT} & \textbf{FC} & \textbf{TRC} & \textbf{PC} & \textbf{TRE} &  \textbf{TT} & \textbf{CP} & \textbf{D} & \textbf{SU} & \textbf{ET} & \textbf{R} & \textbf{AT} & \\

        \midrule
        \citet{holmgren2012tapas} &   S   &   & & & &&&&  & \checkmark &  & && & ABM \\
        \midrule
        \citet{pishvaee2012credibility}   & S    & \checkmark  & & \checkmark & \checkmark & \checkmark & \checkmark &&  & \checkmark & && && Fuzzy CCP \\
        \midrule
        \citet{demir2016green}    & T,O &     & \checkmark & & &&&&  & \checkmark & &&&& SAA \\
        \midrule
        \citet{rezaee2017green}    & S,T   &   & & & &&&&  \checkmark & \checkmark & && && Scenario generation  \\       
        \midrule
        \citet{sun2018time}  &   O  & \checkmark  & \checkmark & & &&&&  & &&&&& Fuzzy CCP \\
        \midrule
        \citet{zhao2018stochastic}  & T,O   &   & & & &&&&  & \checkmark & \checkmark & &&& CCP, SAA \\
        \midrule
        \citet{layeb2018simulation}  & T   &   & \checkmark & & &&&&  & \checkmark & & &&& SBO  \\
        \midrule
        \citet{hruvsovsky2018hybrid}      & T &  & \checkmark & & &&&&  && &&&& ABM, DES \\
        \midrule
        \citet{haddadsisakht2018closed}      & S,T &   & & & &&&& \checkmark & \checkmark & &&&& Probabilistic scenario generation, RO + SP \\
        \midrule
        \citet{tsao2019multi} & S,T     &   & & & && \checkmark & &  & \checkmark & &  &&& RO + fuzzy programming \\
        \midrule
        \citet{sun2019bi}   & S, T &   & & & &&&&  & &&& \checkmark && Fuzzy credibilistic CCP \\
        \midrule
        \citet{baykasouglu2019fuzzy}      & S,T,O & \checkmark  & & & \checkmark && & \checkmark &  & \checkmark &&& && Hybrid CCP + fuzzy interactive resolution based approach  \\
        \midrule
        \citet{sun2020green} & O     & \checkmark  & & &  &&& \checkmark &  &  &&& && Fuzzy set theory, Chance constraints  \\
        \midrule
        \citet{wang2020transportation}     & T & \checkmark  & & &  &&& &  &  &&& && Fuzzy set theory, Fuzzy constraints  \\
        \midrule
        \citet{wang2020integrated} & O     &   & & &  &&& &  &  &&& & \checkmark & Probability distribution \\
        \midrule
        \citet{jiang2020regional} & S,T     &   & & &  &&& &  & \checkmark &&& && RO  \\
        \midrule
        \citet{dai2020distributionally} & T,O  &   & & &  &&& &  & \checkmark &&& && Distributionally robust CCP + approximation technique \\
        \midrule
        \citet{mousavi2020robust}       & S, T &   & & \checkmark & \checkmark & \checkmark && &  & \checkmark &&& && Robust stochastic possibilistic CCP, Fuzzy theory  \\
        \midrule
        \citet{zhang2021low}   & T   &   & \checkmark & &  &&& &  & \checkmark &&& && RO + SP  \\
        \midrule
        \citet{li2022path} &  T    &   & & &  &&& & \checkmark  & \checkmark &&& && RO + SP  \\
        \midrule
        \citet{sun2022green}  & T    &   & \checkmark & &  && \checkmark & &  &  &&& && Interactive fuzzy CCP \\
        \midrule
        \citet{ko2022stochastic} &  S    &   & & &  \checkmark &&& &  &  &&& &&  Scenario generation \\    
        \midrule
        \citet{heinold2023primal} & O     &   & \checkmark & &  &&& &  & \checkmark  &&& && SBO, SP \\  
        \midrule
        \citet{li2023optimum}   & T   &   & \checkmark & &  &&& \checkmark &  & \checkmark &&& && Fuzzy expected value method, fuzzy CCP  \\ 
        \midrule
        \citet{sun2023modeling}   & T   &  \checkmark & & &  &&& &  &  && & &&  Fuzzy CCP \\ 
        \midrule
        \citet{guo2023toward}   & T   &  & & &  &&& &  &  && \checkmark & &&  Scenario generation \\ 
        \bottomrule
        \end{tabularx}
    }
    
    \begin{tablenotes}
        \footnotesize
        \item\textbf{S}:Strategic, \textbf{O}:Operational, \textbf{T}:Tactical,
        \textbf{CA}:Capacity, \textbf{TRT}:Transportation Time, \textbf{FC}:Fixed Cost, \textbf{TRC}:Transportation Cost,
        \textbf{PC}:Production Cost, \textbf{TRE}:Transportation Emission, 
        \textbf{TT}:Transshipment Time, \textbf{CP}:Carbon Price/Tax,  \textbf{D}:Demand, \textbf{SU}:Supply, \textbf{ET}:Economic Transition, \textbf{R}:Risk, \textbf{AT}:Arrival Time,              
        \textbf{ABM}:Agent Based Modelling,
        \textbf{CCP}:Chance Constrained Programming,
        \textbf{SAA}:Sample Average Approximation,
        \textbf{RO}:Robust Optimization,
        \textbf{DES}:Discrete Event Simulation,
        \textbf{SBO}:Simulation-Based Optimization,         
        \textbf{SP}:Stochastic Programming
    \end{tablenotes}
\end{table}

\section{Future Research Directions \& Identified Gaps} \label{FutureResearch}

This section identifies research gaps, proposes potential avenues, and highlights key areas for  exploration.



\noindent \textbf{Network reliability:} Time is a critical factor for both MMT and IMT system reliability \citep{zhao2020study}, which enables smooth coordination of multiple transport services with varying schedules arranged in sequence, directly influencing the overall performance of the system \citep{tyroller2021challenges}. However, the current literature mostly considers the transportation time and uncertainties related to that and neglects transshipment time or handling time, making the state-of-the-art far away from real-world scenarios. Additionally, there is a significant lack of studies examining the interplay between multiple uncertainties, such as demand and travel times at the strategic level, as well as the relationship between demand and carbon pricing or taxes at the operational level. In this regard, investigating multistage stochastic models could enhance uncertainty management in IMT systems, enabling decision-makers to optimize operations by assessing the impact of current decisions on future stages and dynamic factors. While Markov Decision Processes (MDPs) are effective for modeling complex, multi-stage problems, their application in current literature is limited \citep{steimle2021multi}, highlighting the need for further research. Partially observable MDPs show promise for managing uncertainties by operating under incomplete information or limited visibility of the system's state \citep{ccelik2015post}.

The efficiency of terminal operations also plays a crucial role in service network reliability. This issue, commonly studied under AMP, requires further exploration in the context of IMT and MMT decarbonization, as highlighted in Table \ref{tab:ORByProblem}. Real-world terminal operations face numerous challenges, including traffic flow delays, equipment congestion, and labor shortages. To overcome these challenges, a shift towards Precision Scheduled Railroading (PSR) can be beneficial, emphasizing regular trips over maximum train lengths and underscoring the necessity for research into labor scheduling and operational efficiencies that can enhance throughput without compromising environmental goals. Explicitly modeling inventories for container routing management at terminals can also provide valuable insights into the trade-offs between emissions and inventory costs \citep{haddadsisakht2018closed}.  


\noindent \textbf{Network resiliency:} Disruptions due to hub failures, such as non-functionality of transfer nodes, can significantly disrupt the entire IMT network, leading to heightened vulnerabilities in critical infrastructure \citep{ hosseini2021freight}. Until now, these infrastructures have been analyzed to determine the interdependency of critical nodes in terms of cost or operational efficiency \citep{hossain2020modeling}. However, an unexplored research avenue still considers the carbon emissions generated through node interdependency. Further, non-operational transfer nodes create intermodal terminal congestion requiring efficient queuing systems in intermodal routing \citep{sun2023modeling}. In this regard, holistic network visualization is instrumental, which can be achieved by integrating technologies like the Internet of Things \citep{fraga2017towards}, Digital Twins, 
and network initiatives like Virtual Watch Tower \citep{VirtualTower} for real-time re-routing.

Emerging technologies like artificial intelligence (AI), reinforcement learning (RL), and blockchain could further enhance the resiliency of IMT systems. For example, RL has recently been adopted by many OR studies for solving complex problems like vehicle routing and TSP, particularly when traditional methods face limitations in scalability and adaptability \citep{hildebrandt2023opportunities}. In the field of IMT, RL’s ability to handle dynamic, multi-faceted environments and learn from complex interactions makes it a valuable tool. Additionally, greater SC visibility through improved intermodal shipment tracking, GPS systems, telematics technologies, and cloud-based intermodal software can facilitate partnerships among different companies operating in the intermodal sector.

\noindent \textbf{Pricing Mechanisms:} Since RMPs are relatively underexplored in the literature on IMT, the pricing mechanisms for transportation services addressed in this context are also limited. Exploring this gap while considering collaboration among stakeholders can influence economies of scale, providing an interesting avenue for research \citep{zhang2022preference, tao2021energy}. Another extension could involve integrating quantity discounts across different modes \citep{rahiminia2023adopting}. 
Uncertainty surrounding pricing mechanisms, in the context of multi-period dynamic investment and pricing problems \citep{zhang2018joint}, remains largely overlooked in our surveyed literature, highlighting a rich area for exploration. Addressing the influence of strategic decisions (i.e., port location) on pricing issues among governments, carriers, and shippers while considering transportation and processing capacities may yield more realistic insights \citep{wu2023game}. 
Further, exploration of factors that influence pricing, i.e., demand and supply balance, may provide more insights. In this regard, demand predictions are a challenge that can be handled through AI predictions \citep{luo2016revenue} with the integration of behavioral aspects.

\noindent \textbf{Stakeholder Collaboration:}
Collaboration among diverse stakeholders (i.e., governments, logistics providers, infrastructure developers, and policymakers) can harness the potential of IMT despite their varying and often conflicting objectives. Incorporating effective stakeholder collaboration on the optimal use of available assets into optimization models can improve resource utilization and reduce carbon emissions in IMT \citep{guo2020dynamic}. Furthermore, it is essential to investigate competition among similar stakeholders, such as port operators, to explore opportunities for sharing dry ports and thereby reduce competition and fixed costs \citep{xu2018modelling}. Developing integrated global and regional models and collaboration strategies that address trade imbalances can be another exciting research avenue \citep{castrellon2023assessing}. Another issue for smooth collaboration is the privatization of railways while addressing national transportation problems. Due to privatization, companies seeking to transport goods must negotiate individually with specific railroads, thereby restricting access to critical infrastructure.

While collaboration among multiple rail providers holds promise for enhancing efficiency, substantial challenges remain, including information sharing, equity, and real-time planning issues. Thus, further research is critical to establish standardized protocols, advance digital technologies, and identify effective government incentives to promote collaboration. Moreover, effective stakeholder collaboration facilitates the implementation of advanced digital technologies, such as digital twins and machine learning algorithms, thereby enhancing real-time planning and improving information exchange across various stages of integrated MMT. Consequently, the advancement of collaborative efforts and digital technologies must progress parallel to maximize their potential benefits.

Integrating shipper heterogeneity into IMT models in case of stakeholder collaboration is crucial, as carriers often transport containers for multiple shippers with diverse preferences, including low cost, speed, sustainability, reliability, and low risk. Future research could focus on IMT and SMT planning considering these varying preferences \citep{wu2023game, wu2023evaluation}. Additionally, exploring preference-based multi-objective optimization models presents an interesting research avenue \citep{zhang2022flexible}. It is essential to recognize that shippers' stated preferences may differ from their actual behaviors. Future studies should analyze shippers’ historical decisions to model these preferences more accurately and learn from their current choices in a dynamic context \citep{zhang2022heterogeneous}. A thorough analysis of behavior is essential to understand and replicate the decision-making processes of both shippers and carriers, and 
a cost-benefit analysis of the assessed policy measures is vital \citep{nassar2023system}.


\noindent \textbf{Alternate Energy Sources and Autonomous Vehicles:} The transition to alternative energy sources in IMT necessitates a deeper exploration of holistic estimation models, such as LCA frameworks, to evaluate transportation and transshipment emissions \citep{xu2018modelling} while considering uncertainty in emissions rates across different transport modes throughout the entire transportation process. This is particularly relevant with the rise of technologies like 100\% battery-powered locomotives and hybrid vehicles, as well as the increasing utilization of hydrogen and electric vehicles. Strategic studies are essential for understanding the required investments versus potential energy efficiency improvements, thereby enabling more precise financial planning \citep{zhou2018capacitated}. Research should also focus on developing energy infrastructure, such as charging stations, and assess the viability of existing resources to support long-haul electric trucks and trains \citep{farahani2018decision}. Moreover, future research can assess the development and integration of high-capacity 1 MW mega chargers in public locations for longer-distance transport, as well as the establishment of battery swapping stations.

Integrating autonomous vehicles, such as self-driven trucks and drones, within intermodal logistics networks can significantly enhance efficiency and accessibility while promoting environmental sustainability. Drones (UAVs) are a pivotal extension of existing ground transport systems, particularly in last-mile delivery for urban freight. Operating from intermodal hubs, these UAVs can swiftly pick up goods and deliver them to nearby areas, facilitating seamless transitions between transport modes. Additionally, drones can monitor, inspect, and maintain logistics operations, offering real-time tracking and data-gathering capabilities. However, multi-modal drone delivery's sustainability and emissions aspects have yet to be fully explored. The potential for UAVs to support cold chain logistics further exemplifies their versatility, combining last-mile de-consolidation with end-haul transportation. Furthermore, studies evaluating the trade-offs between traditional and autonomous vehicles focusing on their respective characteristics and benefits are essential.

\noindent \textbf{Container De/Consolidation and Repositioning:} Other significant yet frequently overlooked aspects in the IMT literature are freight consolidation and deconsolidation, as well as the repositioning of empty containers.  Considering consolidation and delivery links together can possibly yield more instructive solutions to improve the overall sustainability of the transportation system \citep{lu2019sustainable}. Research should also delve into the impacts of freight consolidation on smaller shipments as well as shipment packaging \citep{harris2018impact}, as this can influence emissions and overall efficiency. Developing GT models incorporating cargo containerization technology within the shipping market is another compelling area for exploration \citep{feng2023multimodal}. Another often overlooked aspect in the current literature is the role of container size. The dimensions of containers directly influence freight consolidation and deconsolidation strategies, impacting both transportation costs and carbon footprints. By neglecting container size, existing studies miss critical dynamics related to load optimization and space utilization, leading to less effective sustainability outcomes. Regarding repositioning, incorporating factors such as container shortages and transport network contingencies can significantly enhance the evaluation of ECR strategies. Strategies like trip-sharing, which combines bulk cargo transportation with ECR to reduce costs for both parties, are particularly promising in this context. 






\noindent \textbf{Underutilized OR Techniques:}
The use of some OR techniques, such as LR \citep{hoen2014switching} and BD \citep{haddadsisakht2018closed}, in IMT literature remains limited. This could be mainly due to the challenges posed by the consolidation and deconsolidation of goods, which disrupt the block structure essential for these traditional optimization algorithms and complicate large-scale linear programming formulations. Consequently, there has been a noticeable inclination in the literature toward meta-heuristics and hybrid solution techniques, primarily because these methods are more adept at managing the inherent complexities of IMT and offer enhanced flexibility. However, we argue that further exploration and assessment of OR techniques within the context of decarbonizing the IMT sector is essential, particularly given the demonstrated benefits of these methods. This includes column generation for handling pricing mechanisms \citep{agarwal2008ship}, BD for managing economies of scale \citep{de2009benders}, and partial BD for solving network design problems \citep{crainic2021partial}.



A persistent challenge in optimization is developing effective solution methods for large problem instances. While exact solution techniques can be applied to small theoretical cases in IMT problems, they often fall short for larger, more complex scenarios. Various meta-heuristic approaches have been introduced and are typically faster, but it's unclear which algorithms best suit specific real-world problems, as the literature lacks clarity on efficiency. The absence of standardized reference problems complicates comparisons, with most studies focusing on unique variants. Establishing benchmark datasets could improve evaluation, but reliance on proprietary datasets hinders meaningful comparisons. 
Another challenge is the computational time and scalability associated with complex problems. To tackle this issue, researchers have increasingly turned to techniques like Neural Combinatorial Optimization, which have emerged as a recent trend \citep{deineko2024learning,aghazadeh2023reinforcement} to enhance network reliability and resiliency. However, their application is scarcely explored in the decarbonization of IMT.

\section{Conclusion} \label{Conclusions}
In this study, we thoroughly review studies from the OR field aimed at contributing to the decarbonization of the freight sector taking advantage of the benefits displayed by IMT, MMT, and SMT. We begin by describing the evolution of the field in terms of relevant topics and OR techniques used to studied them. The first attempts in the field were made to recognize the benefits of IMT and MMT over other transportation systems, followed by enhancements in network design, estimation models, and routing accomplished by improvements on OR tools. 

We recognize an increasing trend in the number of publications in the last few years, resembling the pressing need to reduce emissions from the sector where scholars, practitioners and government making use of OR tools that simultaneously evaluate the impact of multiple strategies or carbon footprint policies. We also provide an integrated analysis of the literature to provide insights about aspects in the field such as modality mix, and planning decision level. Finally, we identify research gaps that could give direction to future research in the field.

\section*{Declaration of interests}
The authors have no competing interests to declare that are relevant to the content of this article.
\newpage

\appendix
\renewcommand{\thefigure}{\Alph{section}.\arabic{figure}}
\renewcommand{\thetable}{\Alph{section}.\arabic{table}}
\setcounter{figure}{0}
\setcounter{table}{0}
\renewcommand{\thesection}{Appendix \Alph{section}} 
\section{Review Methodology} \label{Review_Methodology}

\begin{figure}[!htp]
    \centering
    
    \begin{tikzpicture}[
        font=\small, 
        node distance=0.25cm and 0.2cm, 
        block/.style={
            arrow box, 
            draw, 
            text width=1.75cm, 
            text centered, 
            fill=white!15, 
            minimum height=1em,
            arrow box arrows={east:0.5cm},
            drop shadow
        },
        line/.style={draw, -Latex, thick},
        last block/.style={
            block,
            arrow box arrows={east:0cm}
        },    scale=1, 
        every node/.style={transform shape} 
    ]    

    \node [block] (step1) {Review \\ questions};
    \node [block, right=of step1] (step2) {Literature \\collection };
    \node [block, right=of step2] (step3) { Filter $\&$ \\snowball};
    \node [block, right=of step3] (step4) {Classify  \\ literature};
    \node [block, right=of step4] (step5) {Analysis  $\&$ synthesis};
    \node [last block, right=of step5] (step6) {Future directions};
    \end{tikzpicture}
    \caption{Systematic Literature Review Process Flowchart}
    \label{fig:review-process}
\end{figure}
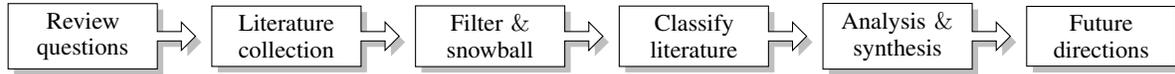

In this section, we discuss the details of the methodology we employ to identify studies related to IMT and decarbonization.
We adapt a systematic literature review (SLR) methodology as illustrated in Fig.~\ref{fig:review-process}, previously utilized in domain-specific studies (see \citet{agamez2017intermodal}).

First, in response to the gaps identified in the  literature, we develop three key research questions. 
\begin{enumerate}
    \item What mitigation strategies have been studied to reduce emissions in the IMT sector?
    \item Which OR techniques, along with their characterization, have been applied in IMT decarbonization?
    \item What are the promising further research directions in IMT decarbonization?
\end{enumerate}

Second, we define our search strategy and databases. We establish the search string as shown below by combining relevant terms aligned with our research questions \citep{meyer2020decafreighttrans}: \textit{``intermodal transportation"} OR \textit{``multimodal transportation"} OR \textit{``freight transportation"} AND (\textit{green} OR \textit{decarboni}* OR \textit{``environmental sustainab*”} OR ((\textit{fuel} OR \textit{“carbon emission”} OR \textit{“carbon dioxide emission”} OR \textit{emission} OR \textit{``greenhouse gas”} OR \textit{ghg})) AND (\textit{reduc*} OR \textit{improve*} OR \textit{efficiency})).

Our research primarily utilizes \textit{Web of Science}, \textit{Scopus}, \textit{IEEE}, and \textit{ABI Inform (ProQuest)} databases, chosen due to their extensive coverage in our field \citep{emodi2022transportdeca}. We focus our search on the titles, keywords, and abstracts of English-language publications from 2010 to 2024. This period aligns with the rising importance of modeling techniques in IMT systems and the growing emphasis on freight decarbonization to mitigate environmental impacts.

Third, we establish inclusion criteria to filter publications that use an OR technique, involve at least two modes of transport, and incorporate a decarbonization strategy.
We identify a total of 89 publications meeting these criteria. We then apply a backward snowballing technique, as described by \cite{wnuk2018knowledge}, using references from these publications to find additional relevant papers.  This process increases our total number of publications to a combined total of 133. Finally, the literature is analyzed according to the chronological evolution of the domain, types of modes involved, sustainability strategies, and the OR technique applied.

\newpage
\section{Timeline} \label{Timeline}

\renewcommand\footnoterule

\begin{tikzpicture}
  \draw (-1,0) -- (13,0) -- (13, -4) -- (-3,-4) -- (-3, -8) -- (12,-8); 

  \foreach \x in {-1,2,5,8,11} {
    \draw (\x cm,3pt) -- (\x cm,-3pt);
    \draw (\x cm,-3.9) -- (\x cm,-4.1); 
    \draw (\x cm,-7.9) -- (\x cm,-8.1); 
  }

  \draw (-1,0) node[below=3pt] {2010} 
        node[above=3pt, text width=3.5cm, align=center]  {\scriptsize • Emission costs$^{1}$ \\ 
         • Mode-specific emission costs$^{2}$\\
         • Hub and spoke$^{3}$\\
         };
  
  \draw (2,0) node[below=3pt, text width=3.5cm, align=center] {\scriptsize  
            • Dry ports' impact$^{4}$ \\ 
            • Carbon cap$^{5}$ \\
            • Cap and trade$^{6}$\\ }
        node[above=3pt] {2011};

  \draw (5,0) node[below=3pt] {2012} 
        node[above=3pt, text width=3.5cm, align=center] {\scriptsize  
            • Uncertainty in Supply Chain Network Design$^{7,8}$ \\
            • Agent based model$^{8}$ \\
            • Social impacts$^{9}$\\  };
  
  \draw (8,0) node[below=3pt, text width=3.5cm, align=center] {\scriptsize 
            • Vehicle-specific emission cost$^{10}$ \\
            • Discrete event simulation$^{11}$ \\}
        node[above=3pt] {2013};

  \draw (11,0) node[below=3pt] {2014} node[above=3pt, text width=3.5cm, align=center] {\scriptsize  
  • Decision support system$^{12}$ \\
  • Bi level program$^{13}$ \\
  •Optimization+simulation$^{14}$ \\
  };

  \draw (-1,-4) node[below=3pt] {2019} node[above=3pt, text width=3.5cm, align=center] {\scriptsize  • System dynamics$^{28}$ \\
  • Fuzzy+stochastic$^{29}$ \\
  • Fuzzy+robust$^{30}$ \\
  • Real-time planning$^{31}$ \\
  };
  \draw (2,-4) node[below=3pt, text width=3.5cm, align=center] {\scriptsize  • Stochastic+robust$^{26}$\\
  • Port competition$^{27}$\\} node[above=3pt] {2018};
  \draw (5,-4) node[below=3pt] {2017} node[above=3pt, text width=3.5cm, align=center] {\scriptsize  • Uncertainty in carbon price  $^{25}$\\ };
  \draw (8,-4) node[below=3pt, text width=3.5cm, align=center] {\scriptsize 
  • Uncertainty in Service Network Design$^{20}$\\
  • Empty container repositioning$^{21}$\\
  • Transshipment$^{22}$ \\
  • Synchromodality$^{23}$\\
  • Consolidation$^{24}$\\
  } node[above=3pt] {2016};
  \draw (11,-4) node[below=3pt] {2015} node[above=3pt, text width=3.5cm, align=center] {\scriptsize  • Container rerouting$^{15}$ \\
  • Governmental tax$^{16}$ \\
  • Governmental subsidies$^{17}$\\
  • Game theory$^{18}$\\
  •	Mesoscopic emission estimation model$^{19}$\\};

  \draw (-1,-8) node[below=3pt] {2020} node[above=3pt, text width=3.5cm, align=center] {\scriptsize • Shipment matching$^{32}$ \\
  • Supplier risk$^{33}$ \\
  • Eco-labels $^{34}$ \\};
  
  \draw (2,-8) node[below=3pt, text width=3.5cm, align=center] {\scriptsize  • Fourth-party logistic providers$^{35}$ \\
  • Offline planning+online replanning$^{36}$ \\
  • Transshipment emission in estimation models$^{37}$ \\
  } node[above=3pt] {2021};
  \draw (5,-8) node[below=3pt] {2022} node[above=3pt, text width=3.5cm, align=center] {\scriptsize  • Electric mode emissions$^{38}$ \\
  • Shipper preference$^{39}$ \\
  • Rise in synchromodal research$^{39,40,41}$ \\
 };
  \draw (8,-8) node[below=3pt, text width=3.5cm, align=center] {\scriptsize  • Blockchain$^{42}$ 
  • Carbon peak$^{43}$\\
  • Rise in Game theory modelling$^{44,45,46}$\\
  
  } node[above=3pt] {2023};
  \draw (11,-8) node[below=3pt] {2024} node[above=3pt, text width=3.5cm, align=center] {\scriptsize  • Impact of electrification and biofuel$^{47}$ \\
  • Regional demand$^{48}$ \\};

  \draw[->] (11, -8) -- (12, -8);
\end{tikzpicture}
\vspace{0.1 cm}

\begin{justify}

\scriptsize 
1\textsubscript{\text{\cite{bauer2010minimizing}},} 
2\textsubscript{\text{\cite{chang2010optimization}},} 
3\textsubscript{\text{\cite{comer2010marine}},} 
4\textsubscript{\text{\cite{henttu2011financial}},} 
5\textsubscript{\text{\cite{zhang2011mode}},}
6\textsubscript{\text{\cite{chaabane2011designing}},}
7\textsubscript{\text{\cite{pishvaee2012credibility}},}
8\textsubscript{\text{\cite{holmgren2012tapas}},}
9\textsubscript{\text{\cite{sawadogo2012reducing}},}
10\textsubscript{\text{\cite{le2013model}},}
11\textsubscript{\text{\cite{lattila2013hinterland}},}
12\textsubscript{\text{\cite{kengpol2014development}},}
13\textsubscript{\text{\cite{chen2014design}},}
14\textsubscript{\text{\cite{liotta2014optimization}},}
15\textsubscript{\text{\cite{rodrigues2015uk}},}
16\textsubscript{\text{\cite{duan2015carbon}},}
17\textsubscript{\text{\cite{bouchery2015cost}},}
18\textsubscript{\text{\cite{wang2015effects}},}
19\textsubscript{\text{\cite{kirschstein2015ghg}},}
20\textsubscript{\text{\cite{demir2016green}},}
21\textsubscript{\text{\cite{lam2016market}},}
22\textsubscript{\text{\cite{rudi2016freight}},}
23\textsubscript{\text{\cite{zhang2016synchromodal}},}
24\textsubscript{\text{\cite{baykasouglu2016multi}},}
25\textsubscript{\text{\cite{rezaee2017green}},}
26\textsubscript{\text{\cite{haddadsisakht2018closed}},}
27\textsubscript{\text{\cite{xu2018modelling}},}
28\textsubscript{\text{\cite{choi2019system}},}
29\textsubscript{\text{\cite{baykasouglu2019fuzzy}},}
30\textsubscript{\text{\cite{tsao2018seaport}},}
31\textsubscript{\text{\cite{resat2019discrete}},}
32\textsubscript{\text{\cite{guo2020dynamic}},}
33\textsubscript{\text{\cite{kabadurmus2020sustainable}},}
34\textsubscript{\text{\cite{heinold2020emission}},}
35\textsubscript{\text{\cite{wu2021research}},}
36\textsubscript{\text{\cite{hruvsovsky2021real}},}
37\textsubscript{\text{\cite{tao2021energy}},}
38\textsubscript{\text{\cite{guo2022modeling}},}
39\textsubscript{\text{\cite{zhang2022preference}},}
40\textsubscript{\text{\cite{zhang2022flexible}},}
41\textsubscript{\text{\cite{zhang2022collaborative}},}
42\textsubscript{\text{\cite{oudani2023combined}},}
43\textsubscript{\text{\cite{zuo2023using}},}
44\textsubscript{\text{\cite{rahiminia2023adopting}},}
45\textsubscript{\text{\cite{wu2023game}},}
46\textsubscript{\text{\cite{shams2023game}},}
47\textsubscript{\text{\cite{ghisolfi2024dynamics}},}
48\textsubscript{\text{\cite{derpich2024pursuing}}}
\end{justify}

\vspace{1em}


\begin{thebibliography}{}

\bibitem[Abu~Aisha et~al., 2020]{abu2020optimization}
Abu~Aisha, T., Ouhimmou, M., and Paquet, M. (2020).
\newblock {Optimization of container terminal layouts in the seaport—Case of port of Montreal}.
\newblock {\em Sustainability}, 12(3):1165.

\bibitem[Agalianos et~al., 2020]{agalianos2020discrete}
Agalianos, K., Ponis, S., Aretoulaki, E., Plakas, G., and Efthymiou, O. (2020).
\newblock {Discrete event simulation and digital twins: review and challenges for logistics}.
\newblock {\em Procedia Manufacturing}, 51:1636--1641.

\bibitem[Agamez-Arias and Moyano-Fuentes, 2017]{agamez2017intermodal}
Agamez-Arias, A.-d.-M. and Moyano-Fuentes, J. (2017).
\newblock {Intermodal transport in freight distribution: A literature review}.
\newblock {\em Transport Reviews}, 37(6):782--807.

\bibitem[Agarwal and Ergun, 2008]{agarwal2008ship}
Agarwal, R. and Ergun, {\"O}. (2008).
\newblock {Ship scheduling and network design for cargo routing in liner shipping}.
\newblock {\em Transportation Science}, 42(2):175--196.

\bibitem[Aghazadeh and Wang, 2023]{aghazadeh2023reinforcement}
Aghazadeh, H. and Wang, X. (2023).
\newblock {Reinforcement Learning for Intermodal Transportation Planning with Time Windows and Limited Cargo Capacity}.
\newblock In {\em Proceedings of the 16th ACM SIGSPATIAL International Workshop on Computational Transportation Science}, pages 28--31.

\bibitem[Amaran et~al., 2016]{amaran2016simulation}
Amaran, S., Sahinidis, N.~V., Sharda, B., and Bury, S.~J. (2016).
\newblock {Simulation optimization: a review of algorithms and applications}.
\newblock {\em Annals of Operations Research}, 240:351--380.

\bibitem[Aminzadegan et~al., 2022]{Aminzadegan2022decarintermodal}
Aminzadegan, S., Shahriari, M., Mehranfar, F., and Abramović, B. (2022).
\newblock {Factors affecting the emission of pollutants in different types of transportation: A literature review}.
\newblock {\em Energy Reports}, 8:2508--2529.

\bibitem[Archetti et~al., 2022]{Archetti20221ORMultimodal}
Archetti, C., Peirano, L., and Speranza, M.~G. (2022).
\newblock {Optimization in multimodal freight transportation problems: A Survey}.
\newblock {\em European Journal of Operational Research}, 299(1):1--20.

\bibitem[ARPA-E, 2023]{ARPAE2023}
ARPA-E (2023).
\newblock { Intermodal Freight Transportation System}.
\newblock \url{https://www.arpa-e.energy.gov/technologies/exploratory-topics/intermodal-freight)}.
\newblock (Accessed on 04/29/2024).

\bibitem[Banister et~al., 2011]{banister2011transportation}
Banister, D., Anderton, K., Bonilla, D., Givoni, M., and Schwanen, T. (2011).
\newblock {Transportation and the environment}.
\newblock {\em Annual review of environment and resources}, 36(1):247--270.

\bibitem[Bauer et~al., 2010]{bauer2010minimizing}
Bauer, J., Bekta{\c{s}}, T., and Crainic, T.~G. (2010).
\newblock {Minimizing greenhouse gas emissions in intermodal freight transport: an application to rail service design}.
\newblock {\em Journal of the Operational Research Society}, 61:530--542.

\bibitem[Baykaso{\u{g}}lu and Subulan, 2016]{baykasouglu2016multi}
Baykaso{\u{g}}lu, A. and Subulan, K. (2016).
\newblock {A multi-objective sustainable load planning model for intermodal transportation networks with a real-life application}.
\newblock {\em Transportation Research Part E: Logistics and Transportation Review}, 95:207--247.

\bibitem[Baykaso{\u{g}}lu and Subulan, 2019]{baykasouglu2019fuzzy}
Baykaso{\u{g}}lu, A. and Subulan, K. (2019).
\newblock {A fuzzy-stochastic optimization model for the intermodal fleet management problem of an international transportation company}.
\newblock {\em Transportation Planning and Technology}, 42(8):777--824.

\bibitem[Bouchery and Fransoo, 2015]{bouchery2015cost}
Bouchery, Y. and Fransoo, J. (2015).
\newblock {Cost, carbon emissions and modal shift in intermodal network design decisions}.
\newblock {\em International Journal of Production Economics}, 164:388--399.

\bibitem[Boussa{\"\i}d et~al., 2013]{boussaid2013survey}
Boussa{\"\i}d, I., Lepagnot, J., and Siarry, P. (2013).
\newblock {A survey on optimization metaheuristics}.
\newblock {\em Information sciences}, 237:82--117.

\bibitem[Cachon and Netessine, 2006]{cachon2006game}
Cachon, G.~P. and Netessine, S. (2006).
\newblock {Game theory in supply chain analysis}.
\newblock {\em Models, methods, and applications for innovative decision making}, pages 200--233.

\bibitem[Caris et~al., 2013]{caris2013decision}
Caris, A., Macharis, C., and Janssens, G.~K. (2013).
\newblock {Decision support in intermodal transport: A new research agenda}.
\newblock {\em Computers in industry}, 64(2):105--112.

\bibitem[Castrellon et~al., 2023]{castrellon2023assessing}
Castrellon, J.~P., Sanchez-Diaz, I., Roso, V., Altuntas-Vural, C., Rogerson, S., Sant{\'e}n, V., and Kalahasthi, L.~K. (2023).
\newblock {Assessing the eco-efficiency benefits of empty container repositioning strategies via dry ports}.
\newblock {\em Transportation Research Part D: Transport and Environment}, 120:103778.

\bibitem[{\c{C}}elik et~al., 2015]{ccelik2015post}
{\c{C}}elik, M., Ergun, {\"O}., and Keskinocak, P. (2015).
\newblock {The post-disaster debris clearance problem under incomplete information}.
\newblock {\em Operations Research}, 63(1):65--85.

\bibitem[Chaabane et~al., 2011]{chaabane2011designing}
Chaabane, A., Ramudhin, A., and Paquet, M. (2011).
\newblock {Designing supply chains with sustainability considerations}.
\newblock {\em Production Planning \& Control}, 22(8):727--741.

\bibitem[Chaabane et~al., 2012]{chaabane2012design}
Chaabane, A., Ramudhin, A., and Paquet, M. (2012).
\newblock {Design of sustainable supply chains under the emission trading scheme}.
\newblock {\em International journal of production economics}, 135(1):37--49.

\bibitem[Chang et~al., 2010]{chang2010optimization}
Chang, Y.-T., Lee, P. T.-W., Kim, H.-J., and Shin, S.-H. (2010).
\newblock {Optimization model for transportation of container cargoes considering short sea shipping and external cost: South Korean case}.
\newblock {\em Transportation Research Record}, 2166(1):99--108.

\bibitem[Chasek et~al., 2016]{chasek2016getting}
Chasek, P.~S., Wagner, L.~M., Leone, F., Lebada, A.-M., and Risse, N. (2016).
\newblock {Getting to 2030: Negotiating the post-2015 sustainable development agenda}.
\newblock {\em Review of European, Comparative \& International Environmental Law}, 25(1):5--14.

\bibitem[Chen et~al., 2014]{chen2014design}
Chen, K., Yang, Z., and Notteboom, T. (2014).
\newblock {The design of coastal shipping services subject to carbon emission reduction targets and state subsidy levels}.
\newblock {\em Transportation Research Part E: Logistics and Transportation Review}, 61:192--211.

\bibitem[Chen et~al., 2023]{chen2023subsidy}
Chen, Z., Zhang, Z., Bian, Z., Dai, L., and Hu, H. (2023).
\newblock {Subsidy policy optimization of multimodal transport on emission reduction considering carrier pricing game and shipping resilience: A case study of Shanghai port}.
\newblock {\em Ocean \& Coastal Management}, 243:106760.

\bibitem[Choi et~al., 2019]{choi2019system}
Choi, B., Park, S.-i., and Lee, K.-D. (2019).
\newblock {A system dynamics model of the modal shift from road to rail: Containerization and imposition of taxes}.
\newblock {\em Journal of Advanced Transportation}, 2019(1):7232710.

\bibitem[Choong et~al., 2002]{choong2002empty}
Choong, S.~T., Cole, M.~H., and Kutanoglu, E. (2002).
\newblock {Empty container management for intermodal transportation networks}.
\newblock {\em Transportation Research Part E: Logistics and Transportation Review}, 38(6):423--438.

\bibitem[Comer et~al., 2010]{comer2010marine}
Comer, B., Corbett, J.~J., Hawker, J.~S., Korfmacher, K., Lee, E.~E., Prokop, C., and Winebrake, J.~J. (2010).
\newblock {Marine vessels as substitutes for heavy-duty trucks in Great Lakes freight transportation}.
\newblock {\em Journal of the Air \& Waste Management Association}, 60(7):884--890.

\bibitem[Commission, 2011]{ec2011roadmap}
Commission, E.-E. (2011).
\newblock {Roadmap to a Single European Transport Area-Towards a competitive and resource efficient transport system}.
\newblock {\em White Paper, Communication}, 144.

\bibitem[Craig et~al., 2013]{craig2013estimating}
Craig, A.~J., Blanco, E.~E., and Sheffi, Y. (2013).
\newblock {Estimating the CO2 intensity of intermodal freight transportation}.
\newblock {\em Transportation Research Part D: Transport and Environment}, 22:49--53.

\bibitem[Crainic, 2000]{crainic2000service}
Crainic, T.~G. (2000).
\newblock {Service network design in freight transportation}.
\newblock {\em European journal of operational research}, 122(2):272--288.

\bibitem[Crainic et~al., 2021]{crainic2021partial}
Crainic, T.~G., Hewitt, M., Maggioni, F., and Rei, W. (2021).
\newblock {Partial benders decomposition: general methodology and application to stochastic network design}.
\newblock {\em Transportation Science}, 55(2):414--435.

\bibitem[Crainic and Laporte, 1997]{crainic1997planning}
Crainic, T.~G. and Laporte, G. (1997).
\newblock {Planning models for freight transportation}.
\newblock {\em European journal of operational research}, 97(3):409--438.

\bibitem[Crainic et~al., 2018]{crainic2018simulation}
Crainic, T.~G., Perboli, G., and Rosano, M. (2018).
\newblock {Simulation of intermodal freight transportation systems: A taxonomy}.
\newblock {\em European Journal of Operational Research}, 270(2):401--418.

\bibitem[Dai and Yang, 2020]{dai2020distributionally}
Dai, Q. and Yang, J. (2020).
\newblock {A distributionally robust chance-constrained approach for modeling demand uncertainty in green port-hinterland transportation network optimization}.
\newblock {\em Symmetry}, 12(9):1492.

\bibitem[de~Almeida~Rodrigues et~al., 2023]{de2023flexible}
de~Almeida~Rodrigues, T., de~Miranda~Mota, C.~M., Ojiako, U., Chipulu, M., Marshall, A., and Dweiri, F. (2023).
\newblock {A flexible cost model for seaport-hinterland decisions in container shipping}.
\newblock {\em Research in Transportation Business \& Management}, 49:101016.

\bibitem[De~Camargo et~al., 2009]{de2009benders}
De~Camargo, R.~S., de~Miranda~Jr, G., and Luna, H. P.~L. (2009).
\newblock {Benders decomposition for hub location problems with economies of scale}.
\newblock {\em Transportation Science}, 43(1):86--97.

\bibitem[Deineko et~al., 2024]{deineko2024learning}
Deineko, E., Jungnickel, P., and Kehrt, C. (2024).
\newblock {Learning-Based Optimisation for Integrated Problems in Intermodal Freight Transport: Preliminaries, Strategies, and State of the Art}.
\newblock {\em Applied Sciences}, 14(19):8642.

\bibitem[Demir et~al., 2014]{demir2014review}
Demir, E., Bekta{\c{s}}, T., and Laporte, G. (2014).
\newblock {A review of recent research on green road freight transportation}.
\newblock {\em European Journal of Operational Research}, 237(3):775--793.

\bibitem[Demir et~al., 2016]{demir2016green}
Demir, E., Burgholzer, W., Hru{\v{s}}ovsk{\`y}, M., Ar{\i}kan, E., Jammernegg, W., and Van~Woensel, T. (2016).
\newblock {A green intermodal service network design problem with travel time uncertainty}.
\newblock {\em Transportation Research Part B: Methodological}, 93:789--807.

\bibitem[Demir et~al., 2019]{demir2019green}
Demir, E., Hru{\v{s}}ovsk{\`y}, M., Jammernegg, W., and Van~Woensel, T. (2019).
\newblock {Green intermodal freight transportation: bi-objective modelling and analysis}.
\newblock {\em International Journal of Production Research}, 57(19):6162--6180.

\bibitem[Derpich et~al., 2024]{derpich2024pursuing}
Derpich, I., Duran, C., Carrasco, R., Moreno, F., Fernandez-Campusano, C., and Espinosa-Leal, L. (2024).
\newblock {Pursuing Optimization Using Multimodal Transportation System: A Strategic Approach to Minimizing Costs and CO2 Emissions}.
\newblock {\em Journal of Marine Science and Engineering}, 12(6):976.

\bibitem[Duan and Heragu, 2015]{duan2015carbon}
Duan, X. and Heragu, S. (2015).
\newblock {Carbon emission tax policy in an intermodal transportation network}.
\newblock In {\em Proceedings of the IIE Annual Conference, Nashville, TN, USA}, volume~30.

\bibitem[Emodi et~al., 2022]{emodi2022transportdeca}
Emodi, N.~V., Okereke, C., Abam, F.~I., Diemuodeke, O.~E., Owebor, K., and Nnamani, U.~A. (2022).
\newblock {Transport sector decarbonisation in the Global South: A systematic literature review}.
\newblock {\em Energy Strategy Reviews}, 43:100925.

\bibitem[EPA, 2024]{epa_ghg}
EPA (2024).
\newblock {Global Greenhouse Gas Overview - U.S. Environmental Protection Agency}.
\newblock \url{https://www.epa.gov/ghgemissions/global-greenhouse-gas-overview}.
\newblock (Accessed on 07/08/2024).

\bibitem[Facanha and Horvath, 2006]{horvath2006environmental}
Facanha, C. and Horvath, A. (2006).
\newblock {Environmental Assessment of Freight Transportation in the US}.
\newblock {\em The International Journal of Life Cycle Assessment}, 11:229--239.

\bibitem[Fahimnia et~al., 2015]{fahimnia2015tactical}
Fahimnia, B., Sarkis, J., Choudhary, A., and Eshragh, A. (2015).
\newblock {Tactical supply chain planning under a carbon tax policy scheme: A case study}.
\newblock {\em International Journal of Production Economics}, 164:206--215.

\bibitem[Fahimnia et~al., 2013]{fahimnia2013impact}
Fahimnia, B., Sarkis, J., Dehghanian, F., Banihashemi, N., and Rahman, S. (2013).
\newblock {The impact of carbon pricing on a closed-loop supply chain: an Australian case study}.
\newblock {\em Journal of Cleaner Production}, 59:210--225.

\bibitem[Fallahi et~al., 2023]{fallahi2023game}
Fallahi, N., Hafezalkotob, A., Raissi, S., and Ghezavati, V. (2023).
\newblock {A game theoretic approach to sustainable freight transportation: competition between green, non-green and semi-green transportation networks under government sustainable policies}.
\newblock {\em Environment, Development and Sustainability}, pages 1--48.

\bibitem[Farahani et~al., 2018]{farahani2018decision}
Farahani, N.~Z., Noble, J.~S., Klein, C.~M., and Enayati, M. (2018).
\newblock {A decision support tool for energy efficient synchromodal supply chains}.
\newblock {\em Journal of Cleaner Production}, 186:682--702.

\bibitem[Feng et~al., 2023]{feng2023multimodal}
Feng, X., Song, R., Yin, W., Yin, X., and Zhang, R. (2023).
\newblock {Multimodal transportation network with cargo containerization technology: Advantages and challenges}.
\newblock {\em Transport Policy}, 132:128--143.

\bibitem[Fitzgerald et~al., 2011]{fitzgerald2011energy}
Fitzgerald, W.~B., Howitt, O.~J., Smith, I.~J., and Hume, A. (2011).
\newblock {Energy use of integral refrigerated containers in maritime transportation}.
\newblock {\em Energy Policy}, 39(4):1885--1896.

\bibitem[Forkenbrock, 2001]{forkenbrock2001comparison}
Forkenbrock, D.~J. (2001).
\newblock {Comparison of external costs of rail and truck freight transportation}.
\newblock {\em Transportation Research Part A: Policy and Practice}, 35(4):321--337.

\bibitem[Fraga-Lamas et~al., 2017]{fraga2017towards}
Fraga-Lamas, P., Fern{\'a}ndez-Caram{\'e}s, T.~M., and Castedo, L. (2017).
\newblock {Towards the Internet of smart trains: A review on industrial IoT-connected railways}.
\newblock {\em Sensors}, 17(6):1457.

\bibitem[Fulzele et~al., 2019]{fulzele2019model}
Fulzele, V., Shankar, R., and Choudhary, D. (2019).
\newblock {A model for the selection of transportation modes in the context of sustainable freight transportation}.
\newblock {\em Industrial Management \& Data Systems}, 119(8):1764--1784.

\bibitem[Galarraga~Gallastegui, 2002]{galarraga2002use}
Galarraga~Gallastegui, I. (2002).
\newblock {The use of eco-labels: a review of the literature}.
\newblock {\em European Environment}, 12(6):316--331.

\bibitem[Gallardo et~al., 2021]{gallardo2021sequential}
Gallardo, P., Murray, R., and Krumdieck, S. (2021).
\newblock {A sequential optimization-simulation approach for planning the transition to the low carbon freight system with case study in the North Island of New Zealand}.
\newblock {\em Energies}, 14(11):3339.

\bibitem[Geoffrion, 1972]{geoffrion1972generalized}
Geoffrion, A.~M. (1972).
\newblock {Generalized benders decomposition}.
\newblock {\em Journal of optimization theory and applications}, 10:237--260.

\bibitem[Ghisolfi et~al., 2024]{ghisolfi2024dynamics}
Ghisolfi, V., Tavasszy, L.~A., Correia, G. H. d. A.~R., Chaves, G. d. L.~D., and Ribeiro, G.~M. (2024).
\newblock {Dynamics of freight transport decarbonization: A simulation study for Brazil}.
\newblock {\em Transportation Research Part D: Transport and Environment}, 127:104020.

\bibitem[Giusti et~al., 2019]{giusti2019synchromodal}
Giusti, R., Manerba, D., Bruno, G., and Tadei, R. (2019).
\newblock {Synchromodal logistics: An overview of critical success factors, enabling technologies, and open research issues}.
\newblock {\em Transportation Research Part E: Logistics and Transportation Review}, 129:92--110.

\bibitem[Goulder and Schein, 2013]{goulder2013carbon}
Goulder, L.~H. and Schein, A.~R. (2013).
\newblock {Carbon taxes versus cap and trade: a critical review}.
\newblock {\em Climate Change Economics}, 4(03):1350010.

\bibitem[Green, 2010]{green2010private}
Green, J.~F. (2010).
\newblock {Private standards in the climate regime: the greenhouse gas protocol}.
\newblock {\em Business and Politics}, 12(3):1--37.

\bibitem[Guo et~al., 2023a]{guo2023integrated}
Guo, L., Du, J., Zheng, J., and He, N. (2023a).
\newblock {Integrated Planning of Feeder Route Selection, Schedule Design, and Fleet Allocation with Multimodal Transport Path Selection Considered}.
\newblock {\em Journal of Marine Science and Engineering}, 11(7):1445.

\bibitem[Guo et~al., 2023b]{guo2023toward}
Guo, T., Liu, P., Wang, C., Xie, J., Du, J., and Lim, M.~K. (2023b).
\newblock {Toward sustainable port-hinterland transportation: A holistic approach to design modal shift policy mixes}.
\newblock {\em Transportation Research Part A: Policy and Practice}, 174:103746.

\bibitem[Guo et~al., 2020]{guo2020dynamic}
Guo, W., Atasoy, B., van Blokland, W.~B., and Negenborn, R.~R. (2020).
\newblock {A dynamic shipment matching problem in hinterland synchromodal transportation}.
\newblock {\em Decision Support Systems}, 134:113289.

\bibitem[Guo et~al., 2022]{guo2022modeling}
Guo, X., He, J., Lan, M., Yu, H., and Yan, W. (2022).
\newblock {Modeling carbon emission estimation for hinterland-based container intermodal network}.
\newblock {\em Journal of Cleaner Production}, 378:134593.

\bibitem[Guo et~al., 2023c]{guo2023carbon}
Guo, X., He, J., Yu, H., and Liu, M. (2023c).
\newblock {Carbon peak simulation and peak pathway analysis for hub-and-spoke container intermodal network}.
\newblock {\em Transportation Research Part E: Logistics and Transportation Review}, 180:103332.

\bibitem[Haddadsisakht and Ryan, 2018]{haddadsisakht2018closed}
Haddadsisakht, A. and Ryan, S.~M. (2018).
\newblock {Closed-loop supply chain network design with multiple transportation modes under stochastic demand and uncertain carbon tax}.
\newblock {\em International journal of production economics}, 195:118--131.

\bibitem[Halim, 2023]{halim2023boosting}
Halim, R.~A. (2023).
\newblock {Boosting intermodal rail for decarbonizing freight transport on Java, Indonesia: A model-based policy impact assessment}.
\newblock {\em Research in Transportation Business \& Management}, 48:100909.

\bibitem[Harris et~al., 2018]{harris2018impact}
Harris, I., Rodrigues, V.~S., Pettit, S., Beresford, A., and Liashko, R. (2018).
\newblock {The impact of alternative routeing and packaging scenarios on carbon and sulphate emissions in international wine distribution}.
\newblock {\em Transportation Research Part D: Transport and Environment}, 58:261--279.

\bibitem[Heinold and Meisel, 2018]{heinold2018emission}
Heinold, A. and Meisel, F. (2018).
\newblock {Emission rates of intermodal rail/road and road-only transportation in Europe: A comprehensive simulation study}.
\newblock {\em Transportation Research Part D: Transport and Environment}, 65:421--437.

\bibitem[Heinold and Meisel, 2020]{heinold2020emission}
Heinold, A. and Meisel, F. (2020).
\newblock {Emission limits and emission allocation schemes in intermodal freight transportation}.
\newblock {\em Transportation Research Part E: Logistics and Transportation Review}, 141:101963.

\bibitem[Heinold et~al., 2023]{heinold2023primal}
Heinold, A., Meisel, F., and Ulmer, M.~W. (2023).
\newblock {Primal-dual value function approximation for stochastic dynamic intermodal transportation with eco-labels}.
\newblock {\em Transportation Science}, 57(6):1452--1472.

\bibitem[Henttu and Hilmola, 2011]{henttu2011financial}
Henttu, V. and Hilmola, O.-P. (2011).
\newblock {Financial and environmental impacts of hypothetical Finnish dry port structure}.
\newblock {\em Research in Transportation Economics}, 33(1):35--41.

\bibitem[Hildebrandt et~al., 2023]{hildebrandt2023opportunities}
Hildebrandt, F.~D., Thomas, B.~W., and Ulmer, M.~W. (2023).
\newblock {Opportunities for reinforcement learning in stochastic dynamic vehicle routing}.
\newblock {\em Computers \& operations research}, 150:106071.

\bibitem[Hoen et~al., 2014]{hoen2014switching}
Hoen, K.~M., Tan, T., Fransoo, J.~C., and van Houtum, G.-J. (2014).
\newblock {Switching transport modes to meet voluntary carbon emission targets}.
\newblock {\em Transportation Science}, 48(4):592--608.

\bibitem[Holmgren et~al., 2012]{holmgren2012tapas}
Holmgren, J., Davidsson, P., Persson, J.~A., and Ramstedt, L. (2012).
\newblock {TAPAS: A multi-agent-based model for simulation of transport chains}.
\newblock {\em Simulation Modelling Practice and Theory}, 23:1--18.

\bibitem[Hossain et~al., 2020]{hossain2020modeling}
Hossain, N. U.~I., El~Amrani, S., Jaradat, R., Marufuzzaman, M., Buchanan, R., Rinaudo, C., and Hamilton, M. (2020).
\newblock {Modeling and assessing interdependencies between critical infrastructures using Bayesian network: A case study of inland waterway port and surrounding supply chain network}.
\newblock {\em Reliability Engineering \& System Safety}, 198:106898.

\bibitem[Hosseini and Al~Khaled, 2021]{hosseini2021freight}
Hosseini, S. and Al~Khaled, A. (2021).
\newblock {Freight flow optimization to evaluate the criticality of intermodal surface transportation system infrastructures}.
\newblock {\em Computers \& Industrial Engineering}, 159:107522.

\bibitem[Hru{\v{s}}ovsk{\`y} et~al., 2018]{hruvsovsky2018hybrid}
Hru{\v{s}}ovsk{\`y}, M., Demir, E., Jammernegg, W., and Van~Woensel, T. (2018).
\newblock {Hybrid simulation and optimization approach for green intermodal transportation problem with travel time uncertainty}.
\newblock {\em Flexible Services and Manufacturing Journal}, 30:486--516.

\bibitem[Hru{\v{s}}ovsk{\`y} et~al., 2021]{hruvsovsky2021real}
Hru{\v{s}}ovsk{\`y}, M., Demir, E., Jammernegg, W., and Van~Woensel, T. (2021).
\newblock {Real-time disruption management approach for intermodal freight transportation}.
\newblock {\em Journal of Cleaner Production}, 280:124826.

\bibitem[Huang et~al., 2022]{huang2022overview}
Huang, J., Cui, Y., Zhang, L., Tong, W., Shi, Y., and Liu, Z. (2022).
\newblock {An Overview of Agent-Based Models for Transport Simulation and Analysis}.
\newblock {\em Journal of Advanced Transportation}, 2022(1):1252534.

\bibitem[Huang, 2016]{huang2016understanding}
Huang, Y. (2016).
\newblock {Understanding China's Belt \& Road initiative: motivation, framework and assessment}.
\newblock {\em China Economic Review}, 40:314--321.

\bibitem[Iannone, 2012a]{iannone2012model}
Iannone, F. (2012a).
\newblock {A model optimizing the port-hinterland logistics of containers: The case of the Campania region in Southern Italy}.
\newblock {\em Maritime Economics \& Logistics}, 14(1):33--72.

\bibitem[Iannone, 2012b]{iannone2012private}
Iannone, F. (2012b).
\newblock {The private and social cost efficiency of port hinterland container distribution through a regional logistics system}.
\newblock {\em Transportation Research Part A: Policy and Practice}, 46(9):1424--1448.

\bibitem[IATA, 2024]{air_35}
IATA (2024).
\newblock {Value of Air Cargo- International Air Transport Association}.
\newblock Accessed: 2024-07-07.

\bibitem[IEA, 2018]{international2018co2}
IEA, I. E.~A. (2018).
\newblock {\em {CO2 emissions from fuel combustion}}.
\newblock OECD.

\bibitem[Ji and Luo, 2017]{ji2017hybrid}
Ji, S.-f. and Luo, R.-j. (2017).
\newblock {A hybrid estimation of distribution algorithm for multi-objective multi-sourcing intermodal transportation network design problem considering carbon emissions}.
\newblock {\em Sustainability}, 9(7):1133.

\bibitem[Jiang et~al., 2020]{jiang2020regional}
Jiang, J., Zhang, D., Meng, Q., and Liu, Y. (2020).
\newblock {Regional multimodal logistics network design considering demand uncertainty and CO2 emission reduction target: A system-optimization approach}.
\newblock {\em Journal of Cleaner Production}, 248:119304.

\bibitem[Kabadurmus and Erdogan, 2020]{kabadurmus2020sustainable}
Kabadurmus, O. and Erdogan, M.~S. (2020).
\newblock {Sustainable, multimodal and reliable supply chain design}.
\newblock {\em Annals of Operations Research}, 292(1):47--70.

\bibitem[Ke et~al., 2023]{ke2023optimization}
Ke, H., Xu, G., Li, C., Gao, J., Xiao, X., Wu, X., and Yan, Q. (2023).
\newblock {Optimization of China’s freight transportation structure based on adaptive genetic algorithm under the background of carbon peak}.
\newblock {\em Environmental Science and Pollution Research}, 30(36):85087--85101.

\bibitem[Kengpol et~al., 2014]{kengpol2014development}
Kengpol, A., Tuammee, S., and Tuominen, M. (2014).
\newblock {The development of a framework for route selection in multimodal transportation}.
\newblock {\em The International Journal of Logistics Management}, 25(3):581--610.

\bibitem[Kim et~al., 2009]{kim2009trade}
Kim, N.~S., Janic, M., and Van~Wee, B. (2009).
\newblock {Trade-off between carbon dioxide emissions and logistics costs based on multiobjective optimization}.
\newblock {\em Transportation Research Record}, 2139(1):107--116.

\bibitem[Kim and Van~Wee, 2009a]{kim2009assessment}
Kim, N.~S. and Van~Wee, B. (2009a).
\newblock {Assessment of CO2 emissions for intermodal freight transport systems and truck-only system: case study of Western--Eastern Europe corridor}.
\newblock In {\em 88th Annual Meeting of the Transportation Research Board, Washington, D.C}.

\bibitem[Kim and Van~Wee, 2009b]{kim2009assessmentprevious}
Kim, N.~S. and Van~Wee, B. (2009b).
\newblock {Assessment of CO2 emissions for truck-only and rail-based intermodal freight systems in Europe}.
\newblock {\em Transportation Planning and Technology}, 32(4):313--333.

\bibitem[Kim et~al., 2013]{kim2013multimodal}
Kim, S., Park, M., and Lee, C. (2013).
\newblock {Multimodal freight transportation network design problem for reduction of greenhouse gas emissions}.
\newblock {\em Transportation Research Record}, 2340(1):74--83.

\bibitem[Kirschstein and Meisel, 2015]{kirschstein2015ghg}
Kirschstein, T. and Meisel, F. (2015).
\newblock {GHG-emission models for assessing the eco-friendliness of road and rail freight transports}.
\newblock {\em Transportation Research Part B: Methodological}, 73:13--33.

\bibitem[Ko et~al., 2022]{ko2022stochastic}
Ko, S., Choi, K., Yu, S., and Lee, J. (2022).
\newblock {A stochastic optimization model for sustainable multimodal transportation for bioenergy production}.
\newblock {\em Sustainability}, 14(3):1889.

\bibitem[Kurtulu{\c{s}}, 2023]{kurtulucs2023green}
Kurtulu{\c{s}}, E. (2023).
\newblock {Green Container Shipping Schedule Design Considering Collaboration and Disruption Recovery}.
\newblock {\em IEEE Access}.

\bibitem[Lam and Gu, 2016]{lam2016market}
Lam, J. S.~L. and Gu, Y. (2016).
\newblock {A market-oriented approach for intermodal network optimisation meeting cost, time and environmental requirements}.
\newblock {\em International Journal of Production Economics}, 171:266--274.

\bibitem[L{\"a}ttil{\"a} et~al., 2013]{lattila2013hinterland}
L{\"a}ttil{\"a}, L., Henttu, V., and Hilmola, O.-P. (2013).
\newblock {Hinterland operations of sea ports do matter: Dry port usage effects on transportation costs and CO2 emissions}.
\newblock {\em Transportation Research Part E: Logistics and Transportation Review}, 55:23--42.

\bibitem[Laurent et~al., 2020]{laurent2020carbonroadmap}
Laurent, A.-B., Vallerand, S., Van~der Meer, Y., and D'Amours, S. (2020).
\newblock {CarbonRoadMap: A multicriteria decision tool for multimodal transportation}.
\newblock {\em International Journal of Sustainable Transportation}, 14(3):205--214.

\bibitem[Layeb et~al., 2018]{layeb2018simulation}
Layeb, S.~B., Jaoua, A., Jbira, A., and Makhlouf, Y. (2018).
\newblock {A simulation-optimization approach for scheduling in stochastic freight transportation}.
\newblock {\em Computers \& Industrial Engineering}, 126:99--110.

\bibitem[Le and Lee, 2013]{le2013model}
Le, T. P.~N. and Lee, T.-R. (2013).
\newblock {Model selection with considering the CO 2 emission alone the global supply chain}.
\newblock {\em Journal of Intelligent Manufacturing}, 24:653--672.

\bibitem[Lemar{\'e}chal, 2001]{lemarechal2001lagrangian}
Lemar{\'e}chal, C. (2001).
\newblock {Lagrangian relaxation}.
\newblock {\em Computational combinatorial optimization: optimal or provably near-optimal solutions}, pages 112--156.

\bibitem[Lessmann and Kramer, 2024]{lessmann2024effect}
Lessmann, C. and Kramer, N. (2024).
\newblock {The effect of cap-and-trade on sectoral emissions: Evidence from California}.
\newblock {\em Energy Policy}, 188:114066.

\bibitem[Li and Wang, 2023]{li2023hierarchical}
Li, H. and Wang, Y. (2023).
\newblock {Hierarchical multimodal hub location problem with carbon emissions}.
\newblock {\em Sustainability}, 15(3):1945.

\bibitem[Li et~al., 2023]{li2023optimum}
Li, L., Zhang, Q., Zhang, T., Zou, Y., and Zhao, X. (2023).
\newblock {Optimum Route and Transport Mode Selection of Multimodal Transport with Time Window under Uncertain Conditions}.
\newblock {\em Mathematics}, 11(14):3244.

\bibitem[Li and Zhang, 2020]{li2020integrated}
Li, L. and Zhang, X. (2020).
\newblock {Integrated optimization of railway freight operation planning and pricing based on carbon emission reduction policies}.
\newblock {\em Journal of Cleaner Production}, 263:121316.

\bibitem[Li and Sun, 2022]{li2022path}
Li, M. and Sun, X. (2022).
\newblock {Path optimization of low-carbon container multimodal transport under uncertain conditions}.
\newblock {\em Sustainability}, 14(21):14098.

\bibitem[Li et~al., 2019]{li2019carbon}
Li, X., Kuang, H., and Hu, Y. (2019).
\newblock {Carbon mitigation strategies of port selection and multimodal transport operations—A case study of northeast China}.
\newblock {\em Sustainability}, 11(18):4877.

\bibitem[Liang et~al., 2021]{liang2021multi}
Liang, X., Liu, X., Mei, L., and Zhang, D. (2021).
\newblock {Multi-objective green multimodal transport path optimization with the participation of high-speed railway}.
\newblock In {\em 2021 6th International Conference on Transportation Information and Safety (ICTIS) (ICTIS)}, pages 1072--1077. IEEE.

\bibitem[Liao et~al., 2009]{liao2009comparing}
Liao, C.-H., Tseng, P.-H., and Lu, C.-S. (2009).
\newblock {Comparing carbon dioxide emissions of trucking and intermodal container transport in Taiwan}.
\newblock {\em Transportation Research Part D: Transport and Environment}, 14(7):493--496.

\bibitem[Lin et~al., 2017]{lin2017modeling}
Lin, B., Liu, C., Wang, H., and Lin, R. (2017).
\newblock {Modeling the railway network design problem: A novel approach to considering carbon emissions reduction}.
\newblock {\em Transportation Research Part D: Transport and Environment}, 56:95--109.

\bibitem[Liotta et~al., 2014]{liotta2014optimization}
Liotta, G., Kaihara, T., and Stecca, G. (2014).
\newblock {Optimization and simulation of collaborative networks for sustainable production and transportation}.
\newblock {\em IEEE Transactions on Industrial Informatics}, 12(1):417--424.

\bibitem[Liotta et~al., 2015]{liotta2015optimisation}
Liotta, G., Stecca, G., and Kaihara, T. (2015).
\newblock {Optimisation of freight flows and sourcing in sustainable production and transportation networks}.
\newblock {\em International Journal of Production Economics}, 164:351--365.

\bibitem[Liu et~al., 2019]{liu2019system}
Liu, P., Liu, C., Du, J., and Mu, D. (2019).
\newblock {A system dynamics model for emissions projection of hinterland transportation}.
\newblock {\em Journal of Cleaner Production}, 218:591--600.

\bibitem[Liu, 2023]{liu2023multimodal}
Liu, S. (2023).
\newblock {Multimodal transportation route optimization of cold chain container in time-varying network considering carbon emissions}.
\newblock {\em Sustainability}, 15(5):4435.

\bibitem[Liu et~al., 2014]{liu2014global}
Liu, Z., Meng, Q., Wang, S., and Sun, Z. (2014).
\newblock {Global intermodal liner shipping network design}.
\newblock {\em Transportation Research Part E: Logistics and Transportation Review}, 61:28--39.

\bibitem[Lu et~al., 2019]{lu2019sustainable}
Lu, Y., Lang, M., Yu, X., and Li, S. (2019).
\newblock {A sustainable multimodal transport system: the two-echelon location-routing problem with consolidation in the euro--China expressway}.
\newblock {\em Sustainability}, 11(19):5486.

\bibitem[Luo et~al., 2016]{luo2016revenue}
Luo, T., Gao, L., and Ak{\c{c}}ay, Y. (2016).
\newblock {Revenue management for intermodal transportation: the role of dynamic forecasting}.
\newblock {\em Production and Operations Management}, 25(10):1658--1672.

\bibitem[Ma et~al., 2018]{ma2018optimization}
Ma, Q., Wang, W., Peng, Y., and Song, X. (2018).
\newblock {An optimization approach to the intermodal transportation network in fruit cold chain, considering cost, quality degradation and carbon dioxide footprint}.
\newblock {\em Polish Maritime Research}, 25(1):61--69.

\bibitem[Macharis and Bontekoning, 2004]{macharis2004opportunities}
Macharis, C. and Bontekoning, Y.~M. (2004).
\newblock {Opportunities for OR in intermodal freight transport research: A review}.
\newblock {\em European Journal of Operational Research}, 153(2):400--416.

\bibitem[Macharis et~al., 2010]{macharis2010decision}
Macharis, C., Van~Hoeck, E., Pekin, E., and Van~Lier, T. (2010).
\newblock {A decision analysis framework for intermodal transport: Comparing fuel price increases and the internalisation of external costs}.
\newblock {\em Transportation Research Part A: Policy and Practice}, 44(7):550--561.

\bibitem[Maia and Couto, 2013]{maia2013strategic}
Maia, L.~C. and Couto, A. (2013).
\newblock {Strategic rail network optimization model for freight transportation}.
\newblock {\em Transportation Research Record}, 2378(1):1--12.

\bibitem[Maiyar and Thakkar, 2019]{maiyar2019modelling}
Maiyar, L.~M. and Thakkar, J.~J. (2019).
\newblock {Modelling and analysis of intermodal food grain transportation under hub disruption towards sustainability}.
\newblock {\em International Journal of Production Economics}, 217:281--297.

\bibitem[Maneengam, 2020]{maneengam2020bi}
Maneengam, A. (2020).
\newblock {A bi-objective programming model for multimodal transportation routing problem of bulk cargo transportation}.
\newblock In {\em 2020 IEEE 7th International Conference on Industrial Engineering and Applications (ICIEA)}, pages 890--894. IEEE.

\bibitem[Mangan et~al., 2008]{mangan2008port}
Mangan, J., Lalwani, C., and Fynes, B. (2008).
\newblock {Port-centric logistics}.
\newblock {\em The International Journal of Logistics Management}, 19(1):29--41.

\bibitem[McKinnon et~al., 2015]{mckinnon2015green}
McKinnon, A., Browne, M., Whiteing, A., and Piecyk, M. (2015).
\newblock {\em {Green logistics: Improving the environmental sustainability of logistics}}.
\newblock Kogan Page Publishers.

\bibitem[McKinnon and Piecyk, 2010]{mckinnon2010measuring}
McKinnon, A.~C. and Piecyk, M. (2010).
\newblock {Measuring and managing CO2 emissions in European chemical transport}.
\newblock Technical report, Logistics Research Centre Heriot-Watt University.

\bibitem[Merrina et~al., 2007]{merrina2007intermodal}
Merrina, A., Sparavigna, A., and Wolf, R. (2007).
\newblock {The intermodal networks: A survey on intermodalism}.
\newblock {\em World Review of Intermodal Transportation Research}, 1(3):286--299.

\bibitem[Metcalf, 2021]{metcalf2021carbon}
Metcalf, G.~E. (2021).
\newblock {Carbon taxes in theory and practice}.
\newblock {\em Annual Review of Resource Economics}, 13(1):245--265.

\bibitem[Meyer, 2020]{meyer2020decafreighttrans}
Meyer, T. (2020).
\newblock {Decarbonizing road freight transportation – A bibliometric and network analysis}.
\newblock {\em Transportation Research Part D: Transport and Environment}, 89:102619.

\bibitem[Moreno-Camacho et~al., 2019]{moreno2019sustainability}
Moreno-Camacho, C.~A., Montoya-Torres, J.~R., Jaegler, A., and Gondran, N. (2019).
\newblock {Sustainability metrics for real case applications of the supply chain network design problem: A systematic literature review}.
\newblock {\em Journal of cleaner production}, 231:600--618.

\bibitem[Morrison et~al., 2016]{morrison2016branch}
Morrison, D.~R., Jacobson, S.~H., Sauppe, J.~J., and Sewell, E.~C. (2016).
\newblock {Branch-and-bound algorithms: A survey of recent advances in searching, branching, and pruning}.
\newblock {\em Discrete Optimization}, 19:79--102.

\bibitem[Mostert et~al., 2017]{mostert2017road}
Mostert, M., Caris, A., and Limbourg, S. (2017).
\newblock {Road and intermodal transport performance: the impact of operational costs and air pollution external costs}.
\newblock {\em Research in Transportation Business \& Management}, 23:75--85.

\bibitem[Mostert et~al., 2018]{mostert2018intermodal}
Mostert, M., Caris, A., and Limbourg, S. (2018).
\newblock {Intermodal network design: a three-mode bi-objective model applied to the case of Belgium}.
\newblock {\em Flexible Services and Manufacturing Journal}, 30:397--420.

\bibitem[Mousavi~Ahranjani et~al., 2020]{mousavi2020robust}
Mousavi~Ahranjani, P., Ghaderi, S.~F., Azadeh, A., and Babazadeh, R. (2020).
\newblock {Robust design of a sustainable and resilient bioethanol supply chain under operational and disruption risks}.
\newblock {\em Clean technologies and environmental policy}, 22:119--151.

\bibitem[M{\"u}ller-Merbach, 1981]{muller1981heuristics}
M{\"u}ller-Merbach, H. (1981).
\newblock {Heuristics and their design: a survey}.
\newblock {\em European Journal of Operational Research}, 8(1):1--23.

\bibitem[Nassar et~al., 2023]{nassar2023system}
Nassar, R.~F., Ghisolfi, V., Annema, J.~A., van Binsbergen, A., and Tavasszy, L.~A. (2023).
\newblock {A system dynamics model for analyzing modal shift policies towards decarbonization in freight transportation}.
\newblock {\em Research in Transportation Business \& Management}, page 100966.

\bibitem[Oudani, 2023]{oudani2023combined}
Oudani, M. (2023).
\newblock {A combined multi-objective multi criteria approach for blockchain-based synchromodal transportation}.
\newblock {\em Computers \& Industrial Engineering}, 176:108996.

\bibitem[Pan et~al., 2017]{pan2017multimodal}
Pan, Y., Li, X., Zhang, M., Zhou, M., and Duan, Y. (2017).
\newblock {Multimodal transportation network optimization with environmental and economic performance considered: An ongoing research}.
\newblock In {\em 2017 International Conference on Service Systems and Service Management}, pages 1--6. IEEE.

\bibitem[Pedinotti-Castelle et~al., 2022]{pedinotti2022freight}
Pedinotti-Castelle, M., Pineau, P.-O., Vaillancourt, K., and Amor, B. (2022).
\newblock {Freight transport modal shifts in a TIMES energy model: Impacts of endogenous and exogenous modeling choice}.
\newblock {\em Applied Energy}, 324:119724.

\bibitem[Pishvaee et~al., 2012]{pishvaee2012credibility}
Pishvaee, M.~S., Torabi, S.~A., and Razmi, J. (2012).
\newblock {Credibility-based fuzzy mathematical programming model for green logistics design under uncertainty}.
\newblock {\em Computers \& Industrial Engineering}, 62(2):624--632.

\bibitem[Pizzol, 2019]{pizzol2019deterministic}
Pizzol, M. (2019).
\newblock {Deterministic and stochastic carbon footprint of intermodal ferry and truck freight transport across Scandinavian routes}.
\newblock {\em Journal of cleaner production}, 224:626--636.

\bibitem[Puterman, 1990]{puterman1990markov}
Puterman, M.~L. (1990).
\newblock {Markov decision processes}.
\newblock {\em Handbooks in operations research and management science}, 2:331--434.

\bibitem[Qi et~al., 2022]{qi2022transport}
Qi, Y., Harrod, S., Psaraftis, H.~N., and Lang, M. (2022).
\newblock {Transport service selection and routing with carbon emissions and inventory costs consideration in the context of the Belt and Road Initiative}.
\newblock {\em Transportation Research Part E: Logistics and Transportation Review}, 159:102630.

\bibitem[Qian et~al., 2021]{qian2021decision}
Qian, Q., Li, D., Gan, M., and Yao, Z. (2021).
\newblock {Decision analysis of the optimal freight structure at provincial level in China}.
\newblock {\em Environmental Science and Pollution Research}, 28(39):54972--54985.

\bibitem[Qiu et~al., 2015]{qiu2015bilevel}
Qiu, X., Lam, J. S.~L., and Huang, G.~Q. (2015).
\newblock {A bilevel storage pricing model for outbound containers in a dry port system}.
\newblock {\em Transportation Research Part E: Logistics and Transportation Review}, 73:65--83.

\bibitem[Qu et~al., 2016]{qu2016sustainability}
Qu, Y., Bekta{\c{s}}, T., and Bennell, J. (2016).
\newblock {Sustainability SI: multimode multicommodity network design model for intermodal freight transportation with transfer and emission costs}.
\newblock {\em Networks and Spatial Economics}, 16:303--329.

\bibitem[Rahiminia et~al., 2023]{rahiminia2023adopting}
Rahiminia, S., Mehrabi, A., Pourseyed~Aghaee, M., and Jamili, A. (2023).
\newblock {Adopting a Bi-level Optimization Method for the Freight Transportation Problem: A Multi-objective Programming Approach}.
\newblock {\em Transportation Research Record}, 2677(2):490--504.

\bibitem[Ramudhin et~al., 2008]{ramudhin2008carbon}
Ramudhin, A., Chaabane, A., Kharoune, M., and Paquet, M. (2008).
\newblock {Carbon market sensitive green supply chain network design}.
\newblock In {\em 2008 IEEE International Conference on Industrial Engineering and Engineering Management}, pages 1093--1097.

\bibitem[Resat and Turkay, 2015]{resat2015design}
Resat, H.~G. and Turkay, M. (2015).
\newblock {Design and operation of intermodal transportation network in the Marmara region of Turkey}.
\newblock {\em Transportation Research Part E: Logistics and Transportation Review}, 83:16--33.

\bibitem[Resat and Turkay, 2019]{resat2019discrete}
Resat, H.~G. and Turkay, M. (2019).
\newblock {A discrete-continuous optimization approach for the design and operation of synchromodal transportation networks}.
\newblock {\em Computers \& Industrial Engineering}, 130:512--525.

\bibitem[Rezaee et~al., 2017]{rezaee2017green}
Rezaee, A., Dehghanian, F., Fahimnia, B., and Beamon, B. (2017).
\newblock {Green supply chain network design with stochastic demand and carbon price}.
\newblock {\em Annals of operations research}, 250:463--485.

\bibitem[Ricci and Black, 2005]{ricci2005social}
Ricci, A. and Black, I. (2005).
\newblock {The social costs of intermodal freight transport}.
\newblock {\em Research in Transportation Economics}, 14:245--285.

\bibitem[Ritchie, 2024]{owid-global-aviation-emissions}
Ritchie, H. (2024).
\newblock {What share of global CO$_2$ emissions come from aviation?}
\newblock {\em Our World in Data}.
\newblock https://ourworldindata.org/global-aviation-emissions.

\bibitem[Rodrigues et~al., 2014]{rodrigues2014assessing}
Rodrigues, V.~S., Beresford, A., Pettit, S., Bhattacharya, S., and Harris, I. (2014).
\newblock {Assessing the cost and CO2e impacts of rerouteing UK import containers}.
\newblock {\em Transportation Research Part A: Policy and Practice}, 61:53--67.

\bibitem[Rodrigues et~al., 2015]{rodrigues2015uk}
Rodrigues, V.~S., Pettit, S., Harris, I., Beresford, A., Piecyk, M., Yang, Z., and Ng, A. (2015).
\newblock {UK supply chain carbon mitigation strategies using alternative ports and multimodal freight transport operations}.
\newblock {\em Transportation Research Part E: Logistics and Transportation Review}, 78:40--56.

\bibitem[Roso, 2007]{roso2007evaluation}
Roso, V. (2007).
\newblock {Evaluation of the dry port concept from an environmental perspective: A note}.
\newblock {\em Transportation Research Part D: Transport and Environment}, 12(7):523--527.

\bibitem[Rudi et~al., 2016]{rudi2016freight}
Rudi, A., Fr{\"o}hling, M., Zimmer, K., and Schultmann, F. (2016).
\newblock {Freight transportation planning considering carbon emissions and in-transit holding costs: a capacitated multi-commodity network flow model}.
\newblock {\em EURO Journal on Transportation and Logistics}, 5(2):123--160.

\bibitem[Santos et~al., 2015]{santos2015impact}
Santos, B.~F., Limbourg, S., and Carreira, J.~S. (2015).
\newblock {The impact of transport policies on railroad intermodal freight competitiveness--The case of Belgium}.
\newblock {\em Transportation Research Part D: Transport and Environment}, 34:230--244.

\bibitem[Sawadogo et~al., 2012]{sawadogo2012reducing}
Sawadogo, M., Anciaux, D., and Daniel, R. (2012).
\newblock {Reducing intermodal transportation impacts on society and environment by path selection: a multiobjective shortest path approach}.
\newblock {\em IFAC Proceedings Volumes}, 45(6):505--513.

\bibitem[Shams et~al., 2023]{shams2023game}
Shams, H.-R., Tamannaei, M., and Zarei, H. (2023).
\newblock {A game-theoretic approach to designing carbon regulations in a duopoly freight transportation market: Road and multimodal road-rail competitive systems}.
\newblock {\em Environmental Science and Pollution Research}, 30(51):111284--111308.

\bibitem[Shepherd, 2014]{shepherd2014review}
Shepherd, S. (2014).
\newblock {A review of system dynamics models applied in transportation}.
\newblock {\em Transportmetrica B: Transport Dynamics}, 2(2):83--105.

\bibitem[Shoukat and Xiaoqiang, 2023]{shoukat2023sustainable}
Shoukat, R. and Xiaoqiang, Z. (2023).
\newblock {Sustainable logistics network optimization from dry ports to seaport: a case study from Pakistan}.
\newblock {\em Transportation Research Record}, 2677(3):302--318.

\bibitem[Soysal et~al., 2014]{soysal2014modelling}
Soysal, M., Bloemhof-Ruwaard, J.~M., and Van Der~Vorst, J.~G. (2014).
\newblock {Modelling food logistics networks with emission considerations: The case of an international beef supply chain}.
\newblock {\em International Journal of Production Economics}, 152:57--70.

\bibitem[Steimle et~al., 2021]{steimle2021multi}
Steimle, L.~N., Kaufman, D.~L., and Denton, B.~T. (2021).
\newblock {Multi-model Markov decision processes}.
\newblock {\em IISE Transactions}, 53(10):1124--1139.

\bibitem[Sun, 2020]{sun2020green}
Sun, Y. (2020).
\newblock {Green and Reliable Freight Routing Problem in the Road-Rail Intermodal Transportation Network with Uncertain Parameters: A Fuzzy Goal Programming Approach}.
\newblock {\em Journal of Advanced Transportation}, 2020(1):7570686.

\bibitem[Sun et~al., 2018]{sun2018time}
Sun, Y., Hru{\v{s}}ovsk{\`y}, M., Zhang, C., Lang, M., et~al. (2018).
\newblock {A time-dependent fuzzy programming approach for the green multimodal routing problem with rail service capacity uncertainty and road traffic congestion}.
\newblock {\em Complexity}, 2018.

\bibitem[Sun and Lang, 2015]{sun2015modeling}
Sun, Y. and Lang, M. (2015).
\newblock {Modeling the Multicommodity Multimodal Routing Problem with Schedule-Based Services and Carbon Dioxide Emission Costs}.
\newblock {\em Mathematical Problems in Engineering}, 2015(1):406218.

\bibitem[Sun et~al., 2019]{sun2019bi}
Sun, Y., Li, X., Liang, X., and Zhang, C. (2019).
\newblock {A bi-objective fuzzy credibilistic chance-constrained programming approach for the hazardous materials road-rail multimodal routing problem under uncertainty and sustainability}.
\newblock {\em Sustainability}, 11(9):2577.

\bibitem[Sun et~al., 2023a]{sun2023study}
Sun, Y., Lu, Z., Lian, F., and Yang, Z. (2023a).
\newblock {Study of channel upgrades and ship choices of river-shipping of port access-transportation}.
\newblock {\em Transportation Research Part D: Transport and Environment}, 119:103733.

\bibitem[Sun et~al., 2023b]{sun2023modeling}
Sun, Y., Sun, G., Huang, B., and Ge, J. (2023b).
\newblock {Modeling a Carbon-Efficient Road--Rail Intermodal Routing Problem with Soft Time Windows in a Time-Dependent and Fuzzy Environment by Chance-Constrained Programming}.
\newblock {\em Systems}, 11(8):403.

\bibitem[Sun et~al., 2022]{sun2022green}
Sun, Y., Yu, N., and Huang, B. (2022).
\newblock {Green road--rail intermodal routing problem with improved pickup and delivery services integrating truck departure time planning under uncertainty: An interactive fuzzy programming approach}.
\newblock {\em Complex \& Intelligent Systems}, 8(2):1459--1486.

\bibitem[Tao and Wu, 2021]{tao2021energy}
Tao, X. and Wu, Q. (2021).
\newblock {Energy consumption and CO2 emissions in hinterland container transport}.
\newblock {\em Journal of cleaner production}, 279:123394.

\bibitem[Tavasszy et~al., 2017]{tavasszy2017intermodality}
Tavasszy, L., Behdani, B., and Konings, R. (2017).
\newblock {Intermodality and synchromodality}.
\newblock In {\em Ports and Networks}, pages 251--266. Routledge.

\bibitem[Thore and Iannone, 2010]{thore2010economic}
Thore, S. and Iannone, F. (2010).
\newblock {An economic logistics model for the multimodal inland distribution of maritime containers}.
\newblock {\em An Economic Logistics Model for the Multimodal Inland Distribution of Maritime Containers}, pages 1000--1046.

\bibitem[Tiwari et~al., 2021]{tiwari2021freight}
Tiwari, S., Wee, H.~M., Zhou, Y., and Tjoeng, L. (2021).
\newblock {Freight consolidation and containerization strategy under business as usual scenario \& carbon tax regulation}.
\newblock {\em Journal of Cleaner Production}, 279:123270.

\bibitem[Tran et~al., 2017]{tran2017container}
Tran, N.~K., Haasis, H.-D., and Buer, T. (2017).
\newblock {Container shipping route design incorporating the costs of shipping, inland/feeder transport, inventory and CO2 emission}.
\newblock {\em Maritime Economics \& Logistics}, 19(4):667--694.

\bibitem[Tsao and Linh, 2018]{tsao2018seaport}
Tsao, Y.-C. and Linh, V.~T. (2018).
\newblock {Seaport-dry port network design considering multimodal transport and carbon emissions}.
\newblock {\em Journal of Cleaner Production}, 199:481--492.

\bibitem[Tsao and Thanh, 2019]{tsao2019multi}
Tsao, Y.-C. and Thanh, V.-V. (2019).
\newblock {A multi-objective mixed robust possibilistic flexible programming approach for sustainable seaport-dry port network design under an uncertain environment}.
\newblock {\em Transportation Research Part E: Logistics and Transportation Review}, 124:13--39.

\bibitem[Tyroller, 2021]{tyroller2021challenges}
Tyroller, L. (2021).
\newblock {Challenges to intermodal transportation: a case study}.
\newblock Master's thesis, Hanken School of Economics.

\bibitem[VWT, 2022]{VirtualTower}
VWT (2022).
\newblock Virtual watchtower network.
\newblock Accessed on October 4th, 2024.

\bibitem[Wang et~al., 2020a]{wang2020transportation}
Wang, C.-N., Dang, T.-T., Le, T.~Q., and Kewcharoenwong, P. (2020a).
\newblock {Transportation optimization models for intermodal networks with fuzzy node capacity, detour factor, and vehicle utilization constraints}.
\newblock {\em Mathematics}, 8(12):2109.

\bibitem[Wang et~al., 2021]{wang2021bi}
Wang, D., Zhou, L., Zhang, H., and Liang, X. (2021).
\newblock {A bi-level model for green freight transportation pricing strategy considering enterprise profit and carbon emissions}.
\newblock {\em Sustainability}, 13(12):6514.

\bibitem[Wang et~al., 2015]{wang2015effects}
Wang, M., Liu, K., Choi, T.-M., and Yue, X. (2015).
\newblock {Effects of carbon emission taxes on transportation mode selections and social welfare}.
\newblock {\em IEEE Transactions on Systems, Man, and Cybernetics: Systems}, 45(11):1413--1423.

\bibitem[Wang et~al., 2020b]{wang2020modelling}
Wang, Q.-Z., Chen, J.-M., Tseng, M.-L., Luan, H.-M., and Ali, M.~H. (2020b).
\newblock {Modelling green multimodal transport route performance with witness simulation software}.
\newblock {\em Journal of Cleaner Production}, 248:119245.

\bibitem[Wang et~al., 2020c]{wang2020integrated}
Wang, W., Xu, X., Jiang, Y., Xu, Y., Cao, Z., and Liu, S. (2020c).
\newblock {Integrated scheduling of intermodal transportation with seaborne arrival uncertainty and carbon emission}.
\newblock {\em Transportation Research Part D: Transport and Environment}, 88:102571.

\bibitem[Wei and Dong, 2019]{wei2019import}
Wei, H. and Dong, M. (2019).
\newblock {Import-export freight organization and optimization in the dry-port-based cross-border logistics network under the Belt and Road Initiative}.
\newblock {\em Computers \& Industrial Engineering}, 130:472--484.

\bibitem[Wieberneit, 2008]{wieberneit2008service}
Wieberneit, N. (2008).
\newblock {Service network design for freight transportation: a review}.
\newblock {\em OR spectrum}, 30(1):77--112.

\bibitem[Winebrake et~al., 2008]{winebrake2008assessing}
Winebrake, J.~J., Corbett, J.~J., Falzarano, A., Hawker, J.~S., Korfmacher, K., Ketha, S., and Zilora, S. (2008).
\newblock {Assessing energy, environmental, and economic tradeoffs in intermodal freight transportation}.
\newblock {\em Journal of the Air \& Waste Management Association}, 58(8):1004--1013.

\bibitem[Wnuk and Garrepalli, 2018]{wnuk2018knowledge}
Wnuk, K. and Garrepalli, T. (2018).
\newblock {Knowledge management in software testing: a systematic snowball literature review}.
\newblock {\em e-Informatica Software Engineering Journal}, 12(1):51--78.

\bibitem[Wu et~al., 2021]{wu2021research}
Wu, J., Wang, Y., Li, W., and Wu, H. (2021).
\newblock {Research on Green Transport Mode of Chinese Bulk Cargo Based on Fourth-Party Logistics}.
\newblock {\em Journal of Advanced Transportation}, 2021(1):6142226.

\bibitem[Wu and Zhang, 2023a]{wu2023evaluation}
Wu, Y. and Zhang, R. (2023a).
\newblock {Evaluation Model for a Port Hinterland Intermodal Freight Network Considering Environmental Impacts and Capacity Constraints}.
\newblock {\em Transportation Research Record}, 2677(2):462--478.

\bibitem[Wu and Zhang, 2023b]{wu2023game}
Wu, Y. and Zhang, R. (2023b).
\newblock {Game-theoretical method toward dry port multilevel location considering capacity constraints and shippers’ choice behavior}.
\newblock {\em Transportation Research Record}, 2677(11):64--82.

\bibitem[Xie et~al., 2022]{xie2022research}
Xie, F.-J., Feng, R.-C., and Zhou, X.-Y. (2022).
\newblock {Research on the Optimization of Cross-Border Logistics Paths of the “Belt and Road” in the Inland Regions}.
\newblock {\em Journal of Advanced Transportation}, 2022(1):5776334.

\bibitem[Xu et~al., 2018]{xu2018modelling}
Xu, X., Zhang, Q., Wang, W., Peng, Y., Song, X., and Jiang, Y. (2018).
\newblock {Modelling port competition for intermodal network design with environmental concerns}.
\newblock {\em Journal of Cleaner Production}, 202:720--735.

\bibitem[Yang et~al., 2023]{yang_multi-objective_2023}
Yang, L., Zhang, C., and Wu, X. (2023).
\newblock {Multi-Objective Path Optimization of Highway-Railway Multimodal Transport Considering Carbon Emissions}.
\newblock {\em Applied Sciences}, 13(8):4731.

\bibitem[Yang et~al., 2021]{yang2021coastal}
Yang, Z., Xin, X., Chen, K., and Yang, A. (2021).
\newblock {Coastal container multimodal transportation system shipping network design—toll policy joint optimization model}.
\newblock {\em Journal of Cleaner Production}, 279:123340.

\bibitem[Yin et~al., 2021]{yin2021interrelations}
Yin, C., Ke, Y., Chen, J., and Liu, M. (2021).
\newblock {Interrelations between sea hub ports and inland hinterlands: Perspectives of multimodal freight transport organization and low carbon emissions}.
\newblock {\em Ocean \& Coastal Management}, 214:105919.

\bibitem[Yin et~al., 2024]{yin2024low}
Yin, C., Zhang, Z.-A., Fu, X., and Ge, Y.-E. (2024).
\newblock {A low-carbon transportation network: Collaborative effects of a rail freight subsidy and carbon trading mechanism}.
\newblock {\em Transportation Research Part A: Policy and Practice}, 184:104066.

\bibitem[Yu et~al., 2018]{yu2018achievement}
Yu, S., Zheng, S., and Li, X. (2018).
\newblock {The achievement of the carbon emissions peak in China: The role of energy consumption structure optimization}.
\newblock {\em Energy Economics}, 74:693--707.

\bibitem[Zhang et~al., 2017]{zhang2017multimodal}
Zhang, D., He, R., Li, S., and Wang, Z. (2017).
\newblock {A multimodal logistics service network design with time windows and environmental concerns}.
\newblock {\em PLoS One}, 12(9):e0185001.

\bibitem[Zhang et~al., 2018]{zhang2018joint}
Zhang, D., Zhan, Q., Chen, Y., and Li, S. (2018).
\newblock {Joint optimization of logistics infrastructure investments and subsidies in a regional logistics network with CO2 emission reduction targets}.
\newblock {\em Transportation Research Part D: Transport and Environment}, 60:174--190.

\bibitem[Zhang et~al., 2024]{zhang2024sparrow}
Zhang, H., Huang, Q., Ma, L., and Zhang, Z. (2024).
\newblock {Sparrow search algorithm with adaptive t distribution for multi-objective low-carbon multimodal transportation planning problem with fuzzy demand and fuzzy time}.
\newblock {\em Expert Systems with Applications}, 238:122042.

\bibitem[Zhang et~al., 2011]{zhang2011mode}
Zhang, J., Ding, H.~W., Wang, X.~Q., Yin, W.~J., Zhao, T.~Z., and Dong, J. (2011).
\newblock {Mode choice for the intermodal transportation considering carbon emissions}.
\newblock In {\em Proceedings of 2011 ieee international conference on service operations, logistics and informatics}, pages 297--301.

\bibitem[Zhang et~al., 2015]{zhang2015freight}
Zhang, M., Janic, M., and Tavasszy, L.~A. (2015).
\newblock {A freight transport optimization model for integrated network, service, and policy design}.
\newblock {\em Transportation Research Part E: Logistics and Transportation Review}, 77:61--76.

\bibitem[Zhang and Pel, 2016]{zhang2016synchromodal}
Zhang, M. and Pel, A. (2016).
\newblock {Synchromodal hinterland freight transport: Model study for the port of Rotterdam}.
\newblock {\em Journal of Transport Geography}, 52:1--10.

\bibitem[Zhang et~al., 2013]{zhang2013optimization}
Zhang, M., Wiegmans, B., and Tavasszy, L. (2013).
\newblock {Optimization of multimodal networks including environmental costs: a model and findings for transport policy}.
\newblock {\em Computers in Industry}, 64(2):136--145.

\bibitem[Zhang and Chen, 2023]{zhang2023research}
Zhang, T. and Chen, N. (2023).
\newblock {Research on Multimodal Transportation Path Optimization Based on Improved Genetic Algorithm}.
\newblock In {\em 2023 7th International Conference on Transportation Information and Safety (ICTIS)}, pages 2034--2043. IEEE.

\bibitem[Zhang et~al., 2022a]{zhang2022acyclic}
Zhang, W., Uhan, N.~A., Dessouky, M., and Toriello, A. (2022a).
\newblock {Acyclic mechanism design for freight consolidation}.
\newblock {\em Transportation Science}, 56(3):571--584.

\bibitem[Zhang et~al., 2021]{zhang2021low}
Zhang, X., Jin, F.-Y., Yuan, X.-M., and Zhang, H.-Y. (2021).
\newblock {Low-carbon multimodal transportation path optimization under dual uncertainty of demand and time}.
\newblock {\em Sustainability}, 13(15):8180.

\bibitem[Zhang et~al., 2022b]{zhang2022preference}
Zhang, Y., Atasoy, B., and Negenborn, R.~R. (2022b).
\newblock {Preference-based multi-objective optimization for synchromodal transport using adaptive large neighborhood search}.
\newblock {\em Transportation Research Record}, 2676(3):71--87.

\bibitem[Zhang et~al., 2022c]{zhang2022flexible}
Zhang, Y., Guo, W., Negenborn, R.~R., and Atasoy, B. (2022c).
\newblock {Synchromodal transport planning with flexible services: Mathematical model and heuristic algorithm}.
\newblock {\em Transportation Research Part C: Emerging Technologies}, 140:103711.

\bibitem[Zhang et~al., 2022d]{zhang2022collaborative}
Zhang, Y., Heinold, A., Meisel, F., Negenborn, R.~R., and Atasoy, B. (2022d).
\newblock {Collaborative planning for intermodal transport with eco-label preferences}.
\newblock {\em Transportation Research Part D: Transport and Environment}, 112:103470.

\bibitem[Zhang et~al., 2022e]{zhang2022heterogeneous}
Zhang, Y., Li, X., van Hassel, E., Negenborn, R.~R., and Atasoy, B. (2022e).
\newblock {Synchromodal transport planning considering heterogeneous and vague preferences of shippers}.
\newblock {\em Transportation Research Part E: Logistics and Transportation Review}, 164:102827.

\bibitem[Zhao et~al., 2020]{zhao2020study}
Zhao, J., Zhu, X., and Wang, L. (2020).
\newblock {Study on scheme of outbound railway container organization in rail-water intermodal transportation}.
\newblock {\em Sustainability}, 12(4):1519.

\bibitem[Zhao et~al., 2018]{zhao2018stochastic}
Zhao, Y., Xue, Q., and Zhang, X. (2018).
\newblock {Stochastic empty container repositioning problem with CO2 emission considerations for an intermodal transportation system}.
\newblock {\em Sustainability}, 10(11):4211.

\bibitem[Zhong et~al., 2023]{zhong2023system}
Zhong, H., Chen, W., and Gu, Y. (2023).
\newblock {A system dynamics model of port hinterland intermodal transport: A case study of Guangdong-Hong Kong-Macao Greater Bay Area under different carbon taxation policies}.
\newblock {\em Research in Transportation Business \& Management}, 49:100987.

\bibitem[Zhou et~al., 2018]{zhou2018capacitated}
Zhou, M., Duan, Y., Yang, W., Pan, Y., and Zhou, M. (2018).
\newblock {Capacitated multi-modal network flow models for minimizing total operational cost and CO2e emission}.
\newblock {\em Computers \& Industrial Engineering}, 126:361--377.

\bibitem[Zuo et~al., 2023]{zuo2023using}
Zuo, D., Liang, Q., Zhan, S., Huang, W., Yang, S., and Wang, M. (2023).
\newblock {Using energy consumption constraints to control the freight transportation structure in China (2021--2030)}.
\newblock {\em Energy}, 262:125512.

\end{thebibliography}

\newpage
\begin{landscape}
\section{Comprehensive Analysis} \label{Comprehensive_Analysis}

    \Centering

    \footnotesize

    \begin{longtable}{p{0.6in} p{1.65in} p{0.7cm} p{0.75cm} p{0.25cm} p{0.25cm} p{0.25cm} p{0.25cm} p{0.25cm} p{0.25cm} p{0.25cm} p{0.25cm} p{0.25cm} p{0.25cm} p{0.25cm} p{0.25cm} p{0.25cm} p{0.45in} p{1.55in}}
        \caption{Summary of Solution Techniques}\\
     \label{tab:Solution}\\
        \toprule
        \multirow{2}{4em}{Classification} & \multirow{2}{12em}{Publications} & \multirow{2}{3em}{S/M/B} & \multirow{2}{4em}{DE/ST} & \multicolumn{13}{c}{Objective Function}  & \multirow{2}{6em}{Solution Technique } & \multirow{2}{12em}{Specific Technique and Tool}\\  
        \cline{5-17} \\ 
        & & & & TC & EC & HC & PC & IC & FC & O & TT & HT & WT & TE & HE & FE & & \\
        \midrule
        SNDP-TDP & \cite{comer2010marine} &S & DE & \checkmark &  &\checkmark &  &  &  & & \checkmark & \checkmark& & \checkmark &  & & DSS & Dijkstra's algorithm/ArcGIS  \\
        \midrule
        SNDP-SSP  & \cite{bauer2010minimizing} &S & DE  & \checkmark & \checkmark & &  &  &  &\checkmark &  & & &  &  &  
        & E   & Branch and bound/CPLEX  \\
        \midrule
        SNDP-TDP & \cite{chang2010optimization} &S & DE  & \checkmark &\checkmark  &\checkmark &  & \checkmark &  &  &  & & &  &  &  & E  & Branch and bound/CPLEX  \\
        \midrule
        NDP-HSP & \cite{henttu2011financial} &S & DE  &\checkmark  & \checkmark & &  &  &\checkmark  &\checkmark &  & & &  &  &  & E & Gravitational models \\
        \midrule
        SNDP-TDP & \cite{zhang2011mode} &S & DE  &\checkmark  & \checkmark & \checkmark& \checkmark &  &  & &  & & &  &  &  & E & Branch and bound/CPLEX  \\
        \midrule
        SCND & \cite{chaabane2011designing} &M & DE  &\checkmark  &  & &  &  &\checkmark  & \checkmark&  & & &  \checkmark&  &  & E& Branch and bound/CPLEX  \\
        \midrule                             
        NDP-HSP & \cite{iannone2012private} &S & DE  &  & \checkmark & &  &  &  & \checkmark&  & & &  &  &  & E& Branch and bound/CPLEX  \\
        \midrule                                       
        SNDP-TDP & \cite{sawadogo2012reducing} &M & DE  &\checkmark  &  & \checkmark& \checkmark & \checkmark &  & \checkmark&\checkmark  & & & \checkmark &  &  &  MH& Ant Colony Algorithm  \\
        \midrule                                  
        SCND & \cite{holmgren2012tapas} &S  & ST  & \checkmark & \checkmark &\checkmark &  \checkmark&  \checkmark&  \checkmark& \checkmark&  & & &  &  &  &  SIM&  Agent Based Simulation \\
        \midrule                                          
        NDP-HSP & \cite{lattila2013hinterland} &S,M & DE  & \checkmark &  &\checkmark &  &  &  & &  & & &  \checkmark&\checkmark  &  &  SIM& DES  \\
        \midrule                                          
        SCND & \cite{chaabane2012design} &M & DE  &  &  & &  &\checkmark  & \checkmark & \checkmark&  & & &\checkmark  &  &  & E  & Branch and bound/LINGO  \\
        \midrule                                       
        SCND & \cite{pishvaee2012credibility} &M & ST  & \checkmark &  & &  &  & \checkmark & &  & & &\checkmark &  &  & HE & Credibility-based Fuzzy CC  \\
        \midrule                                      
        NDP-CEP & \cite{maia2013strategic} & M & DE  &  \checkmark&  & &  &  &  & &  & & & \checkmark &  &  & HE &  Greedy Algorithm \\
        \midrule                                    
        NDP-CEP& \cite{kim2013multimodal} & B & DE  & \checkmark & \checkmark & &  &  &  & & \checkmark & & &  &  &  & MH & GA-based Algorithm  \\
        \midrule                                        
        NDP-LP & \cite{zhang2013optimization} & B & DE  & \checkmark &\checkmark  &\checkmark &  &  &  & \checkmark&  & & &  &  &  & HYB & All-or-Nothing, GA  \\
        \midrule                              
        SCND & \cite{le2013model} & S & DE  &\checkmark  &  & &  &\checkmark  &  & \checkmark&  & & & \checkmark  &  &  & E &  Branch and bound/LINGO  \\
        \midrule                                    
        SCND & \cite{fahimnia2013impact} & S & DE  & \checkmark & \checkmark  &\checkmark & \checkmark & \checkmark & \checkmark &\checkmark &  & & &  &  &  & E& Branch and bound/CPLEX  \\
        \midrule                                        
        SCND & \cite{liotta2014optimization} & S & DE  & \checkmark & \checkmark & &  &  &  &\checkmark &  & & &  &  &  &  HYB &  DES, Branch and bound/CS \\        
        \midrule
        SNDP-SSP & \cite{kengpol2014development} & M & DE  & \checkmark  &  & &  &  &  & & \checkmark & & &\checkmark  &  &  &  DSS & AHP, DEA, Spreadsheet Software \\
        \midrule
        SNDP-RMP & \cite{hoen2014switching} & S & DE  & \checkmark & \checkmark & \checkmark &  &  &  & &  & & &  &  &  & O &  Lagrangian Relaxation \\
        \midrule
        PS & \cite{chen2014design} & B & DE  & \checkmark &  & &  &  &  & \checkmark &  & & & \checkmark  &  &  & MH & GA \& Frank-Wolfe Algorithm  \\
        \midrule
        SCND & \cite{soysal2014modelling} & M & DE & \checkmark &  & &  & \checkmark &  & &  & & & \checkmark  &  &  &E & Branch and bound/CPLEX  \\
        \midrule
        
        \\

        \multirow{2}{4em}{Classification} & \multirow{2}{12em}{Publications} & \multirow{2}{3em}{S/M/B} & \multirow{2}{4em}{DE/ST} & \multicolumn{13}{c}{Objective Function}  & \multirow{2}{6em}{Solution Technique } & \multirow{2}{12em}{Specific Technique and Tool}\\  
        \cline{5-17} \\ 
        & & & & TC & EC & HC & PC & IC & FC & O & TT & HT & WT & TE & HE & FE & & \\
        \midrule        
        NDP-CEP & \cite{rodrigues2014assessing} & S  & DE  & \checkmark &  & &  &  &  & &  & & & \checkmark  &  &  &E & Spreadsheet Software   \\
        \midrule
        NDP-LP & \cite{bouchery2015cost} & S & DE  & \checkmark  &  & &  &  &  & &  & &   &\checkmark &&& E & Branch and Bound/CS \\
        \midrule
        SNDP-TDP & \cite{duan2015carbon} & S/M & DE  &\checkmark  &\checkmark  & & \checkmark &  &  & &  & &   &&&& E & Branch and Bound/CS \\
        \midrule
        PS & \cite{wang2015effects} & B  & DE   & \checkmark  & \checkmark  & &  &   & & \checkmark  &  &  &&&&&GT  & Backward Induction \\
        \midrule
        SNDP-SSP & \cite{sun2015modeling} & S & DE  & \checkmark & \checkmark & \checkmark&  &\checkmark  &   &  &  &  &&&& & E & Branch and bound/LINGO  \\
        \midrule
        SCND & \cite{liotta2015optimisation} & S & DE  & \checkmark & \checkmark & &  &  &  \checkmark &\checkmark  &  &  &&&& & E & Branch and bound/CPLEX  \\
        \midrule
        SCND & \cite{fahimnia2015tactical} & M & DE   & \checkmark &  & & \checkmark & \checkmark &  \checkmark &\checkmark  &  &  &&\checkmark&& & MH & Nested Integrated Cross-Entropy  \\
        \midrule
        SNDP-SSP & \cite{zhang2015freight} & B & DE  & \checkmark & \checkmark & &  &  &  \checkmark &  &  & \checkmark &&&& & HYB   & GA, All or Nothing Algorithm  \\
        \midrule
        SNDP-SSP & \cite{qu2016sustainability} & S/M & DE  &\checkmark  & \checkmark & \checkmark&  &  & \checkmark & &  & & &  &  &  & E  & Branch and bound/CPLEX  \\
        \midrule
        SNDP-SSP & \cite{rudi2016freight} & S & DE  & \checkmark &  & \checkmark&  & \checkmark &  & & \checkmark & & & \checkmark &  &   & DSS & MCDM/GAMS+MS Excel  \\
        \midrule   
        SNDP-AMP & \cite{lam2016market} & M & DE  & \checkmark &  & \checkmark&  & \checkmark & \checkmark &\checkmark &\checkmark  & \checkmark& &  &    &  &  E& Branch and bound/CPLEX  \\
        \midrule
        SNDP-SSP & \cite{demir2016green} & S & ST  &  \checkmark&\checkmark  &\checkmark &\checkmark  &  &  & &  & & &  &  &  & E  & Branch and bound/CPLEX  \\
        \midrule
        SNDP-SSP & \cite{baykasouglu2016multi} & S/M & DE  & \checkmark &  & &  &  &  & \checkmark& \checkmark &\checkmark & \checkmark&\checkmark  &  &  & E  & Compromise and fuzzy GP/LINGO  \\
        \midrule   
        SNDP-TDP & \cite{pan2017multimodal} & M & DE  &\checkmark &  & &  &  &  & &  & & & \checkmark &  &  & E  & Branch and bound/CPLEX  \\
       \midrule   
       SNDP-TDP & \cite{ji2017hybrid} & M & DE  &  \checkmark&\checkmark  & \checkmark&  &  &  & & \checkmark & & &  &  &  &   MH   & HEDA \& LS  \\
        \midrule 
       NDP-HSP & \cite{tran2017container} & S & DE  & \checkmark & \checkmark & \checkmark &  &  & \checkmark & &  & & &  &  &  &  HE & Brute-force and Greedy algorithm  \\
        \midrule
        NDP-LP & \cite{mostert2018intermodal} & M & DE & \checkmark &  & \checkmark&  &  &  &\checkmark &  & & & \checkmark & \checkmark & & E & Branch and bound/CPLEX  \\
        \midrule
        NDP+SNDP & \cite{zhang2017multimodal} & S & DE &\checkmark  & \checkmark &\checkmark &\checkmark  &  & \checkmark & &  & & &  &  &  &  MH & GA  \\
        \midrule
        NDP-CEP & \cite{lin2017modeling} & B & DE  & \checkmark & \checkmark & & \checkmark  &  &  & \checkmark &  & & &  &  &  & E  & Branch and bound/LINGO  \\
        \midrule
        SCND & \cite{rezaee2017green} & S & ST  & \checkmark & \checkmark & & \checkmark &  &  & \checkmark  &  & & &  &  &  & E& Branch and bound/CPLEX  \\
        \midrule 
        NDP-LP & \cite{tsao2018seaport} & B & DE  & \checkmark & \checkmark & \checkmark &  & \checkmark & \checkmark & \checkmark &  & & &  &  &  & GT & Continuous Approximation Model \\
        \midrule
        SNDP-SSP & \cite{sun2018time} & S/M & ST  & \checkmark & \checkmark & \checkmark& \checkmark & \checkmark &  & &  & & &  &  &  & E & Fuzzy CCP/LINGO  \\
        \midrule
        SNDP-TDP & \cite{harris2018impact} & S & DE  & \checkmark &  & &  &  &  & &  & & & \checkmark &  &  &  E & MS Excel Optimization \\
        \midrule
        SNDP-AMP& \cite{zhao2018stochastic} & S & ST  &  &\checkmark  & &  &  &  &\checkmark &  & & &  &  &  &  MH& Two-phase Tabu Search \\
        \midrule
        \\
        
        \multirow{2}{4em}{Classification} & \multirow{2}{12em}{Publications} & \multirow{2}{3em}{S/M/B} & \multirow{2}{4em}{DE/ST} & \multicolumn{13}{c}{Objective Function}  & \multirow{2}{6em}{Solution Technique } & \multirow{2}{12em}{Specific Technique and Tool}\\  
        \cline{5-17} \\ 
        & & & & TC & EC & HC & PC & IC & FC & O & TT & HT & WT & TE & HE & FE & & \\
        \midrule
        
        SNDP-TDP& \cite{ma2018optimization} & M & DE  &\checkmark  &  &\checkmark &  &\checkmark  &  & \checkmark& \checkmark  & & &  &  &  & E & Branch and bound/ CPLEX \\
        \midrule           
        SNDP-TDP& \cite{zhou2018capacitated} & S/M & DE  &\checkmark  &  &\checkmark &  &  &  & &  & & &  &  &  & E& Branch and bound/ CPLEX \\
        \midrule
        SNDP-SSP & \cite{layeb2018simulation} & S & ST  & \checkmark & \checkmark & &\checkmark  &  &  & &  & & &  &  &  & HYB & DES/Arena, Search Algorithm /Optquest \\
        \midrule
        SNDP-SSP & \cite{hruvsovsky2018hybrid} & S & ST  & \checkmark & \checkmark & & \checkmark & \checkmark &  & &  & & &  &  &  & HYB & Branch and bound/CPLEX, ABM+DES/ANYLOGIC\\
        \midrule
        NDP-LP & \cite{xu2018modelling} & B & DE  & \checkmark & \checkmark & &  &  &  & \checkmark & \checkmark & & &  &  &  & GT &  Nash Equilibrium\\
        \midrule
        PS & \cite{zhang2018joint} & B & DE  & \checkmark &  & &  &  &  & &\checkmark  & & &\checkmark  &  &  & MH & GFWHA \\
        \midrule
        SCND & \cite{haddadsisakht2018closed} & S & ST  & \checkmark & \checkmark & &  &  & \checkmark  & \checkmark& & & &  &  &  & O & Benders Decomposition \\
        \midrule
        SCND & \cite{farahani2018decision} & NA & DE  & \checkmark & \checkmark & &  &  &  & &  & & &  &  &  &DSS & Ms Excel \\
        \midrule
        SNDP-SSP& \cite{demir2019green} & M & DE  & \checkmark &\checkmark  &\checkmark &  &  &  & &  & & & &  &  & E  & Branch and bound/CPLEX \\
        \midrule
        SNDP-TDP & \cite{fulzele2019model} & M & DE  & \checkmark &  & &  &  &  & &  & & &  \checkmark&  &  & DSS & MCDM, Branch and bound/LINGO \\
        \midrule
        SNDP-TDP & \cite{maiyar2019modelling} & S & DE  & \checkmark &\checkmark  & &  &  &\checkmark  & \checkmark&  & & &  &  &  &  MH& PSO + DE/MATLAB \\  
        \midrule
        SNDP-TDP & \cite{li2019carbon} & S & DE  & \checkmark &\checkmark &\checkmark &  &  &  & &  & & &  &  &  & DSS & MCDM, Branch and bound/CPLEX \\
        \midrule        
        SNDP-SSP & \cite{lu2019sustainable} & S & DE  &  &  & &  &  &  & &  & & &  \checkmark&  &  &  MH & DE, CW \& LS \\
        \midrule
        SNDP-AMP & \cite{wei2019import} & M & DE  & \checkmark &  & \checkmark&  &\checkmark  &  &\checkmark &\checkmark  & \checkmark& &  &  &  & MH & Adaptive-weight GA \\
        \midrule
        PS & \cite{liu2019system} & S & DE  & \checkmark  & \checkmark  & &  &  &  & \checkmark &  & & &  &  &  & SIM  & SD \\
        \midrule
        PS & \cite{choi2019system} & S & DE  & \checkmark & \checkmark & &  &  &  & &  & & &  &  &  & SIM & SD/VENSIM PLE \\
        \midrule
        NDP-LP & \cite{tsao2019multi} & M & ST  & \checkmark & \checkmark & \checkmark&  & \checkmark & \checkmark &\checkmark &  & & &  &  & \checkmark  & E  & Branch and bound/LINGO \\
        \midrule
        SNDP-SSP & \cite{sun2019bi} & M & ST  & \checkmark &  &\checkmark & \checkmark & \checkmark &  & &  & & &  &  &  & E& Fuzzy CCP, Branch and Bound/LINGO\\
        \midrule
        SNDP-AMP & \cite{baykasouglu2019fuzzy} & M & ST  & \checkmark &  & &  &  &\checkmark  &\checkmark & \checkmark & \checkmark& \checkmark& \checkmark &  &  & E&  Hybrid CCP, Branch and Bound/LINGO \\
        \midrule
        SNDP-SSP & \cite{resat2019discrete} & M & DE  &\checkmark  &  & \checkmark& \checkmark &  &  & & \checkmark & \checkmark& \checkmark& \checkmark &  &  & O & Logic-cut Method, Branch and Bound/CPLEX, DICOPT \\
        \midrule
        SCND & \cite{kabadurmus2020sustainable} & S & DE  &  \checkmark& \checkmark  & &  &  & \checkmark & \checkmark &  & & &  &  &  & E& Branch and bound/CPLEX \\
        \midrule
        
        \\
    
    \multirow{2}{4em}{Classification} & \multirow{2}{12em}{Publications} & \multirow{2}{3em}{S/M/B} & \multirow{2}{4em}{DE/ST} & \multicolumn{13}{c}{Objective Function}  & \multirow{2}{6em}{Solution Technique } & \multirow{2}{12em}{Specific Technique and Tool}\\  
    \cline{5-17} \\ 
    & & & & TC & EC & HC & PC & IC & FC & O & TT & HT & WT & TE & HE & FE & & \\
        \midrule
        SNDP-SSP & \cite{sun2020green} & M & ST  & \checkmark &  &\checkmark & \checkmark &\checkmark  &  & &  & & & \checkmark &  &  & E & Branch and Bound/LINGO \\
        \midrule
        
        SNDP-AMP & \cite{abu2020optimization} & M & DE  & \checkmark &  & &  & \checkmark &  & &  & & &\checkmark  &  &  & E & Branch and bound/CPLEX \\
        \midrule
        SNDP-RMP & \cite{li2020integrated} & B & DE  & \checkmark&  & \checkmark&  &  &  & &  & & &  &  &  &  HYB& NSGA-III + Descent algorithm / CVX Toolbox MATLAB\\
        \midrule
        
        SNDP-AMP & \cite{wang2020transportation} & S & ST  & \checkmark & \checkmark & \checkmark&  &  &\checkmark  & &  & & &  &  &  & E  & Branch and bound/CPLEX \\
        \midrule
        SNDP-TDP & \cite{wang2020modelling} & S/M & DE  &\checkmark  &  &\checkmark &  &  &  & \checkmark& \checkmark & & &\checkmark  & \checkmark &  & SIM & DES/Witness \\
        \midrule
        SNDP-TDP & \cite{laurent2020carbonroadmap} & M & DE  & \checkmark &  & &  &  &  & & \checkmark & & & \checkmark &  &  & DSS& Label-setting Algorithm/C++, Web-based (C\#, Javascript, HTML) \\
        \midrule
        SNDP-TDP & \cite{maneengam2020bi} & M & DE  &  \checkmark & & &  &  &  & &  & & & \checkmark & \checkmark &  & E & Branch and bound/Gurobi \\
        \midrule
        SCND & \cite{wang2020integrated} & S & ST  & \checkmark &  & & \checkmark & \checkmark &  & \checkmark &  & & &  &  &  &  MH & Improved GA/MATLAB\\
        \midrule
        SNDP-SSP & \cite{heinold2020emission} & M & DE  &  \checkmark&  & \checkmark&  &  &  & &  & & & \checkmark & \checkmark &  & E & Branch and bound/CPLEX \\
        \midrule
        PS & \cite{jiang2020regional} & B & ST  &\checkmark&  &\checkmark &  &  &  &\checkmark & \checkmark & \checkmark& &  &  &  & E & Branch and bound/CPLEX \\
        \midrule
        SNDP-TDP & \cite{dai2020distributionally} & M & ST  &  \checkmark&  & \checkmark&  & \checkmark &  & &  & & &\checkmark  & \checkmark &  & E & Branch and bound/CPLEX \\
        \midrule
        SCND & \cite{mousavi2020robust} & S & ST  &  \checkmark&  &\checkmark &  &  \checkmark& \checkmark & \checkmark&  & & &  &  &  & E& Branch and bound/GAMS+CPLEX \\
        \midrule
        SNDP-SSP & \cite{guo2020dynamic} & S & DE  &  \checkmark& \checkmark & \checkmark& \checkmark &  \checkmark&  & &  & & &  &  &  & E & Branch and bound/CPLEX \\
        \midrule
        SNDP-TDP & \cite{yin2021interrelations} & M & DE  & \checkmark &  & \checkmark&  &  &  & &\checkmark  &\checkmark & &\checkmark  &\checkmark  &  & E & Branch and bound/MATLAB \\
        \midrule
        PS & \cite{yang2021coastal} & S & DE  &  &  & &  &  &  & &  & & & \checkmark &  &  & HE & Active Set Algorithm \\
        \midrule
        PS & \cite{qian2021decision} & S & DE  &  &  & &  &  &  & \checkmark &  & & & \checkmark &  &  & DSS & DEA-SBM \\
        \midrule
        SNDP-AMP & \cite{tiwari2021freight} & S & DE  & \checkmark &  & &  &  &  & &  & & &  &  &  & E & Branch and bound/LINGO \\
        \midrule
        SNDP-SSP & \cite{wu2021research} & S & DE  & \checkmark & \checkmark & \checkmark&  &  &  & \checkmark&  & & &  &  &  & MH & PSO/MATLAB \\
        \midrule
        NDP+SNDP & \cite{gallardo2021sequential} & S & DE  &  &  & &  &  &  & &\checkmark  & & & \checkmark &  &  & HYB & Dijkstra’s Algorithm/ ArcGIS+ Python 2, DES/Anylogic \\
        \midrule
        SNDP-RMP & \cite{wang2021bi} & B & DE  &  \checkmark& \checkmark & &  &  &  & &  & & & &  &  &  HE & Sensitivity Analysis Algorithm \\
        \midrule
        SNDP-SSP & \cite{zhang2021low} & S & ST  & \checkmark & \checkmark & &  &  &  & \checkmark&  & & &  &  &  & HYB & CA-GA + MC Sampling \\
        \midrule
        
        \\

        
    \multirow{2}{4em}{Classification} & \multirow{2}{12em}{Publications} & \multirow{2}{3em}{S/M/B} & \multirow{2}{4em}{DE/ST} & \multicolumn{13}{c}{Objective Function}  & \multirow{2}{6em}{Solution Technique } & \multirow{2}{12em}{Specific Technique and Tool}\\  
    \cline{5-17} \\ 
    & & & & TC & EC & HC & PC & IC & FC & O & TT & HT & WT & TE & HE & FE & & \\        
        \midrule
        SNDP-SSP & \cite{liang2021multi} & M & DE  & \checkmark & \checkmark & \checkmark&\checkmark  &\checkmark  &  & &  & & &  &  &  & MH & NSGA-III \\
        \midrule
        SNDP-SSP & \cite{hruvsovsky2021real} & S & DE  & \checkmark &  & &  &  &  & \checkmark &  & & & \checkmark &  &  & HYB & Branch and bound/CPLEX, ABM+DES/ANYLOGIC \\
        \midrule   
        SNDP-TDP & \cite{li2022path} & S & ST  &  \checkmark & \checkmark & \checkmark&  &  &  &\checkmark &  & & &  &  &  & MH & FAGSO \\
        \midrule
        SNDP-SSP & \cite{sun2022green} & M & ST  & \checkmark &  & \checkmark&  & \checkmark &  & &  & & & \checkmark &  &  & E & Branch and bound/LINGO \\
        \midrule
        SNDP-TDP & \cite{xie2022research} & M & DE  &\checkmark  &\checkmark & &  &  & \checkmark & &\checkmark  & \checkmark & &  &  &  &  MH & NSGA-II \\
        \midrule
        PS & \cite{pedinotti2022freight} & S & DE &  &  & &  &  &  & \checkmark&  & & & \checkmark &  &  & DSS & NATEM \\
        \midrule   
        SNDP-AMP & \cite{ko2022stochastic} & M & DE/ST  & \checkmark & \checkmark & \checkmark &  &  &  & &  & & & \checkmark &  &  & E  & Random Vectorization of Cost Parameters \\    
        \midrule
        SNDP-SSP & \cite{zhang2022collaborative} & M & DE  & \checkmark & & \checkmark& \checkmark & \checkmark &  & &  & & & \checkmark &  &  & MH  & ALNS \\
        \midrule
        SNDP-SSP & \cite{qi2022transport} & S & DE  &  \checkmark &  \checkmark& &  & \checkmark &  & &  & & &  &  &  & E  & Branch and Bound/CPLEX \\
        \midrule
        SNDP-SSP & \cite{zhang2022preference} & M & DE  & \checkmark &  & \checkmark&&  &  &\checkmark &  \checkmark& & &  \checkmark&  &  &  MH & ALNS \\
        \midrule
        SNDP-SSP & \cite{zhang2022flexible} & S & DE  & \checkmark & \checkmark & \checkmark& \checkmark & \checkmark &  & \checkmark&  & & &  &  &  &  MH & ALNS \\
        \midrule
        SNDP-SSP & \cite{zhang2022heterogeneous} & M & DE  & \checkmark & \checkmark & \checkmark& \checkmark & \checkmark &  & \checkmark&  & & &  &  &  & MH  & ALNS \\
        \midrule
        PS & \cite{ke2023optimization} & M & DE  & \checkmark &  & &  &  &  & &  & & & \checkmark &  &  &  MH & Adaptive GA/MATLAB \\
        \midrule
        SNDP-SSP & \cite{heinold2023primal} & S & ST  & \checkmark &  & & \checkmark &  &  & &  & & &  &  &  & O & Value Function Approximation \\
        \midrule
        SNDP-SSP & \cite{li2023optimum} & M & ST  & \checkmark &  &\checkmark & \checkmark & \checkmark &  & &\checkmark  & \checkmark& & \checkmark &  &  & E & Branch and bound/MATLAB \\
        \midrule
        SNDP-SSP & \cite{sun2023modeling} & S/M & ST  & \checkmark & \checkmark & \checkmark&  & \checkmark &  & &  & & &  &  &  & E& Branch and bound/LINGO \\
        \midrule
        PS & \cite{chen2023subsidy} & B & DE  &  &  & &  &  &  & &  & & & \checkmark &  &  &  GT & Nash Equilibrium, Continuous Optimization Algorithm/ MATLAB \\
        \midrule
        PS & \cite{wu2023evaluation} & S & DE  & \checkmark &\checkmark  & & \checkmark &  &  & \checkmark&  & & &  &  &  & O & Augmented Lagrangian Multiplier Method \\
        \midrule
        NDP+SNDP & \cite{sun2023study} & B & DE  & \checkmark & \checkmark & &  &  & \checkmark & \checkmark&  & & &  &  &  & E & Frank-Wolfe Method/MATLAB \\
        \midrule
        SNDP-SSP & \cite{shoukat2023sustainable} & M & DE  & \checkmark &  & &  &  &  & &  & & & \checkmark &  &  &  MH & GA/MATLAB \\
        \midrule
        SNDP-TDP & \cite{feng2023multimodal} & S & DE  &\checkmark  &  &\checkmark &  & \checkmark &  & \checkmark&  & & &  &  &  & E & Branch and bound/CPLEX \\
        \midrule
        SNDP-TDP & \cite{liu2023multimodal} & S & DE  & \checkmark & \checkmark & &  & \checkmark &  & \checkmark&  & & &  &  &  &  MH & Hummingbird Evolutionary GA/ MATLAB \\
        \midrule
        PS & \cite{nassar2023system} & M & DE  & \checkmark&  & &  &  & \checkmark &\checkmark &  & & & \checkmark &  &  & SIM & SD \\
        \midrule      
        \\

    \multirow{2}{4em}{Classification} & \multirow{2}{12em}{Publications} & \multirow{2}{3em}{S/M/B} & \multirow{2}{4em}{DE/ST} & \multicolumn{13}{c}{Objective Function}  & \multirow{2}{6em}{Solution Technique } & \multirow{2}{12em}{Specific Technique and Tool}\\  
    \cline{5-17} \\ 
    & & & & TC & EC & HC & PC & IC & FC & O & TT & HT & WT & TE & HE & FE & & \\
        \midrule                       
        SNDP-TDP & \cite{zhang2023research} & M & DE  & \checkmark &  & \checkmark&\checkmark  &  &  & \checkmark&\checkmark  &\checkmark &\checkmark &  &  &  & MH & Modified GA/MATLAB \\
        \midrule
        
        SNDP-SSP & \cite{oudani2023combined} & M & DE  &\checkmark  &  &\checkmark &  &  &  & & \checkmark & \checkmark& & \checkmark &  &  & DSS & MARCOS, VIKOR, and Possibility Degree MCDM Approaches  \\
        \midrule        
        PS & \cite{shams2023game} & S & DE  & \checkmark & \checkmark & &  &  &  & \checkmark &  & & & \checkmark &  &  & GT& Stackelberg Game, Nash Game \\
        \midrule
        SNDP-RMP & \cite{rahiminia2023adopting} & B & DE  & \checkmark &  & &  &  &  & \checkmark& \checkmark & & &  \checkmark&  &  &GT & Karush–Kuhn–Tucker Approach \\
        \midrule
        NDP-LP & \cite{li2023hierarchical} & M & DE  & \checkmark & \checkmark &\checkmark &  &  & \checkmark & & \checkmark & \checkmark& &  &  &  & MH & Adaptive GA \\
        \midrule
        PS & \cite{halim2023boosting} & M & DE  & \checkmark &  & \checkmark &  &  &  & & \checkmark & \checkmark & & \checkmark &  &  & DSS & Gravity Model, Discrete Choice Model, Route Choice Model, ASIF  \\
        \midrule
        SNDP-TDP & \cite{yang_multi-objective_2023} & M & DE  & \checkmark &  & &  &  &  & & \checkmark & & & \checkmark &  &  &  MH & Improved Fuzzy Adaptive GA/MATLAB \\
        \midrule
        SNDP-SSP & \cite{de2023flexible} & S & DE  &  & \checkmark & \checkmark&  & \checkmark &  & \checkmark&  & & &  &  &  & SIM & DES \\
        \midrule
        SNDP-AMP & \cite{guo2023integrated} & S & DE  &  &  & &  &  &  & \checkmark 
            &  & & &  &  &  & MH & PSO/CPLEX \\
        \midrule
        SNDP-SSP & \cite{kurtulucs2023green} & M & DE  & \checkmark 
        & \checkmark & \checkmark&  &  &  & \checkmark  & 
        & & & &  &  & E & Branch and bound/Gurobi \\
        \midrule
        PS & \cite{zhong2023system} & S & DE  &  & \checkmark & &  &  &  & &  & & &  &  &  & SIM & SD/Vensim PLE \\
        \midrule
        PS & \cite{guo2023toward} & M & ST  & \checkmark & \checkmark &  &  &  & \checkmark & \checkmark  & \checkmark&\checkmark & &  &  &  & SIM & SD, Monte Carlo \\
        \midrule
        NDP-LP & \cite{wu2023game} & B &  DE & \checkmark & \checkmark & &  &  &  & &  & & &  &  &  & GT & Complete Enumeration Algorithm, GA \\
        \midrule
        SNDP-AMP & \cite{castrellon2023assessing} & M &  DE & \checkmark &  & &  &  &  & &  & & & \checkmark &  &  & SIM & ABM+DES \\
        \bottomrule
    \begin{threeparttable}
\blankfootnote{
\footnotesize\textbf{NA:} No Applicable, \textbf{S:} Single-objective Optimization, \textbf{M:} Multi-objective Optimization, \textbf{B:} Bi-level Optimization, \textbf{DE:} Deterministic, \textbf{ST:} Stochastic, \textbf{TC:} Transportation Cost, \textbf{E:} Exact, \textbf{MH:} Meta-heuristic, \textbf{HE:} Heuristic, \textbf{SIM:} Simulation, \textbf{HYB:} Hybrid, \textbf{DSS:} Decision Support System, \textbf{O:} Others, \textbf{DES}: Discrete Event Simulation, \textbf{CS:} Commercial Solver, \textbf{CC:} Chance Constrained Programming,  \textbf{AHP:} Analytic Hierarchy Process, \textbf{DEA:} Data Envelopment Analysis, \textbf{MCDM:} Multi-criteria decision making, \textbf{GP:} Goal Programming, \textbf{HEDA:} Hybrid Estimation of Distribution Algorithm, \textbf{LS:} Local Search, \textbf{GFWHA:} Genetic and Frank–Wolfe Hybrid Algorithm, \textbf{PSO}: Particle Swarm Optimization, \textbf{DE}: Differential Evolution, \textbf{CW}: Clarke–Wright Savings Algorithm, \textbf{MC:} Monte Carlo, \textbf{NSGA:} Non-dominated Sorting GA, \textbf{SBM:} Slacks-Based Measure, \textbf{CA-GA:} Catastrophic Adaptive GA, \textbf{ASIF:} Activity - Structure - Intensity - Factor,  \textbf{FAGSO:} Fireworks Algorithm with Gravitational Search Operator, \textbf{NATEM:} North American TIMES Energy Model, \textbf{ALNS:} Adaptive Large Neighborhood Search,  \textbf{SD:} System Dynamics
} 
     \end{threeparttable}
\end{longtable}

\end{landscape}

\newpage

\end{document}